\documentclass[reqno,11pt]{amsart}
\usepackage{amsfonts}
\usepackage{graphicx}
\usepackage{amsfonts,amsmath, amssymb}
\usepackage{marginnote}
\usepackage{color}
\usepackage{xcolor}

\numberwithin{equation}{section}

\newcommand{\dis}{\displaystyle}

\newcommand{\R}{\mathbb{R}}

\newtheorem{theorem}{Theorem}[section]

\newtheorem{lemma}[theorem]{Lemma}
\newtheorem{proposition}[theorem]{Proposition}
\newtheorem{remark}[theorem]{Remark}
\newtheorem{definition}[theorem]{Definition}

\def\v{\varepsilon}

\def\t{\theta}

\def\g{\gamma}
\def\d{\delta}
\def\l{\lambda}

\def\r{\rho}

\def\f{\frac}
\def\dd{{\rm d}}
\def\M{{\mathcal{M}}}

\usepackage{tikz}

\begin{document}

\title[Compressible Euler-Poisson Equations with Doping Profile]{Global Solutions of the Compressible Euler-Poisson Equations
for Plasma with Doping Profile for Large Initial Data of Spherical Symmetry}

\author[G.-Q. Chen]{Gui-Qiang G. Chen}
\address{Gui-Qiang G.Chen:\, Mathematical Institute, University of Oxford, Oxford OX2 6GG, UK;
School of Mathematical Sciences, Fudan University, Shanghai 200433, China;
AMSS, Chinese Academy of Sciences, Beijing 100190, China}
\email{chengq@maths.ox.ac.uk}

\author[L. He]{Lin He}
\address{Lin He:\, College of Mathematics, Sichuan University, Chengdu 610064, China;
Academy of Mathematics and Systems Sciences, Chinese Academy of Sciences, Beijing 100190, China}
\email{lin\_he@scu.edu.cn}
\author[Y. Wang]{Yong Wang}
\address{Yong Wang:\, Academy of Mathematics and Systems Sciences, Chinese Academy of Sciences, Beijing 100190, China;
School of Mathematical Sciences, University of Chinese Academy of Sciences, Beijing 100049, China}
\email{yongwang@amss.ac.cn}

\author[D.~F. Yuan]{Difan Yuan}
\address{Difan Yuan: \,School of Mathematical Sciences,  Beijing Normal University and Laboratory of Mathematics and Complex Systems,
Ministry of Education, Beijing 100875, China; Mathematical Institute, University of Oxford, Oxford OX2 6GG, UK}
\email{yuandf@amss.ac.cn}

\begin{abstract}
We establish the global-in-time existence of
solutions of finite relative-energy for the multidimensional compressible
Euler-Poisson equations for plasma with doping profile for large initial data of spherical symmetry.
Both the total initial energy and the initial mass are allowed to
be {\it unbounded},
and the doping profile is allowed to be of large variation.
This is achieved by adapting a class of degenerate density-dependent viscosity terms,
so that a rigorous proof of the inviscid limit of global weak solutions of
the Navier-Stokes-Poisson equations with the density-dependent viscosity terms to
the corresponding global solutions of the Euler-Poisson equations for plasma with doping profile
can be established.
New difficulties arise when tackling the non-zero varied doping profile,
which have been overcome by establishing some novel estimates for the electric field terms
so that the neutrality assumption on the initial data is avoided.
In particular, we prove that no concentration is formed in the inviscid limit
for the finite relative-energy solutions of the compressible Euler-Poisson equations
with large doping profiles in plasma physics.
\end{abstract}

\keywords{Euler-Poisson equations, Navier-Stokes-Poisson equations, multidimension, plasma, doping profile, compressible,  compactness framework, inviscid limit, large data,
relative-energy, infinite mass, spherical symmetry,  {\it a priori} estimate, higher integrability, approximate solutions}
\date{\today}
\maketitle


\thispagestyle{empty}
\section{Introduction}
The Euler-Poisson equations have been intensively studied due to their vast relevance
to modeling physical phenomena, especially in semiconductor modeling and plasma physics;
see \cite{Dafermos,Guo2011,Jackson,Makino1990,Markowich}
and the references cited therein.
In this paper, we are concerned with the global-in-time existence of solutions for
the multidimensional (M-D) compressible Euler-Poisson equations (CEPEs)
for plasma with doping profile for large initial data of spherical symmetry in $\R^N:$
\begin{align}\label{1.1}
\begin{cases}
\partial_t \rho+\mbox{div} \mathcal{M}=0,\\
\partial_t \mathcal{M}+\mbox{div}(\frac{\mathcal{M}\otimes \mathcal{M}}{\rho})+\nabla p+\rho\nabla\Phi=\mathbf{0},\\
\Delta\Phi =d(\mathbf{x})-\rho,
\end{cases}
\end{align}
for $(t,\mathbf{x})\in\R_+\times\mathbb{\R}^N$ with $N\geq3,$
where $\r$ is the density, $p$ is the pressure, and $\mathcal{M}\in\mathbb{\R}^N$ represents the momentum,
$d(\mathbf{x})$ is the doping profile with
$\lim\limits_{|\mathbf{x}|\rightarrow \infty}d(\mathbf{x})=\rho_{\ast}>0$
so that the background state varies with the space variable, and
$\Phi$ is the self-consistent electric field potential function.
When $\rho>0,$ $U=\frac{\mathcal{M}}{\rho}\in \R^N$ is the velocity.
The constitutive pressure-density relation for polytropic perfect gases is
\begin{equation}\nonumber
p=p(\r)=\kappa\r^{\gamma},
\end{equation}
where $\gamma>1$ is the adiabatic exponent and $\kappa=\frac{(\gamma-1)^2}{4\gamma}$ (without loss of generality).

We are concerned with the Cauchy problem \eqref{1.1} with the Cauchy initial data:
\begin{align}\label{initial}
(\rho,\mathcal{M})(0,\mathbf{x})=(\rho_0,\mathcal{M}_0)(\mathbf{x})\longrightarrow (\rho_{\ast},\mathbf{0})
\qquad\, \text{as $|\mathbf{x}|\rightarrow\infty$},
\end{align}
subject to the asymptotic condition:
\begin{align}\label{asym}
\Phi(t,\mathbf{x})\longrightarrow0\qquad\, \text{as $|\mathbf{x}|\rightarrow\infty$}.
\end{align}
Here $(\rho_{\ast},\mathbf{0})$ is the constant far-field state for $(\rho, \mathcal{M})$,
for which the far-field velocity has been transformed
to be zero by the Galilean invariance of CEPEs \eqref{1.1}, without loss of generality.
Since a global solution of CEPEs \eqref{1.1} usually contains the vacuum states $\{(t,\mathbf{x}):\rho(t,\mathbf{x})=0\}$ where the fluid
velocity $U(t,\mathbf{x})$ is not well-defined, we use the physical variables
such as momentum $\mathcal{M}(t,\mathbf{x})$
or $\frac{\mathcal{M}(t,\mathbf{x})}{\sqrt{\rho(t,\mathbf{x})}}$, instead of $U(t,\mathbf{x})$,
when the vacuum states occur.
These choices of variables will be shown to be globally well-defined.
It is challenging to prove the existence of global solutions of problem \eqref{1.1}--\eqref{asym},
due to the possible formation of shock waves and blowups
in a finite time \cite{Chen1998,Deng,Guo1998,Makino1992,Makino1990}.
To solve this problem, we consider the inviscid limit of the solutions of
the compressible Navier-Stokes-Poisson equations (CNSPEs)
with delicately designed density-dependent viscosity terms in $\R^N:$
\begin{align}\label{nsp}
\begin{cases}
\partial_t \rho+\mbox{div} \mathcal{M}=0,\\
\partial_t \mathcal{M}+\mbox{div}\big(\frac{\mathcal{M}\otimes \mathcal{M}}{\rho}\big)+\nabla p+\rho\nabla\Phi
=\varepsilon \mbox{div}\big(\mu(\rho) D(\frac{\mathcal{M}}{\rho})\big)+\varepsilon\nabla\big(\lambda(\rho)\mbox{div}(\frac{\mathcal{M}}{\rho})\big),\\
\Delta\Phi =d(\mathbf{x})-\rho,
\end{cases}
\end{align}
where $D(\frac{\mathcal{M}}{\rho})=\frac{1}{2}\big(\nabla(\frac{\mathcal{M}}{\rho})+(\nabla(\frac{\mathcal{M}}{\rho}))^\top\big)$ is the stress tensor,
$\varepsilon>0$ is the inverse of the Reynolds number,
and the Lam$\text{\'{e}}$ (shear and bulk) viscosity coefficients $\mu(\rho)$ and $\lambda(\rho)$ depend on the density, which may vanish on the vacuum
and satisfy
$$
\mu(\rho)\geq0,\quad\mu(\rho)+N\lambda(\rho)\geq0 \qquad\,\,\, \text{for $\rho\geq0$}.
$$
Formally, for the inviscid limit $\varepsilon\rightarrow0+,$ CNSPEs \eqref{nsp} converges to CEPEs \eqref{1.1}.
However, the rigourous mathematical proof of this limit has been one of the most challenging problems
in compressible fluid dynamics;
see  Chen-Feldman \cite{Chen2018}, Dafermos \cite{Dafermos}, and the references cited therein.

The system of CEPEs \eqref{1.1} is a fundamental prototype of hyperbolic balance laws arising in either the two-fluid theory
in plasma physics or in the theory of self-consistent gravitational gaseous star.
In the classical two-fluid model describing plasma dynamics for electron-ion fluids, the self-consistent electrostatic Newtonian potential satisfies
the Poisson equation for
the repulsive Coulomb interaction.
A plasma consists of rapidly-moving charged particles.
Various evidence shows that more than ninety percents of the matter in the visible universe
is of the form of plasma, for instance,  the interior of stars, sparse intergalactic plasma,  nebulas, neon signs, etc.
The existence/nonexistence of weak solutions and the stability/instability
of the CEPEs for plasma are among the key challenging
problems for nuclear fusion. Specifically, charged particles are accelerated at high speed by particle accelerators to emit the electromagnetic energy.

For CEPEs \eqref{1.1} with constant doping profile,
Guo \cite{Guo19981} first constructed a global smooth irrotational solution without shock
by using the dispersive Klein-Gordon effect and adapting Shatah's normal form method.
It is due to the enhanced dispersive effects induced by the repulsive electrostatic interaction, which is different from the pure compressible Euler equations for neutral gas.
Guo-Strauss \cite{Guo2005} proved the existence of global smooth solutions
near a given steady state of the Euler-Poisson equations, for which
the steady state and the doping profile are permitted to be of large variation
while the initial velocity must be small.
Later on, Guo-Pausader \cite{Guo2011} constructed global smooth irrotational solutions with small amplitude for
the Euler-Poisson equations for ion dynamics.
Jang \cite{Jang2012} constructed two-dimensional (2-D) global smooth solutions
for spherically symmetric flows with small perturbed initial data.
A family of smooth solutions was constructed by \cite{Jang2014}.
Smooth irrotational solutions in $\R^2$ were constructed independently
by Ionescu-Pausader \cite{Ionescu} and Li-Wu \cite{Li2014}.
Such a surprising and subtle dispersive property has also been identified
and exploited in other two-fluid models,
which leads to the persistence of global smooth solutions and absence of shock formation.
Guo-Ionescu-Pausader \cite{Guo2016} first proved the global stability of a constant neutral background
in the sense that irrotational, smooth, and localized perturbations of a constant background
with small amplitude lead to global smooth solutions in $\R^3$
for the Euler-Maxwell two-fluid systems.
These results can be applied equally well to other plasma models such as
the Euler-Poisson systems for two-fluids.

As indicated earlier, due to the hyperbolic structure of the Euler-Poisson equations, for large initial data,
the smooth solutions of \eqref{1.1} may form shock waves and blowup in finite time; see \cite{Chen1998,Deng,Guo1998,Makino1992,Makino1990}.
For the compressible Euler equations, we refer
to \cite{Chen1988,Chen1990,Chen2010,Chen2018b,Chen2020,Coutand2012,Ding1989,Diperna1983,Huang2002,Jang2015,LeFloch,Lions1996,Lions1994}
and the references therein.
For the global existence of solutions of the compressible Navier-Stokes equations,
we refer to \cite{Feireisl,Jiang,Lions1998} for the case of constant viscosity
and \cite{Li2015,Vasseur}
for the case of density-dependent viscosity and the reference therein.
In \cite{Li2015,Vasseur}, the BD entropy plays essential roles to obtain the derivative estimate of the density,
which was first identified
by Bresch-Desjardins \cite{Bresch2002}  (also see \cite{Bresch2003} and \cite{Bresch20072}).
The global existence of weak solutions of the compressible Navier-Stokes equations
with spherical symmetry is established in
\cite{Guo-Jiu-Xin-2}
for $\gamma\in(1,3)$ in a finite region under Dirichlet boundary conditions.

Regarding to the inviscid gases as viscous gases with vanishing real physical viscosity
can date back to the pioneer paper by Stokes \cite{Stokes} and the important contributions
of Rankine \cite{Rankine}, Hugoniot \cite{Hugoniot},
and Rayleigh \cite{Rayleigh} ({\it cf.} Dafermos \cite{Dafermos}).
However, the first convergence analysis of the inviscid limit
from the one-dimensional (1-D) barotropic Navier-Stokes to Euler equations
was  rigorously proved by Gilbarg \cite{D. Gilbarg}, in which he established the mathematical existence and vanishing viscous limit of the Navier-Stokes shock layers.
For the convergence analysis confined in the framework of piecewise smooth solutions,
see \cite{GMWZ,Hoff1989,Xin}
and the references cited therein. In Chen-Perepelitsa \cite{Chen2010},
the vanishing real physical viscosity limit of smooth solutions for the 1-D
Navier-Stokes equations to the corresponding finite relative-energy solution
of the Euler equations was first established for $\rho_{\ast}>0$.
In \cite{Chen2015,Chen2018b}, the existence theory of global finite-energy entropy solution of
the Euler equations
with spherical symmetry and large initial data
was established for $\rho_{\ast}=0$ through constructing artificial viscosity approximate smooth solutions (also see \cite{Schrecker}).
Chen-Wang \cite{Chen2020} proved the existence of global finite relative-energy solutions for the
compressible Euler equations with spherical symmetry
through taking a vanishing physical viscosity limit even for $\rho_*>0$.
Most recently, Chen-He-Wang-Yuan \cite{Chen2021} established the global existence of finite-energy solutions of the M-D Euler-Poisson equations
for both compressible gaseous stars and plasmas with large initial data of spherical symmetry.
For further work related with CEPEs with doping profile,
we refer to \cite{Bae,Li2002,HLiu2020,Marcati1995} and the references cited therein.

In this paper, we are concerned with the spherically symmetric solutions of \eqref{1.1}
for plasma with doping profile.
Denote
\begin{align}\label{1.2}
\quad\r(t,\mathbf{x})=\r(t,r),\quad \mathcal{M}(t,\mathbf{x})=m(t,r)\f{\mathbf{x}}{r},
\quad \Phi(t,\mathbf{x})=\Phi(t,r)
\qquad\,\, \text{for $r=|\mathbf{x}|$},
\end{align}
subject to the initial conditions:
\begin{align}\label{initial2}
(\rho,\mathcal{M})(0,\mathbf{x})=(\rho_0,\mathcal{M}_0)(\mathbf{x})=(\rho_0(r),m_0(r)\frac{\mathbf{x}}{r})
\longrightarrow (\rho_{\ast},\mathbf{0}) \qquad\, \text{as $|\mathbf{x}|\longrightarrow\infty$},
\end{align}
and the asymptotic condition:
\begin{align}\label{asym2}
\Phi(t,\mathbf{x})=\Phi(t,r)\longrightarrow0 \qquad\, \text{as $|\mathbf{x}|\rightarrow\infty$}.
\end{align}
It is not necessary to impose the initial data for $\Phi$,
since $\Phi(0,\mathbf{x})$ can be determined by the initial density and the boundary condition \eqref{asym2}.

We establish the inviscid limit of the corresponding spherically symmetric solutions of CNSPEs \eqref{nsp}
with the adapted class of degenerate density-dependent viscosity terms and approximate initial data of form \eqref{initial2}.
For spherically symmetric solutions of form \eqref{1.2}, systems \eqref{1.1} and \eqref{nsp} become
\begin{align}\label{euler}
\begin{cases}
\r_t+m_r+\frac{N-1}{r}m=0,\\[1mm]
m_t+(\frac{m^2}{\rho}+p)_r+\frac{N-1}{r}\frac{m^2}{\rho}+\rho\Phi_r=0,\\[1mm]
\Phi_{rr}+\frac{N-1}{r}\Phi_r=d(r)-\rho,
\end{cases}
\end{align}
and
 \begin{align}\label{1.3}
\begin{cases}
\r_t+m_r+\frac{N-1}{r} m=0,\\[1mm]
m_t+(\frac{m^2}{\rho}+p)_r+\frac{N-1}{r} \frac{m^2}{\rho}+\rho\Phi_r
=\v\Big((\mu+\lambda)\big((\frac{m}{\rho})_r+\frac{N-1}{r}\frac{m}{\rho}\big)\Big)_r
-\v\frac{N-1}{r}\mu_r\frac{m}{\rho},\\[1mm]
\Phi_{rr}+\frac{N-1}{r}\Phi_r=d(r)-\rho,
\end{cases}
\end{align}
respectively.

For the potential function $\Phi,$ we need to assume the boundary condition:
\begin{equation}\nonumber
\lim_{r\rightarrow0+}r^{N-1}\Phi_r(t,r)=0\qquad\,\mbox{for all $t\geq0$}.
\end{equation}
This leads to
\begin{equation}\nonumber
r^{N-1}\Phi_r(t,r)=-\int^{r}_0\big(\rho(t,y)-d(y)\big)\,y^{N-1}{\rm d}y.
\end{equation}

The study of spherically symmetric solutions can date back to the 1950s and
has been motivated by many important physical problems such as the stellar dynamics including gaseous stars
and supernova formation \cite{Courant1948,Merle,Merle2}.
A fundamental open problem is whether the concentration is formed at the origin for
both the compressible Euler equations and the Euler-Poisson equations,
which has been answered for the solutions as the inviscid limits in \cite{Chen2021,Chen2020} (without doping profile).
In this paper, we prove that no concentration (no delta measure) is formed at the origin, for CEPEs \eqref{1.1} even with doping profile for plasma.

The main difficulties in this paper are twofolds:

\smallskip
(i) Since the initial data satisfy \eqref{initial}, the initial mass is allowed to be {\it unbounded}:
$$
\int_{\R^N}\rho_0(\mathbf{x})\,\dd  \mathbf{x}=\omega_N\int^{\infty}_0\rho_0(r)\,r^{N-1}\dd  r=\infty,
$$
where $\omega_N:=\frac{2\pi^{\frac{N}{2}}}{\Gamma(\frac{N}{2})}$ denotes the surface area of unit sphere in $\R^N$.
In addition, the initial total energy is ${\it unbounded}$, while initial total
relative-energy is finite.  We overecome this difficulty by fully taking advantage of \eqref{mc} via
the conservation of relative mass. Note that the finite ``relative-mass'' depends on $b.$

\smallskip
(ii) Since the compressible Euler-Poisson equations involve the non-zero varied doping profile,
new difficulties occur due to the electric field potential terms $\Phi$ depending on the doping profile.
To solve these new difficulties, we make a more careful analysis in the estimates of
the electric field potential terms.
In order to control these terms, our key observation is to use the boundedness of the initial total relative-energy
to control the electric potential,
especially to tackle the non-zero doping profile terms.
It is worth mentioning that, differently from \cite{Chen2020},
when proving the lower bound of the density,
we require the higher exponent of integrability for the derivatives of the density;
see Lemma \ref{lem6.7}, especially \eqref{7.79}.
When we prove the higher integrability of the velocity, we need the extra weight for
the electric field terms; see \eqref{M1}--\eqref{4.1d}.
These estimates indicate that the estimates of the electric field should be much more involved
in the proof of Lemma \ref{lem6.10} and the follow-up lemmas than those for the case without doping profile in
\cite{Chen2021}.

We remark that no imperative neutrality condition is assumed on the initial data:
$$
\int_{\R^N}(\rho_0(\mathbf{x})-d(\mathbf{x}))\,\dd  \mathbf{x}\neq0
$$
in our analysis. It is natural to avoid this nonphysical assumption when we take into account of the effect of the positive density
at the far fields $(|\mathbf{x}|\rightarrow\infty)$.
This is achieved by constructing approximate solutions without
imposing the extra boundary restriction: $\Phi_r(t,b)=0$, where $b$ is the upper bound of
$|\mathbf{x}|$ for the truncated domain.
To avoid this imperative condition, we need to establish some novel estimates
when performing integration by parts in the proof
of BD-type entropy estimates.
The boundary terms are needed to be delicately controlled;
for example, see the last term of the right-hand side of \eqref{6.95}
in Steps 5--6 in the proof
of Lemma \ref{BD}.
The key idea is to make these terms to be controlled by the initial total finite
relative-energy.

The plan of this paper is as follows:
In \S 2, we first introduce the definitions of finite relative-energy solutions of
the Cauchy problem \eqref{1.1}--\eqref{asym}
for CEPEs and then state Main Theorem I: Theorem \ref{existence1} for the global existence of weak solutions.
To prove Theorem \ref{existence1},  the global weak solutions of the Cauchy problem \eqref{nsp} and \eqref{appinitial}--\eqref{phiapproximate}
for CNSPEs are firstly constructed and we analyze their invicid limit,
as stated in Main Theorem II: Theorem \ref{theorem2}.
In \S 3, we give the construction of global approximate smooth solutions
$(\rho^{\v,\delta,b}, m^{\v, \delta,b})$ and
perform the basic energy estimate and the BD-type entropy estimate of $(\rho^{\v,\delta,b}, m^{\v, \delta,b})$, for CNSPEs \eqref{1.3}. In \S 4, we derive the higher integrability of $(\rho^{\v,\delta,b}, m^{\v, \delta,b})$,
uniformly in $b$, for the approximate smooth solutions.
In \S 5, through taking the limit of $(\rho^{\v,\delta,b}, m^{\v, \delta,b})$ as $b\to \infty,$
 we obtain global strong solutions $(\rho^{\v,\delta}, m^{\v, \delta})$ of system \eqref{6.1}
with several uniform bounds in $(\v, \delta)$, and then
we pass to the limit, $\delta\to 0+$, to obtain global existence of spherically symmetric weak solutions of
CNSPEs \eqref{1.3} with several desired uniform bounds
and the $H_{\rm loc}^{-1}$--compactness,
which are important for us to utilize the $L^p$ compensated compactness framework in \S 6
to establish Theorem \ref{theorem2}.
In the appendix,  the approximate initial data with desired estimates are constructed, which
are needed for the construction of the approximate solutions in \S 3.

In this paper, we write $L^p(\Omega),W^{k,p}(\Omega)$, and $H^k(\Omega)$ as the usual Sobolev spaces defined
on $\Omega$ for $p\in[1,\infty]$. We also write $L^p(I; r^{N-1}\dd r)$ or $L^p([0,T)\times I; r^{N-1}\dd r\dd t)$ for the open
interval $I\subseteq \R_+$ with measure $r^{N-1}\dd r$ or $r^{N-1}\dd r\dd t$ respectively,
and $L^p_{\rm loc}([0,\infty); r^{N-1}\dd  r)$ to denote $L^{p}([0,R); r^{N-1}\dd  r)$ for arbitrary fixed $R>0$.

\section{Mathematical Problems and Main Theorems}

In this section, we introduce the definition of weak solutions of finite
relative-energy solutions of the Cauchy problem for CEPEs \eqref{1.1}
in $\R^{N+1}_+:=\R_+\times\R^{N}=(0,\infty)\times \R^{N}$
for $N\geq3$.
We assume that the initial data $(\rho_0,\mathcal{M}_0)(\mathbf{x})$ and the corresponding initial potential function $\Phi_0(\mathbf{x})$ have both finite total initial relative-energy:
\begin{equation}\label{ife}
\mathcal{E}_0:=\int_{\R^N}\Big(\frac{1}{2}\big|\frac{\mathcal{M}_0}{\sqrt{\rho_0}}\big|^2+e(\rho_0,\rho_{\ast})+\frac{1}{2}|\nabla_{\mathbf{x}}\Phi_0|^2\Big)\dd  \mathbf{x}<\infty,
\end{equation}
where $e(\rho,\rho_{\ast})$ is the relative internal energy
with respect to $\rho_{\ast}>0$:
\begin{equation}\label{6.10}
e(\rho,\rho_{\ast})
:=e(\rho)-e(\rho^\ast)-e'(\rho^\ast)(\rho-\rho^\ast)=\frac{\kappa}{\gamma-1} \big(\rho^\g-\r^{\g}_{\ast}-\g \rho_{\ast}^{\g-1} (\rho-\r_{\ast})\big),
\end{equation}
where $e(\rho)=\frac{\kappa}{\gamma-1}\rho^{\gamma}$ denotes the gas internal energy.

We also assume that the doping profile satisfies
\begin{align}
&d(r)\in C([0,\infty)),\quad \lim_{r\rightarrow \infty}d(r)=\rho_{\ast}>0,\label{initialdecay}\\
&\int^{\infty}_0|d(r)-\rho_{\ast}|^jr^{N-1}\dd  r<\infty \qquad\,\,\mbox{for $j=1,2,\frac{\gamma}{\gamma-1}$}.\label{r}
\end{align}

\begin{definition}[Weak Solutions]\label{weakep}
A measurable vector function $(\rho,\mathcal{M},\Phi)$ is said to be a
weak solution of finite relative-energy
of the Cauchy problem \eqref{1.1}--\eqref{asym} provided that
\begin{enumerate}
\item[(i)] $\rho(t,\mathbf{x})\geq0\,$ a.e., and $(\mathcal{M},\frac{\mathcal{M}}{\sqrt{\rho}})(t,\mathbf{x})=\mathbf{0},$ a.e. on the vacuum states $\{(t,\mathbf{x}): \rho(t,\mathbf{x})=0\};$

\smallskip
\item[(ii)] For a.e. $t>0,$ the total relative-energy with respect to the far field state $(\rho_{\ast},\mathbf{0})$ is finite:
\begin{align}\label{finiteenergy}
\int_{\R^N}\Big(\frac{1}{2}\Big|\frac{\mathcal{M}}{\sqrt{\rho}}\Big|^2+e(\rho,\rho_{\ast})+\frac{1}{2}|\nabla_{\mathbf{x}}\Phi|^2\Big)(t,\mathbf{x})\,\dd  \mathbf{x}\leq \mathcal{E}_0.
\end{align}

\item[(iii)] For any $\zeta(t,\mathbf{x})\in C^1_0([0,\infty)\times\R^N),$
\begin{align*}
\int^{\infty}_0\int_{\R^N}\big(\rho\zeta_t+\mathcal{M}\cdot\nabla\zeta\big)\,\dd  \mathbf{x}\dd  t+\int_{\R^N}(\rho\zeta)(0,\mathbf{x})\,\dd  \mathbf{x}=0.
\end{align*}

\item[(iv)] For any
$\mathbf{\psi}(t,\mathbf{x})=(\psi_1,\cdots,\psi_N)(t,\mathbf{x})\in (C^1_0([0,\infty)\times \R^N))^N,$
\begin{align*}
&\int^{\infty}_0\int_{\R^N}
\Big({\mathcal{M}\cdot \mathbf{\partial}_t\psi+\frac{\mathcal{M}}{\sqrt{\rho}}\cdot(\frac{\mathcal{M}}{\sqrt{\rho}}\cdot \nabla)\psi}+p(\rho)\mathrm{div}\psi\Big)\, \dd\mathbf{x}\dd  t
+\int_{\R^N}\mathcal{M}_0(\mathbf{x})\cdot\psi(0,\mathbf{x})\,\dd\mathbf{x}\nonumber\\
&=\int^{\infty}_0\int_{\R^N}\rho\nabla_{\mathbf{x}}\Phi\cdot \psi\,\dd\mathbf{x}\dd  t.
\end{align*}

\item[(v)] For any $\xi(\mathbf{x})\in C^1_0(\R^N),$
\begin{align*}
\int_{\R^N}\nabla_{\mathbf{x}}\Phi(t,\mathbf{x})\cdot\nabla_{\mathbf{x}}\xi(\mathbf{x})\,\dd\mathbf{x}
=\int_{\R^N}\big(\rho(t,\mathbf{x})-d(\mathbf{x})\big)\xi(\mathbf{x})\,\dd \mathbf{x}\qquad \text{ for a.e. $t\geq0$}.
\end{align*}
\end{enumerate}
\end{definition}

\smallskip
\begin{theorem}[Main Theorem I: Existence of Global Solutions of CEPEs for Plasma with Doping Profile and Large Initial Data of Spherical Symmetry]\label{existence1}
Consider the Cauchy problem \eqref{1.1}--\eqref{asym} for CEPEs with large initial data of spherical symmetry
of form \eqref{initial2}--\eqref{asym2}. Assume that the doping profile and
the initial density
satisfy \eqref{ife}, \eqref{initialdecay}--\eqref{r}, and
$\rho_0-d(\mathbf{x})\in (L^{\frac{2N}{N+2}}\cap L^1)(\R^N)$.
Then there exists a globally defined weak solution
$(\rho,\mathcal{M},\Phi)$ of finite relative-energy of
the Cauchy problem \eqref{1.1}--\eqref{asym} and \eqref{initial2}--\eqref{asym2}
with spherical symmetry of form \eqref{1.2} in the sense of {\rm Definition \ref{weakep}}
such that $(\rho,m,\Phi)(t,r)$ is determined by the corresponding system \eqref{euler}
with initial data $(\rho_0,m_0,\Phi_0)(r)$ given in \eqref{initial2}
subject to the asymptotic behavior \eqref{asym2}.
\end{theorem}

To prove Theorem \ref{existence1}, we first construct global weak solutions for CNSPEs \eqref{nsp} with appropriately adapted degenerate density-dependent viscosity terms and approximate initial data:
\begin{equation}\label{appinitial}
(\rho,\mathcal{M},\Phi)|_{t=0}=(\rho^{\varepsilon}_0,\mathcal{M}^{\varepsilon}_0,\Phi^{\varepsilon}_0)(\mathbf{x})\longrightarrow (\rho_0,\mathcal{M}_0,\Phi_0)(\mathbf{x}) \qquad \text{ as $\varepsilon\rightarrow0+$},
\end{equation}
subject to the asymptotic boundary condition:
\begin{equation}\label{phiapproximate}
\Phi^{\varepsilon}(t,\mathbf{x})\rightarrow0\qquad  \text{ as $|\mathbf{x}|\rightarrow\infty$}.
\end{equation}
For simplicity, we consider the viscosity terms with $(\mu,\lambda)=(\rho,0)$ in \eqref{nsp}, and $\varepsilon\in(0,1]$ without loss of generality.

\begin{definition}\label{weaknsp}
A pair $(\rho^{\varepsilon},\mathcal{M}^{\varepsilon},\Phi^{\varepsilon})$ is said to be a weak solution of the Cauchy problem \eqref{nsp} and \eqref{appinitial}--\eqref{phiapproximate} with $(\mu,\lambda)=(\rho,0)$ if the following conditions
hold{\rm :}
\begin{enumerate}
\item[\rm (i)]
$\rho^{\varepsilon}(t,\mathbf{x})\geq0$ a.e.,  $\,(\mathcal{M}^{\varepsilon},\frac{\mathcal{M}^{\varepsilon}}{\sqrt{\rho^{\varepsilon}}})(t,\mathbf{x})=\mathbf{0}$ a.e. on the vacuum states $\{(t,\mathbf{x}): \rho^{\varepsilon}(t,\mathbf{x})=0\},$
\begin{align*}
&\rho^{\varepsilon}\in L^{\infty}(0,T;L^{\gamma}_{\rm loc}(\R^N)),\quad \nabla \sqrt{\rho^{\varepsilon}}\in L^{\infty}(0,T; L^2(\R^{N})),\quad \frac{\mathcal{M}^{\varepsilon}}{\sqrt{\rho^{\varepsilon}}}\in L^{\infty}(0,T; L^2(\R^N)),\\
&\Phi^{\varepsilon}\in L^{\infty}(0,T;L^{\frac{2N}{N-2}}(\R^N)),\quad \nabla\Phi^{\varepsilon}\in L^{\infty}(0,T;L^2(\R^N)).
\end{align*}

\smallskip
\item[\rm (ii)]  For any $t_2\geq t_1\geq0$ and
$\zeta(t,\mathbf{x})\in C^1_0([0,\infty)\times\R^N),$
\begin{align*}
\int_{\R^N}(\rho^{\varepsilon}\zeta)(t_2,\mathbf{x})\,\dd\mathbf{x}-\int_{\R^N}(\rho^{\varepsilon}\zeta)(t_1,\mathbf{x})\,\dd \mathbf{x}=\int^{t_2}_{t_1}\int_{\R^N}(\rho^{\varepsilon}\zeta_t+\mathcal{M}^{\varepsilon}\cdot\nabla\zeta)(t,\mathbf{x})\,\dd  \mathbf{x}\dd  t=0.
\end{align*}

\smallskip
\item[\rm (iii)] For any $\mathbf{\psi}(t,\mathbf{x})=(\psi_1,\cdots,\psi_N)(t,\mathbf{x})\in (C^1_0([0,\infty)\times \R^N))^N,$
\begin{align*}
&\int^{\infty}_0\int_{\R^N}\Big({\mathcal{M}^{\varepsilon}\cdot \partial_t\psi+\frac{\mathcal{M}^{\varepsilon}}{\sqrt{\rho^{\varepsilon}}}\cdot(\frac{\mathcal{M}^{\varepsilon}}{\sqrt{\rho^{\varepsilon}}}\cdot \nabla)\psi}+p(\rho^{\varepsilon})\mathrm{div}\psi\Big)\,\dd   \mathbf{x}\dd   t+\int_{\R^N}\mathcal{M}^{\varepsilon}_0(\mathbf{x})\cdot\psi(0,\mathbf{x})\,\dd \mathbf{x}\nonumber\\
&=-\varepsilon\int^{\infty}_0\int_{\R^N}\Big(\frac{1}{2}\mathcal{M}^{\varepsilon}\cdot (\Delta\psi+\nabla \mathrm{div}\psi)+\frac{\mathcal{M}^{\varepsilon}}{\sqrt{\rho^{\varepsilon}}}\cdot (\nabla\sqrt{\rho^{\varepsilon}}\cdot\nabla)\psi+\nabla\sqrt{\rho^{\varepsilon}}\cdot(\frac{\mathcal{M}^{\varepsilon}}{\sqrt{\rho^{\varepsilon}}}\cdot\nabla)\psi\Big)\,\dd \mathbf{x}\dd t\nonumber\\
&\quad+\int^{\infty}_0\int_{\R^N}(\rho^{\varepsilon}\nabla_{\mathbf{x}}\Phi^{\varepsilon}\cdot \psi)(t,\mathbf{x})\,\dd \mathbf{x}\dd t.
\end{align*}
\item[\rm (iv)] For any $\xi(\mathbf{x})\in C^1_0(\R^N),$
\begin{align*}
\int_{\R^N}\nabla^{\varepsilon}_{\mathbf{x}}\Phi(t,\mathbf{x})\cdot\nabla_{\mathbf{x}}\xi(\mathbf{x})\,\dd \mathbf{x}
=\int_{\R^N}\big(\rho^{\varepsilon}(t,\mathbf{x})-d(\mathbf{x})\big)\xi(\mathbf{x})\,\dd \mathbf{x}\qquad  \text{ for a.e. } t\geq0.
\end{align*}
\end{enumerate}
\end{definition}

We now consider spherical symmetric solutions of form \eqref{1.2}. Then systems \eqref{1.1} and \eqref{nsp} become systems \eqref{euler} and \eqref{1.3}, respectively.
A pair of functions $(\eta(\rho,m),q(\rho,m)),$ or $(\eta(\rho,u),q(\rho,u))$ for $m=\rho u$, is called an entropy-entropy flux pair for the first two equations of
system \eqref{euler} (with $N=1$ and $\Phi=0$) if $(\eta, q)$ is a solution of
$$
\nabla q(\rho, m)=\nabla\eta(\rho, m)\nabla \begin{pmatrix}
m \\
\frac{m^2}{\rho}+p(\rho)\end{pmatrix},
$$
which implies that, for any smooth solution $(\rho(t,r), m(t,r))$ of \eqref{euler},
the following equation holds:
$$
\partial_t\eta(\rho(t,r),m(t,r))+\partial_rq(\rho(t,r),m(t,r))=0.
$$
 Furthermore, $\eta(\rho,m)$ is called a weak entropy if
\begin{equation*}
\eta|_{\rho=0}=0\qquad\text{for any fixed $u=\frac{m}{\rho}$}.
\end{equation*}
We denote $u=\frac{m}{\rho}$ and $m$ alternatively when $\rho>0.$
An entropy $\eta(\r,m)$ is called convex if the Hessian $\nabla^2\eta(\rho,m)$ is nonnegative definite in the region we considered.
From \cite{Lions1994}, it is well-known that any weak entropy pair $(\eta,q)$ can be represented by
\begin{align}\label{weakentropy}
\begin{cases}
\eta^{\psi}(\rho,m)=\eta(\rho,\rho u)=\int_{\R}\chi(\rho;s-u)\psi(s)\,\dd s,\\[1mm]
q^{\psi}(\rho,m)=\eta(\rho,\rho u)=\int_{\R}(\theta s+(1-\theta)u)\chi(\rho;s-u)\psi(s)\,\dd s,
\end{cases}
\end{align}
where the kernel is $\chi(\rho;s-u)=[\rho^{2\theta}-(s-u)^2]_{+}^{\lambda}$ with
$\lambda=\frac{3-\gamma}{2(\gamma-1)}>-\frac{1}{2}$ and $\theta=\frac{\gamma-1}{2}.$
In particular, $\psi(s):=\frac{1}{2}s^2,$ the entropy pair is the mechanical
energy-energy flux pair:
\begin{equation*}
\eta^{*}(\rho,m)=\frac{1}{2}\frac{m^2}{\rho}+e(\rho),\quad q^{*}(\rho,m)=\frac{1}{2}\frac{m^3}{\rho^2}+m e'(\rho).
\end{equation*}

Since $(\rho,m)(t,r)\rightarrow (\rho_{\ast},0)$
with $\rho_{\ast}>0$ as $r\rightarrow\infty$,
We use the relative mechanical energy
\begin{equation*}
\bar{\eta}^{\ast}(\rho,m)=\frac{1}{2}\frac{m^2}{\rho}+e(\rho,\rho_{\ast})
\end{equation*}
with $e(\rho,\rho_{\ast})$ is defined in \eqref{6.10} satisfying
\begin{equation*}
e(\rho,\rho_{\ast})\geq C_{\gamma}\rho(\rho^{\theta}-\rho^{\theta}_{\ast})^2
\end{equation*}
for some constant $C_{\gamma}>0.$

\begin{theorem}[Main Theorem II: Existence and Inviscid Limit for CNSPEs]\label{theorem2}
Consider CNSPEs \eqref{nsp} with $N\geq3$ and the spherically symmetric approximate initial data \eqref{appinitial} satisfying that,
as $\varepsilon\rightarrow0+,$
\begin{equation}\label{ip}
(\rho^{\varepsilon}_0,m^{\varepsilon}_0)(r)\, \rightarrow\,
(\rho_0,m_0)(r) \qquad
\text{in $L^{p}_{\rm loc}([0,\infty); r^{N-1}\dd r)\times L^1_{\rm loc}([0,\infty); r^{N-1}\dd r)$}
\end{equation}
with $p=\max\{\gamma,\frac{2N}{N+2}\}$ and
\begin{align}
&\mathcal{E}^{\varepsilon}_0=\omega_N\int^{\infty}_0\Big(\bar{\eta}^{\ast}(\rho^{\varepsilon}_0,m^{\varepsilon}_0)+\frac{1}{2}|\Phi^{\varepsilon}_{0r}(r)|^2\Big)\,r^{N-1}\dd r\rightarrow \mathcal{E}_0,\label{E0}\\
&\mathcal{E}^{\varepsilon}_1=\varepsilon^2\int^{\infty}_0|(\sqrt{\rho^{\varepsilon}_0})_r|^2\,r^{N-1}\dd r\rightarrow0,\label{E1}
\end{align}
and there exists a constant $C>0$ independent of $\varepsilon\in(0,1]$ such that
\begin{equation}\label{Ebound}
\mathcal{E}^{\varepsilon}_0+\mathcal{E}^{\varepsilon}_1\leq C(\mathcal{E}_0+1)
\end{equation}
for $\mathcal{E}_0$ defined in \eqref{ife}.
Then the following statements hold{\rm :}

\medskip
\noindent
\textbf{Part 1. Existence of Solutions for CNSPEs \eqref{nsp}{\rm :}}
\begin{enumerate}
\item[\rm (i)] For each $\varepsilon>0,$ there exists a
globally-defined spherically symmetric solution
$$(\rho^{\varepsilon},\mathcal{M}^{\varepsilon},\Phi^{\varepsilon})(t,\mathbf{x})=(\rho^{\varepsilon}(t,r),m^{\varepsilon}(t,r)\frac{\mathbf{x}}{r},\Phi^{\varepsilon}(t,r))=(\rho^{\varepsilon}(t,r),\rho^{\varepsilon}(t,r)u^{\varepsilon}(t,r)\frac{\mathbf{x}}{r},\Phi^{\varepsilon}(t,r))$$
of the Cauchy problem of \eqref{nsp} and
\eqref{appinitial}--\eqref{phiapproximate} in the sense of
Definition {\rm \ref{weaknsp}} with
$$
u^{\varepsilon}(t,r)=
\begin{cases}
\frac{m^{\varepsilon}(t,r)}{\rho^{\varepsilon}(t,r)}\quad
& a.e. \mbox{ on $\{(t,r): \rho^{\varepsilon}(t,r)\neq0\}$},\\
0\quad  &a.e.
\mbox{ on $\{(t,r): \rho^{\varepsilon}(t,r)=0\}.$}
\end{cases}
$$
In addition, $(\rho^{\varepsilon},m^{\varepsilon},\Phi^{\varepsilon})(t,r)$ satisfies the following uniform estimates{\rm :}
For $t>0,$
\begin{align}
&\int^{\infty}_0\bar{\eta}^{\ast}(\rho^{\varepsilon},m^{\varepsilon})(t,r)\,r^{N-1}\dd r+\varepsilon\int_{\R^2_+}\rho^{\varepsilon}(s,r)|u^{\varepsilon}(s,r)|^2\,r^{N-3}\dd r\dd s
\nonumber\\
&\qquad
+\frac{1}{2}\int^{\infty}_0|\Phi^{\varepsilon}_{0r}|^2\,r^{N-1}\dd r\leq\frac{\mathcal{E}^{\varepsilon}_0}{\omega_N}\leq C(\mathcal{E}_0+1),
\label{ee}\\[2mm]
&\,\varepsilon^2\int^{\infty}_0|(\sqrt{\rho^{\varepsilon}(t,r)})_r|^2\,r^{N-1}\dd r+\varepsilon\int_{\R^2_+}|((\rho^{\varepsilon}(s,r))^{\frac{\gamma}{2}})_r|^2r^{N-1}\dd r\dd s\leq C(\mathcal{E}_0+1),\label{bd}
\end{align}
and, for any given $T\in(0,\infty)$ and any compactly supported subset $[r_1,r_2]\Subset(0,\infty),$
\begin{align}
&\int^T_0\int^{r_2}_{r_1}(\rho^{\varepsilon}(t,r))^{\gamma+1}\dd r\dd t\leq C(r_1,r_2,T,\mathcal{E}_0),\label{higher}\\
&\int^T_0\int^{r_2}_{r_1}\Big(\rho^{\varepsilon}(t,r)|u^{\varepsilon}(t,r)|^3+(\rho^{\varepsilon}(t,r))^{\gamma+\theta}\Big)\,\dd r\dd t\leq C(r_1,r_2,T,\mathcal{E}_0),\label{highere}
\end{align}
where $C>0$ and $C(r_1,r_2,T,\mathcal{E}_0)>0$ are two constants
independent of $\varepsilon,$ but may depend on $(\gamma,N)$ and $(r_1,r_2,T,\mathcal{E}_0)$ respectively,
and $\R^2_+=:\{(t,r)\,:\, t\in(0,\infty),r\in(0,\infty)\}$.

\smallskip
\item[\rm (ii)] The following energy inequality holds{\rm :}
\begin{align}\label{ei}
&\int_{\R^N}\big(\frac{1}{2}\Big|\frac{\mathcal{M}^{\varepsilon}}{\sqrt{\rho^{\varepsilon}}}\Big|^2+e(\rho^{\varepsilon},\rho_{\ast})+\frac{1}{2}|\nabla_{\mathbf{x}}\Phi^{\varepsilon}|^2\big)(t,\mathbf{x})\dd \mathbf{x}\nonumber\\
&\leq\int_{\R^N}\big(\frac{1}{2}\Big|\frac{\mathcal{M}^{\varepsilon}_0}{\sqrt{\rho^{\varepsilon}_0}}\Big|^2+e(\rho^{\varepsilon}_0,\rho_{\ast})+\frac{1}{2}|\nabla_{\mathbf{x}}\Phi^{\varepsilon}_0|^2\big)(\mathbf{x})\dd \mathbf{x}\qquad  \text{ for $t\geq0$}.
\end{align}

\smallskip
\item[\rm (iii)] Let $(\eta,q)$ be an entropy pair defined by \eqref{weakentropy} for a smooth compactly supported function
$\psi(s)$ on $\R.$
Then, for $\varepsilon\in(0,1],$
\begin{align}
\partial_t\eta^{\psi}(\rho^{\varepsilon},m^{\varepsilon})+\partial_rq^{\psi}(\rho^{\varepsilon},m^{\varepsilon})
\quad \text{ is compact in } H^{-1}_{\rm loc}(\R^2_+),
\end{align}
where $H^{-1}_{\rm loc}(\R^2_+)$ denotes $H^{-1}((0,T)\times \Omega)$ for any $T>0$ and open subset $\Omega\Subset \R_+.$
\end{enumerate}

\medskip
\noindent
\textbf{Part 2. Inviscid Limit and Existence of Global Solutions of CEPEs \eqref{1.1}}{\rm :} For the global weak solutions
$(\rho^{\varepsilon},\mathcal{M}^{\varepsilon},\Phi^{\varepsilon})$
of CNPSEs \eqref{nsp} obtained in {\rm Part 1},
there exist a subsequence $($still denoted$)$ $(\rho^{\varepsilon},m^{\varepsilon},\Phi^{\varepsilon})$ and a vector function $(\rho,m,\Phi)$ such that, as $\varepsilon\rightarrow0+,$
\begin{align}
&(\rho^{\varepsilon},m^{\varepsilon})\rightarrow(\rho,m)(t,r)\,\quad
\mbox{ in $(L^p_{\rm loc}\times L^q_{\rm loc})([0,\infty); r^{N-1}\dd r)$
$\,\,\,$ for $p\in [1,\gamma+1)$  and $q\in[1,\frac{3(\gamma+1)}{\gamma+3})$}, \nonumber\\
&\Phi^{\varepsilon}\rightharpoonup\Phi\quad
\text{ weakly in $L^2(0,T; H^1_{\rm loc}(\R^N))$},\nonumber\\
&\Phi^{\varepsilon}_r(t,r)r^{N-1}
=-\int^r_0\big(\rho^{\varepsilon}(t,y)-d(y)\big)\,y^{N-1}\dd y\nonumber\\
&\quad\quad\quad\quad\quad\rightarrow\,
\Phi_r(t,r)r^{N-1}=-\int^r_0\big(\rho(t,y)-d(y)\big)\,y^{N-1}\dd y  \qquad a.e.\,\, (t,r)\in \R^2_+,\nonumber\\
&\int^T_0\int^{r_2}_0|(\frac{m^{\varepsilon}}{\sqrt{\rho^{\varepsilon}}})(t,r)-(\frac{m}{\sqrt{\rho}})(t,r)|^2\,r^{N-1}\dd r\dd t\rightarrow0
\quad \text{ for any fixed } T,r_2 \in (0,\infty)\nonumber,
\end{align}
and
$(\rho,\mathcal{M},\Phi)(t,\mathbf{x}):=(\rho(t,r),m(t,r)\frac{\mathbf{x}}{r},\Phi(t,r))$ to be a global spherically symmetric solution with finite relative-energy
of problem \eqref{1.1}--\eqref{asym}
in the sense of {\rm Definition \ref{weakep}}.
\end{theorem}

\section{Constructions of Approximate Solutions and Uniform Estimates}
In this section, we construct appropriate approximate solutions
and obtain some uniform estimates of the approximate solutions
of the Cauchy problem for CNSPEs \eqref{1.3} that corresponds to
the initial conditions \eqref{appinitial}--\eqref{phiapproximate}.

\subsection{Constructions of Approximate Solutions}
In this subsection, we consider the approximate CNSPEs in the truncated domains:
\begin{align}\label{6.1}
\begin{cases}
\r_t+(\r u)_r+\frac{N-1}{r}\rho u=0,\\[1mm]
(\r u)_t+(\r u^2+p(\rho))_r+\frac{N-1}{r}\rho u^2-\frac{\rho}{r^{N-1}}\int^r_{\delta}\big(\rho(t,y)-d(y)\big)\,y^{N-1}\dd y\\[0.5mm]
\quad\quad\quad\quad\quad\quad\quad\quad\quad\quad\quad\quad\quad\quad=\v\big((\mu+\lambda)(u_r+\frac{N-1}{r}u)\big)_r-\v\frac{N-1}{r}\mu_ru,
\end{cases}
\end{align}
where $r\in [\delta, b], t>0$,  with $\delta\in(0,1]$ and $b\geq1+\delta^{-1}$. We impose \eqref{6.1} with the following initial data:
\begin{align}\label{6.2}
(\rho, u)(0,r)
=(\r^{\varepsilon,\delta,b}_0, u^{\varepsilon,\delta,b}_0)(r)
\qquad \text{ for $r\in[\delta,b]$},
\end{align}
and the boundary condition:
\begin{align}\label{6.3}
u(t,\delta)=u(t,b)=0 \qquad \text{for $t>0$}.
\end{align}
We take
\begin{equation}\label{3.3a}
\mu(\rho)=\rho+\delta \rho^{\alpha},\quad \lambda(\rho)=\delta(\alpha-1)\rho^{\alpha}
\end{equation}
with $\alpha=\frac{2N-1}{2N}$
so that $\mu(\r)$ and $\lambda(\r)$ satisfy the following BD relation:
\begin{align}\label{miulamuda}
\mu(\r)+\lambda(\r)=\r \mu'(\r).
\end{align}
We first fix parameters $\v>0, \d>0$, and $b\ge 1+\d^{-1}$, and
assume that $(\r^{\varepsilon,\delta,b}_0, u^{\varepsilon,\delta,b}_0)$
are smooth functions such that
\begin{equation}\label{upperlower}
0<({\beta\varepsilon})^{\frac{1}{4}}\leq \r_0^{\varepsilon,\delta,b}\leq (\beta\varepsilon)^{-\frac{1}{2}}<\infty
\end{equation}
for a given small constant $0<\beta\ll 1$.
The construction of approximate initial data functions in \eqref{6.2}
have been established in Appendix A, which satisfy all the properties
in Lemma \ref{lem8.4}.

The existence of global smooth solutions of the initial-boundary value problem \eqref{6.1}--\eqref{6.3}
can be achieved by following similar arguments as in \S 3 and \S 4.1
in Guo-Jiu-Xin \cite{Guo-Jiu-Xin-2}; see also \cite{Jiang}.
Notice that the lower and upper bounds of $\r^{\v,\d,b}$ in \cite{Guo-Jiu-Xin-2} depend on parameters $(\v,\d, b).$
Therefore, the main strategy of this section is to obtain some uniform estimates of $(\r^{\v,\d,b}, u^{\v,\d,b})$ independent of $(\d, b)$
such that we can take both limits $b\rightarrow\infty$ and $\d\rightarrow0+$
to obtain the global existence of weak solution of \eqref{nsp} and \eqref{appinitial}--\eqref{phiapproximate};
see \S 5.

Since $\rho^{\v,\d,b}$ is positive, we may use $u^{\v,\d,b}$ or $m^{\v,\d,b}$ alternatively
,
and drop the superscripts of
solution $(\rho^{\v,\d, b}, u^{\v,\d,b})(t,r)$
and the approximate initial data $(\rho^{\v,\d, b}_0, u^{\v,\d,b}_0)$
for simplicity.
We keep the superscripts when the initial data functions are involved.
After solving \eqref{6.1}--\eqref{6.3},
we can define the potential function $\Phi$ as the solution of
the Poisson equation:
\begin{equation}\label{poisson}
-\Delta\Phi=(\rho-d(\mathbf{x}))\mathbf{I}_{\Omega},
\qquad \lim_{|\mathbf{x}|\rightarrow\infty}\Phi=0,
\end{equation}
with $\Omega=\{\mathbf{x}\in \R^N\,:\,\delta\leq|\mathbf{x}|\leq b\}.$
We can obtain that $\Phi(t,\mathbf{x})=\Phi(t,r)$ with
\begin{align}\label{phib}
\Phi_r(t,r)=\begin{cases}
0 & \text{ for } r\in [0,\delta],\\[1mm]
-\frac{1}{r^{N-1}}\int^r_{\delta}\big(\rho(t,y)-d(y)\big)\,y^{N-1}\dd y
&\text{ for } r\in [\delta,b],\\[1mm]
-\frac{1}{r^{N-1}}\int^b_{\delta}\big(\rho(t,y)-d(y)\big)\,y^{N-1}\dd y
&\text{ for } r\in [b,\infty).
\end{cases}
\end{align}
Using the conservation law of mass, we obtain
\begin{equation}\label{mc}
-\int^b_{\delta}\big(\rho(t,y)-d(y)\big)\,y^{N-1}\dd y
=-\int^b_{\delta}\big(\rho_0(y)-d(y)\big)\,y^{N-1}\dd y:=\frac{M_b}{\omega_N}.
\end{equation}


\subsection{Uniform Estimates of the Approximate Solutions}
The main goal of this section is to obtain some uniform estimates that are independent of $(\v,\d, b)$ so that the two limits
$b\rightarrow\infty$ and $\delta\rightarrow0+$ can be taken in order.

\begin{lemma}[\bf Basic Energy Estimate]\label{bee}
The smooth solution of \eqref{6.1}--\eqref{6.3} satisfies the following
basic energy estimate{\rm :}
\begin{align*}
	&\int^b_{\delta}\Big(\frac{1}{2}\rho u^2+e(\rho,\rho_{\ast})+\frac{1}{2}|\Phi_r|^2\Big)(t,r)\,r^{N-1}\dd r
 +\varepsilon\int^t_0\int^b_{\delta}\rho\Big( u^2_r+\frac{N-1}{r^2} u^2\Big)(s,r)\,r^{N-1}\dd r\dd s\nonumber\\
&\quad+\varepsilon\delta \int^t_0\int^b_{\delta}\rho^{\alpha}\Big(\alpha u^2_r+2(\alpha-1)(N-1)\frac{uu_r}{r}+\big(1+(N-1)(\alpha-1)\big)(N-1)\frac{u^2}{r^2}\Big)(s,r)\,r^{N-1}\dd r\dd s\nonumber\\
&=\int^b_{\delta}\Big(\frac{1}{2}\rho_0 u^2_0+e(\rho_0,\rho_{\ast})+\frac{1}{2}|\Phi_{0r}|^2\Big)(r)\,r^{N-1}\dd r=:\frac{\mathcal{E}^{\varepsilon,\delta,b}_0}{\omega_N}.
\end{align*}
In particular, there exists a constant $c_N>0$ depending only on $N$
such that, for the smooth solution of \eqref{6.1}--\eqref{6.3},
\begin{align}
&\int_{\delta}^{b} \Big(\frac12\rho u^2+e(\rho,\r_{\ast})+\frac{1}{2}|\Phi_r|^2 \Big)(t,r)
 \,r^{N-1} \dd r
 +c_N\v\delta\int^t_0\int^b_{\delta}
 \Big(\rho^{\alpha}u^2_r+\frac{\rho^{\alpha}u^2}{r^2}\Big)(s,r)\,r^{N-1}\dd r\dd s\nonumber\\
&\,\,+\v\int_0^t\int_\delta^b \big(\rho u_r^2+\frac{\rho u^2}{r^2}\big)(s,r)r^{N-1}\,\dd r \dd s
\leq \frac{\mathcal{E}^{\varepsilon,\delta,b}_0}{\omega_N}
\leq C(\mathcal{E}_0+1)\qquad \mbox{ for any $t>0$}\label{6.11}
\end{align}
for some constant $C>0$ independent of $(\varepsilon,\delta,b).$
\end{lemma}

\noindent{\bf Proof.}
Multiplying $\eqref{6.1}_2$ by $r^{N-1}u$ and then integrating by parts,
we have
\begin{align}\label{m1}
&\frac{\dd }{\dd t}\int^b_{\delta}\frac{1}{2}\rho u^2\,r^{N-1}\dd r+\int^{b}_{\delta}p_ru\, r^{N-1}\dd r\nonumber\\
&=-\varepsilon \int^b_{\delta}\Big((\mu+\lambda)(u_r+\frac{N-1}{r}u)(r^{N-1}u)_r-(N-1)\mu(r^{N-2}u^2)_r\Big)\,\dd r\nonumber\\
&\quad+\int^b_{\delta}\rho u\int^r_{\delta}\big(\rho(t,y)-d(y)\big)\,y^{N-1}\dd y\dd r.
\end{align}
For the second term on the left-hand side (LHS) of \eqref{m1}, one has
\begin{equation*}
\begin{split}
\int^b_{\delta}p_ru\,r^{N-1}\dd r
&=\frac{\kappa \gamma}{\gamma-1}\int^b_{\delta}\rho u(\rho^{\gamma-1})_r\,r^{N-1}\dd r=-\frac{\kappa\gamma}{\gamma-1}\int^{b}_{\delta}(r^{N-1}\rho u)_r\rho^{\gamma-1}\,\dd r\\
&=\frac{\kappa}{\gamma-1}\int^b_{\delta}(\rho^{\gamma})_t\,r^{N-1}\dd r=\frac{\kappa}{\gamma-1}\int^b_{\delta}\Big(\rho^{\gamma}-\rho_{\ast}^{\gamma}-\gamma{\rho_{\ast}}^{\gamma-1}(\rho-\rho_{\ast})\Big)_t\,r^{N-1}\dd r\\
&=\frac{\dd }{\dd t}\int^b_{\delta}e(\rho,\rho_{\ast})(t,r)\,r^{N-1}\dd r.
\end{split}
\end{equation*}
For the viscous term, a direct calculation shows that
\begin{align}\label{viscous}
&(\mu+\lambda)(u_r+\frac{N-1}{r})(r^{N-1}u)_r-(N-1)\lambda(r^{N-2}u^2)_r\nonumber\\[1mm]
&=\delta\rho^{\alpha}\Big\{\alpha r^{N-1}u^2_r
+2(\alpha-1)(N-1)r^{N-2}uu_r+\big(N-1+(\alpha-1)(N-1)^2\big)r^{N-3}u^2\Big\}\nonumber\\[1mm]
&\quad+\rho\big(r^{N-1}u^2_r+(N-1)r^{N-3}u^2\big).
\end{align}
For the first term on the right-hand side (RHS) of \eqref{viscous},
noting that $\alpha\in(\frac{N-1}{N},1),$
we calculate its discriminant:
\begin{equation}\nonumber
\begin{split}
\Delta&=4(\alpha-1)^2(N-1)^2-4\alpha\big(N-1+(\alpha-1)(N-1)^2\big)\\
&=4(N-1)^2\big((\alpha-1)^2-\alpha(\alpha-1)\big)-4(N-1)\alpha\\
&=4(N-1)^2(1-\alpha)-4(N-1)\alpha\\
&=4(N-1)^2\big(1-\frac{N}{N-1}\alpha\big)<0 \qquad
\text{ if $\alpha\in (\frac{N-1}{N},1)$},
\end{split}
\end{equation}
so that there exists a constant $c_N>0$ such that
\begin{equation}\nonumber
\begin{split}
&(\mu+\lambda)\big(u_r+\frac{N-1}{r}\big)(r^{N-1}u)_r-(N-1)\mu(r^{N-2}u^2)_r\\
&\geq c_N \delta \rho^{\alpha}\big(r^{N-1}u^2_r+r^{N-3}u^2\big)
+\rho\big(r^{N-1}u^2_r+r^{N-3}u^2\big).
\end{split}
\end{equation}
For the last term on RHS of \eqref{m1}, a direct calculation yields
\begin{equation*}
\begin{split}
&\int^{b}_{\delta}\rho u\int^r_{\delta}
  \big(\rho(t,y)-d(y)\big)\,y^{N-1}\dd y\dd r\\
&=-\int^b_{\delta}\rho u\Phi_r(t,r)\,r^{N-1}\dd r
=\int^{b}_{\delta}(\rho ur^{N-1})_r\Phi\,\dd r\\
&=-\int^b_{\delta}(r^{N-1}\rho)_t\Phi\,\dd r
=-\int^b_{\delta}\big(r^{N-1}(\rho-d(r))\big)_t\Phi\,\dd r\\
&=\int^{b}_{\delta}\big(\Phi_{rr}+\frac{N-1}{r}\Phi_r\big)_t\Phi\,r^{N-1}\dd r
=\int^{b}_{\delta}(r^{N-1}\Phi_{rt})_r\Phi\,\dd r\\
&=-\int^b_{\delta}\Phi_{rt}\Phi_r\,r^{N-1}\dd r
+(r^{N-1}\Phi_{rt}\Phi)\Big|^b_{\delta}
=-\int^{b}_{\delta}\big(\frac{|\Phi_r|^2}{2}\big)_t\,r^{N-1}\dd r\\
&=-\frac{1}{2}\frac{\dd }{\dd t}\int^b_{\delta}|\Phi_r|^2\,r^{N-1}\dd r,
\end{split}
\end{equation*}
where we have used \eqref{phib} and
\begin{equation}\label{phi2}
\begin{split}
\Phi_{rt}(t,\delta)=\Phi_{rt}(t,b)=0\qquad\,\, \mbox{for any $t\geq0$}.
\end{split}
\end{equation}
Combining \eqref{m1}--\eqref{phi2}, we obtain
\begin{align}\label{viscous2}
&\frac{\dd }{\dd t}\Big\{\int^b_{\delta}\Big(\frac{1}{2}\rho u^2+e(\rho,\rho_{\ast})+\frac{1}{2}|\Phi_r|^2\Big)\,r^{N-1}\dd r\Big\}\nonumber\\
&\quad+\varepsilon \int^b_{\delta}
\big(c_N\delta\rho^{\alpha}(r^{N-1}u^2_r+r^{N-3}u^2)+\rho(r^{N-1}u^2_r+r^{N-3}u^2)\big)\dd r\leq0.
\end{align}
Integrating \eqref{viscous2}, we have
\begin{equation*}\label{viscous3}
\begin{split}
&\int^b_{\delta}\Big(\frac{1}{2}\rho u^2+e(\rho,\rho_{\ast})+\frac{1}{2}|\Phi_r|^2\Big)\,r^{N-1}\dd r\\
&\quad+\varepsilon\int^t_0\int^{b}_{\delta}
\big(\rho(r ^{N-1}u^2_r+r^{N-3}u^2)+c_N\delta\rho^{\alpha}(r^{N-1}u^2_r+r^{N-3}u^2)\big)\,\dd r\dd s\\
&\leq\int^b_{\delta}\Big(\frac{1}{2}\rho_0 u^2_0+e(\rho_0,\rho_{\ast})+\frac{1}{2}|\Phi_{0r}|^2\Big)\,r^{N-1}\dd r.
\end{split}
\end{equation*}
Then we conclude Lemma \ref{bee}.
$\hfill\Box$

\begin{remark}
Using \eqref{mc}, we can rewrite the electric field term as
\begin{align}\label{electric1}
\|\nabla\Phi\|^2_{L^2(\R^N)}&=\omega_N\int^{\infty}_0|\Phi_r|^2\,r^{N-1}\dd r\nonumber\\
&=\omega_N\Big(\int^{b}_{\delta}|\Phi_r|^2\,r^{N-1}\dd r
+\int^{\infty}_{b}|\Phi_r|^2\,r^{N-1}\dd r\Big)\nonumber\\
&=\omega_N\Big(\int^{b}_{\delta}|\Phi_r|^2\,r^{N-1}\dd r
+\int^{\infty}_b\frac{M^2_b}{\omega^2_Nr^{N-1}}\,\dd r\Big)\nonumber\\
&=\omega_N\Big(\int^{b}_{\delta}|\Phi_r|^2\,r^{N-1}\dd r
+\frac{1}{N-2}(\frac{M_b}{\omega_N})^2\frac{1}{b^{N-2}}\Big).
\end{align}
On one hand, it follows from \eqref{electric1} that this term has
a good sign, which is important.
On the other hand, we also have the following equality:
\begin{equation*}\label{electric2}
\begin{split}
&\int^b_{\delta}(\rho-d(r))r\int^r_{\delta}\big(\rho(t,y)-d(y)\big)\,y^{N-1}
\dd y\dd r\\
&=-\int^b_{\delta}\big(\rho-d(r)\big)\Phi_r\,r^N\dd r
=\int^b_{\delta}\big(\Phi_{rr}+\frac{N-1}{r}\Phi_r\big)\Phi_r\,r^N\dd r\\
&=\int^b_{\delta}(r^{N-1}\Phi_r)_r\Phi_r\,r\dd r
=\int^b_{\delta}\frac{1}{r^{N-2}}(\frac{1}{2}|r^{N-1}\Phi_r|^2)_r\,\dd r\\
&=-\int^b_{\delta}\frac{1}{2}(r^{-N+2})_r|r^{N-1}\Phi_r|^2\,\dd r
+\frac{1}{2r^{N-2}}\big|r^{N-1}\Phi_r\big|^2\Big|^b_{\delta}\\
&=\frac{N-2}{2}\int^b_{\delta}|\Phi_r|^2\,r^{N-1}\dd r
+\frac{1}{2b^{N-2}}|b^{N-1}\Phi_r(t,b)|^2\\
&=\frac{N-2}{2}\int^b_{\delta}r^{N-1}|\Phi_r|^2\dd r+\frac{1}{2b^{N-2}}(\frac{M_b}{\omega_N})^2\\
&=\frac{N-2}{2}\Big(\int^b_{\delta}|\Phi_r|^2\,r^{N-1}\dd r
+\frac{1}{(N-2)b^{N-2}}(\frac{M_b}{\omega_N})^2\Big)\geq0.
\end{split}
\end{equation*}
\end{remark}

\begin{lemma}[\bf BD Entropy Estimate]\label{BD}
The smooth solution of \eqref{6.1}--\eqref{6.3} satisfies
\begin{align}\label{6.1222}
&\varepsilon^2\int_\delta^b\Big((1+\alpha\delta \rho^{\alpha-1}+\alpha^2\delta^2\rho^{2(\alpha-1)})\frac{\rho^2_r}{\rho}\Big)(t,r)\,r^{N-1}\dd r\nonumber\\
&\,\,+\v \int_0^t\int_\delta^b
\big((1+\alpha\delta\rho^{\alpha-1})\rho^{\gamma-2}\rho^2_r\big)(s,r)\,r^{N-1}\dd r\dd s\leq C(\mathcal{E}_0+1),
\end{align}
where we have used the following relation{\rm :}
\begin{align}\label{bound}
&\sup_{0<\varepsilon,\delta\leq1}\sup_{b\geq1+\delta^{-1}}\big(\mathcal{E}^{\varepsilon,\delta,b}_0
+\mathcal{E}^{\varepsilon,\delta,b}_1\big)\nonumber\\
&\quad+CT\int^b_{\delta}\big(|\rho_{\ast}-d(r)|+|\rho_{\ast}-d(r)|^2+|\rho_{\ast}-d(r)|^{\frac{\gamma}{\gamma-1}}\big)\,r^{N-1}\dd r\leq C(\mathcal{E}_0+1)
\end{align}
with
	\begin{align}\label{6.13-1}
	\mathcal{E}^{\varepsilon,\delta,b}_1
 :=\v^2 \int^b_{\delta}
 \big(1+2\alpha\delta\rho^{\alpha-1}_0+\alpha^2\delta^2\rho^{2\alpha-2}_0\big)|(\sqrt{\rho_0})_r|^2\,r^{N-1} \dd r.
	\end{align}
\end{lemma}

\noindent{\bf Proof.}
We divide the proof into seven steps:

\smallskip
\noindent\textbf{1.}
It is convenient to deal with \eqref{6.1} in the Lagrangian coordinates.
For simplicity, we assume that
\begin{equation}\nonumber
L_b:=\int_{\delta}^b\rho_0(r)\,r^{N-1}\dd r.
\end{equation}
Since
$$
\frac{\dd }{\dd t}\int_{\delta}^{b}\rho(t,r)\,r^{N-1}\dd r
=-\int^{b}_{\delta}(r^{N-1}\rho u)_r(t,r)\,\dd r=0,
$$
we have
\begin{equation}\nonumber
\int_{\delta}^{b}\rho(t,r)\,r^{N-1}\dd r
=\int_{\delta}^b\rho_0(r)\,r^{N-1}\dd r=L_b \qquad
\text{ for all $t>0$}.
\end{equation}
For $r\in[\delta,b]$ and $t\in[0,T]$,
we define the Lagrangian coordinates transformation:
\begin{equation*}\label{l6.7}
x=\int_{\delta}^r \rho(t,y)\,y^{N-1}\dd y,\qquad \tau=t,
\end{equation*}
which translates domain $[0,T]\times[\delta,b]$ into
$[0,T]\times[0,L_b]$ and satisfies
\begin{align*}
\begin{cases}
\frac{\partial x}{\partial r}=\rho r^{N-1},\quad\, \frac{\partial x}{\partial t}=-\rho ur^{N-1},\quad\, \frac{\partial\tau}{\partial r}=0,\quad\, \frac{\partial\tau}{\partial t}=1,\\
\frac{\partial r}{\partial x}=\frac{1}{\rho r^{N-1}},\quad\, \frac{\partial r}{\partial\tau}=u,\quad\, \frac{\partial t}{\partial\tau}=1,\quad\, \frac{\partial t}{\partial x}=0.
\end{cases}
\end{align*}
Applying the Lagrange transformation to \eqref{6.1}, we have
\begin{align}\label{6.9}
\begin{cases}
\rho_\tau+\rho^2(r^{N-1}u)_x=0,\\[1mm]
u_\tau+r^{N-1}p_x=\v r^{N-1}\big(\rho(\mu+\lambda)(r^{N-1}u)_x\big)_x
-\varepsilon(N-1)r^{N-2}\mu_xu
+\frac{1}{r^{N-1}}\Big(x-\int^r_{\delta}d(y)\,y^{N-1}\dd y\Big),
\end{cases}
\end{align}
and the boundary conditions \eqref{6.3} becomes
$$
u(\tau,0)=u(\tau,L_b)=0 \qquad\text{ for $\tau>0$}.
$$

\smallskip
\noindent\textbf{2.}
Multiplying $\eqref{6.9}_1$ by $\mu'(\rho)$
and using \eqref{miulamuda}, we have
\begin{equation}\label{6.91}
\mu_{\tau}+\rho(\mu+\lambda)(r^{N-1}u)_x=0.
\end{equation}
Then we substitute \eqref{6.91} into $\eqref{6.9}_2$ to obtain
\begin{equation}\label{6.92}
u_\tau+r^{N-1}p_x
=-\varepsilon r^{N-1}\mu_{x\tau}-\varepsilon(N-1)r^{N-2}\mu_xu
+\frac{1}{r^{N-1}}\Big(x-\int^r_{\delta}d(y)\,y^{N-1}\dd y\Big).
\end{equation}
Noting that $u=\frac{\partial r}{\partial\tau}$, we have
$$
\varepsilon(N-1)r^{N-2}\mu_xu=\varepsilon(N-1)r^{N-2}\mu_xr_{\tau}=\varepsilon (r^{N-1})_{\tau}\mu_x,
$$
which, together with \eqref{6.92}, yields
\begin{equation}\label{6.93}
(u+\varepsilon r^{N-1}\mu_x)_{\tau}+r^{N-1}p_x
=\frac{1}{r^{N-1}}\Big(x-\int^r_{\delta}d(y)\,y^{N-1}\dd y\Big).
\end{equation}

\noindent\textbf{3.}
Multiplying \eqref{6.93} by $u+\varepsilon r^{N-1}\mu_x,$ we see that
\begin{align}\label{6.94}
&\frac{\dd }{\dd \tau}\int^{L_b}_0
\frac{1}{2}\big|u+\varepsilon r^{N-1}\mu_x\big|^2\,\dd x
+\int^{L_b}_0 p_x\big(u+\varepsilon r^{N-1}\mu_x\big)\,r^{N-1}\dd x\nonumber\\
&=\int^{L_b}_0\frac{1}{r^{N-1}}\Big(x-\int^r_{\delta}d(y)\,y^{N-1}\dd y\Big)
\big(u+\varepsilon r^{N-1}\mu_x\big)\,\dd x.
\end{align}
Integrating \eqref{6.94} over $[0,\tau]$ and pulling it back
to the Eulerian coordinates, we have
\begin{align}\label{6.95}
&\int^b_{\delta}\frac{1}{2}\rho \Big|u+\varepsilon\frac{\mu_r}{\rho}\Big|^2\,r^{N-1}\dd r
+\varepsilon \int^t_0\int^b_{\delta}\frac{p_r\mu_r}{\rho}\,r^{N-1}\dd r\dd s+\int^t_0\int^b_{\delta}up_r\,r^{N-1}\dd r\dd s\nonumber\\
&=\int^b_{\delta}\frac{1}{2}\rho_0
\Big|u_0+\varepsilon\frac{\mu_{0r}}{\rho_0}\Big|^2r^{N-1}\dd r
+\int^t_0\int^b_{\delta}\rho u
\Big(\int^r_{\delta}\big(\rho(s,y)-d(y)\big)\,y^{N-1}\dd y\Big)\dd r\dd s\nonumber\\
&\quad+\varepsilon\int^t_0\int^b_{\delta}\mu_r
\Big(\int^r_{\delta}\big(\rho(s,y)-d(y)\big)\,y^{N-1}\dd y\Big)\dd r\dd s.
\end{align}

\noindent\textbf{4.} A direct calculation shows that
\begin{align}
&\int^t_0\int^{b}_{\delta}up_r\,r^{N-1}\dd r\dd s
=\int^{b}_{\delta}e(\rho,\rho_{\ast})\,r^{N-1}\dd r
-\int^b_{\delta}e(\rho_0,\rho_{\ast})\,r^{N-1}\dd r,\label{6.96}\\[1mm]
&\int^t_0\int^{b}_{\delta}\rho u
\Big(\int^r_{\delta}\big(\rho(s,y)-d(y)\big)\,y^{N-1}\dd y\Big)\dd r\dd s
=-\frac{1}{2}\int^b_{\delta}|\Phi_r|^2\,r^{N-1}\dd r
+\frac{1}{2}\int^b_{\delta}|\Phi_{0r}|^2\,r^{N-1}\dd r.\label{6.97}
\end{align}

\smallskip
\noindent\textbf{5.}
For the last term of RHS of \eqref{6.95},
we integrate by parts to obtain
\begin{align}\label{6.98}
&\varepsilon\int^t_0\int^b_{\delta}\mu_r
\Big(\int^r_{\delta}\big(\rho(s,y)-d(y)\big)\,y^{N-1}\dd y\Big)\dd r\dd s\nonumber\\
&=\varepsilon\int^t_0\int^b_{\delta}\big(\mu(\rho)-\mu(\rho_\ast)\big)_r
\Big(\int^r_{\delta}\big(\rho(s,y)-d(y)\big)\,y^{N-1}\dd y\Big)\dd r\dd s\nonumber\\
&=-\varepsilon\int^t_0\int^b_{\delta}\big(\mu(\rho)-\mu(\rho_{\ast})\big)
\big(\rho(s,r)-d(r)\big)\,r^{N-1}\dd r\dd s\nonumber\\
&\quad+\varepsilon\int^t_0\big(\mu(\rho(s,b))-\mu(\rho_{\ast})\big)
\Big(\int^b_{\delta}\big(\rho(s,y)-d(y)\big)\,y^{N-1} \dd y\Big)\dd s\nonumber\\
&=-\varepsilon\int^t_0\int^b_{\delta}\big(\mu(\rho)-\mu(\rho_{\ast})\big)\big(\rho(s,r)-\rho_{\ast}\big)\,r^{N-1}\dd r\dd s\nonumber\\
&\quad-\varepsilon\int^t_0\int^b_{\delta}\big(\mu(\rho)-\mu(\rho_{\ast})\big)\big(\rho_{\ast}-d(r)\big)\,r^{N-1}\dd r\dd s\nonumber\\
&\quad+\varepsilon\int^t_0\big(\mu(\rho(s,b))-\mu(\rho_{\ast})\big)
\Big(\int^b_{\delta}\big(\rho(s,y)-d(y)\big)\,y^{N-1}\dd y\Big)\dd s\nonumber\\
&\leq-\varepsilon\int^t_0\int^b_{\delta}
\big(\rho_{\ast}-d(r)\big)\big(\mu(\rho(s,r))-\mu(\rho_{\ast})\big)\,r^{N-1}\dd r\dd s\nonumber\\
&\quad+\Big|\frac{M_b}{\omega_N}\Big|\int^t_0\varepsilon\big|\mu(\rho(s,b))-\mu(\rho_{\ast})\big|\,\dd s.
\end{align}
The first term of RHS of \eqref{6.98} becomes
\begin{align}\label{6.99}
&-\varepsilon\int^t_0\int^b_{\delta}\big(\rho_{\ast}-d(r)\big)\big(\mu(\rho(s,r))-\mu(\rho_{\ast})\big)\,r^{N-1}\dd r\dd s\nonumber\\
&\leq C\varepsilon\int^t_0\int^b_{\delta}
\big|\rho_{\ast}-d(r)\big|\,r^{N-1}\dd r\dd s\nonumber\\
&\quad+\varepsilon\int^t_0\int^b_{\delta}
\big|\rho_{\ast}-d(r)\big|\mathbf{1}_{\{\rho(s,r)\geq2\rho_{\ast}\}}
\big|\mu(\rho(s,r))-\mu(\rho_{\ast})\big|\,r^{N-1}\dd r\dd s\nonumber\\
&\leq CT\int^b_{\delta}\big|\rho_{\ast}-d(r)\big|\,r^{N-1}\dd r\nonumber\\
&\quad+C\int^T_0\int^b_{\delta}
\Big(\big|\mu(\rho(s,r))-\mu(\rho_{\ast})\big|^{\gamma}
+\big|\rho_{\ast}-d(r)\big|^{\frac{\gamma}{\gamma-1}}\Big)\mathbf{1}_{\{\rho(s,r)\geq2\rho_{\ast}\}}\,r^{N-1}\dd r\dd s\nonumber\\
&\leq CT\int^b_{\delta}\Big(\big|\rho_{\ast}-d(r)\big|
+\big|\rho_{\ast}-d(r)\big|^{\frac{\gamma}{\gamma-1}}\Big)\,r^{N-1}\dd r
+C\int^T_0\int^b_{\delta}r^{N-1}e(\rho,\rho_{\ast})\,r^{N-1}\dd r\dd s\nonumber\\
&\leq CT\int^b_{\delta}\Big(\big|\rho_{\ast}-d(r)\big|
+\big|\rho_{\ast}-d(r)\big|^{\frac{\gamma}{\gamma-1}}\Big)\,r^{N-1}\dd r
+C\mathcal{E}^{\varepsilon,\delta,b}_0\leq C(1+\mathcal{E}_0),
\end{align}
where $\mathbf{1}_{\{\rho(s,r)\geq2\rho_{\ast}\}}$ is the indicator function of set $\{\rho(s,r)\geq2\rho_{\ast}\}.$

\noindent\textbf{6.}
To bound the last term of RHS of \eqref{6.98}, we have to be very careful, since it involves the boundary value of the density.
\begin{equation}\label{6.991}
\begin{split}
&\varepsilon\int^t_0|\mu(\rho(s,b))-\mu(\rho_{\ast})|\,\dd s\\
&=\varepsilon\int^t_0\mathbf{1}_{\{\rho(s,b)\leq4\rho_{\ast}\}}|\mu(\rho(s,b))-\mu(\rho_{\ast})|\,\dd s
+\varepsilon\int^t_0\mathbf{1}_{\{\rho(s,b)>4\rho_{\ast}\}}|\mu(\rho(s,b))-\mu(\rho_{\ast})|\,\dd s\\
&\leq CT+\varepsilon\int^t_0\mathbf{1}_{\{\rho(s,b)>4\rho_{\ast}\}}|\mu(\rho(s,b))-\mu(\rho_{\ast})|\,\dd s.
\end{split}
\end{equation}
Since
$$
\int^b_{\delta}(\rho(t,r)-\rho_{\ast})\,r^{N-1}\dd r
=\int^b_{\delta}(\rho_0(r)-\rho_{\ast})\,r^{N-1}\dd r,
$$
then
$$
\int^b_{\delta}\rho(t,r)\,r^{N-1}\dd r
=\rho_{\ast}\frac{1}{N}(b^N-\delta^N)
+\int^b_{\delta}(\rho_0(r)-\rho_{\ast})\,r^{N-1}\dd r.
$$
Therefore, there exists $r_0\in[\delta,b]$ such that, for any $t\ge 0$,
\begin{equation}\label{6.992}
\rho(t,r_0(t))
=\rho_{\ast}
+\frac{N}{b^N-\delta^N}\int^b_{\delta}(\rho_0(r)-\rho_{\ast})\,r^{N-1}\dd r
\in [\frac{1}{2}\rho_{\ast},2\rho_{\ast}]  \qquad\text{ if $b>>1$},
\end{equation}
where we have used
$$
\Big|\int^b_{\delta}(\rho_0(r)-\rho_{\ast})\, r^{N-1}\dd r\Big|\leq C < \infty.
$$
Motivated by \cite{Chen2020}, we define
\begin{align*}
&A(t)=\{r\,:\,r\in[\delta,b],\,\,\rho(t,r)\geq2\rho_{\ast}\}\\
&A_1(t)=\{r\in A(t)\,:\, r\geq1\},\qquad\,\, A_2(t)=A(t)\backslash A_1(t).
\end{align*}
Similarly, we also define
\begin{align*}
&B(t)=\{r\,:\,r\in[\delta,b],\,\,0\leq\rho(t,r)\leq\frac{1}{2}\rho_{\ast}\},\\
&B_1(t)=\{r\in B(t)\,:\, r\geq1\},\qquad  B_2(t)=B(t)\backslash B_1(t). 
\end{align*}
It is direct to check that
\begin{equation}\label{eee}
e(\rho,\rho_{\ast})=
\begin{cases}
|\rho-\rho_{\ast}|^2 &\,\,\text{ if } \rho\in[\frac{1}{2}\rho_{\ast},2\rho_{\ast}],\\
|\rho-\rho_{\ast}|^{\gamma} &\,\,\text{ if } \rho\in[0,\frac{1}{2}\rho_{\ast}]\cup[2\rho_{\ast},\infty).\\
\end{cases}
\end{equation}
Then there exists a constant $c(\rho_{\ast})>0$ such that
\begin{equation*}\label{elower}
e(\rho,\rho_{\ast})(t,r)\geq c(\rho_{\ast})>0\qquad \mbox{for any $r\in A(t)\cup B(t)$},
\end{equation*}
which, together with \eqref{6.11}, yields
\begin{equation}\label{6.993}
\mathcal{E}^{\varepsilon,\delta,b}_0
\geq \int^b_{\delta}e(\rho,\rho_{\ast})(t,r)\,r^{N-1}\dd r
\geq c(\rho_{\ast})\int_{A_1(t)\cup B_1(t)}r^{N-1}\dd r\geq c(\rho_{\ast})\big(|A_1(t)|+|B_1(t)|\big).
\end{equation}
Hence, we obtain
\begin{equation}\label{6.994}
|A(t)|+|B(t)|\leq|A_1(t)|+|B_1(t)|+2\leq C(\rho_{\ast})\big(1+\mathcal{E}^{\varepsilon,\delta,b}_0\big).
\end{equation}
Since $\rho(t,r)$ is a continuous function of $(t,r)$, for $s\in \{s\,:\, \rho(s,b)>4\rho_{\ast}\}$,
it follows from \eqref{6.992}--\eqref{6.993} that there exists $\tilde{r}_0(s)\in[1,b]$ such that
\begin{equation}\label{6.995}
\begin{cases}
|b-\tilde{r}_0(s)|\leq C(\rho_{\ast})\big(1+\mathcal{E}^{\varepsilon,\delta,b}_0\big),\\[0.5mm]
\rho(s,r)\geq 2\rho_{\ast}\qquad \mbox{for all $r\in[\tilde{r}_0(s),b]$},\\
\rho(s,\tilde{r}_0(s))=2\rho_{\ast}.
\end{cases}
\end{equation}
Thus, using \eqref{6.995}, we have
\begin{align}\label{6.996}
&\varepsilon\int^t_0\mathbf{1}_{\{\rho(s,b)>4\rho_{\ast}\}}|\mu(\rho(s,b))-\mu(\rho_{\ast})|\,\dd s\nonumber\\
&\leq\varepsilon\int^t_0\mathbf{1}_{\{\rho(s,b)>4\rho_{\ast}\}}|\mu(\rho(s,b))-\mu(\rho(s,\tilde{r}_0(s)))|
\,\dd s+C(\rho_{\ast})T\nonumber\\
&\leq\varepsilon\int^t_0\mathbf{1}_{\{\rho(s,b)>4\rho_{\ast}\}}\int^b_{\tilde{r}_0(s)}|\mu_r|\,\dd r\dd s+C(\rho_{\ast})T\nonumber\\
&\leq\varepsilon\int^t_0\mathbf{1}_{\{\rho(s,b)>4\rho_{\ast}\}}\int^b_{\tilde{r}_0(s)}|\rho_r|
 \big(1+\alpha\delta\rho^{\alpha-1}\big)\,\dd r\,\dd s+C(\rho_{\ast})T\nonumber\\
&\leq C(\rho_{\ast})\varepsilon\int^t_0\mathbf{1}_{\{\rho(s,b)>4\rho_{\ast}\}}\int^b_{\tilde{r}_0(s)}|\rho_r|\dd r\dd s+C(\rho_{\ast})T\nonumber\\
&\leq C(\rho_{\ast})\varepsilon\int^t_0\mathbf{1}_{\{\rho(s,b)>4\rho_{\ast}\}}
\Big(\int^b_{\tilde{r}_0(s)}\rho^{\gamma-2}|\rho_r|^2\,\dd r\Big)^{\frac{1}{2}}
\Big(\int^b_{\tilde{r}_0(s)}\rho^{2-\gamma}\,\dd r\Big)^{\frac{1}{2}}\dd s+C(\rho_{\ast})T.
\end{align}
For $\gamma\in(1,2],$
$$
\rho^{2-\gamma}\leq C(1+\rho^{\gamma}),
$$
which, together with \eqref{6.995}, yields
\begin{align}\label{6.997}
\int^b_{\tilde{r}_0(s)}\rho^{2-\gamma}\,\dd r
&\leq C\int^b_{\tilde{r}_0(s)}\big(1+\rho^{\gamma}\big)\,\dd r
\leq C\int^b_{\tilde{r}_0(s)}\big(1+e(\rho,\rho_{\ast})\big)\,\dd r\nonumber\\
&\leq C(\rho_{\ast})\big(1+\mathcal{E}^{\varepsilon,\delta,b}_0\big)
+C\int^b_{\tilde{r}_0(s)}e(\rho,\rho_{\ast})\,r^{N-1}\dd r\nonumber\\
&\leq C(\rho_{\ast})\big(1+\mathcal{E}^{\varepsilon,\delta,b}_0\big).
\end{align}
For $\gamma\in(2,\infty),$ it follows from \eqref{6.995} that
\begin{equation}\label{6.998}
\int^b_{\tilde{r}_0(s)}\rho^{2-\gamma}\,\dd r
\leq C(\rho_{\ast})|b-\tilde{r}_0(s)|\leq C(\rho_{\ast})\big(1+\mathcal{E}^{\varepsilon,\delta,b}_0\big).
\end{equation}
Combining \eqref{6.996}--\eqref{6.998}, we obtain
\begin{align}\label{6.999}
&\varepsilon\int^T_0\mathbf{1}_{\{\rho(s,b)>4\rho_{\ast}\}}|\mu(\rho(s,b))-\mu(\rho_{\ast})|\,\dd s\nonumber\\
&\leq\frac{\kappa\gamma\varepsilon}{8}\int^t_0\int^b_{\tilde{r}_0(s)}\rho^{\gamma-2}|\rho_r|^2\,\dd r\dd s+C(\rho_{\ast})T\big(1+\mathcal{E}^{\varepsilon,\delta,b}_0\big)\nonumber\\
&\leq\frac{\kappa\gamma\varepsilon}{8}\int^t_0\int^b_{\delta}r^{N-1}\rho^{\gamma-2}|\rho_r|^2\,\dd r\dd s+C(\rho_{\ast})T\big(1+\mathcal{E}^{\varepsilon,\delta,b}_0\big),
\end{align}
which, together with \eqref{6.991}, implies that
\begin{equation}\label{7.12}
\varepsilon\int^t_0|\mu(\rho(s,b))-\mu(\rho_{\ast})|\,\dd s
\leq\frac{\kappa\gamma\varepsilon}{8}
\int^T_0\int^b_{\delta}\rho^{\gamma-2}|\rho_r|^2\,r^{N-1}\dd r\dd s
+C(\rho_{\ast})T\big(1+\mathcal{E}^{\varepsilon,\delta,b}_0\big).
\end{equation}

\smallskip
\noindent\textbf{7.}
Substituting \eqref{6.99} and \eqref{7.12} into \eqref{6.98}, we have
\begin{equation}\label{7.222}
\begin{split}
&\varepsilon\int^t_0\int^b_{\delta}\mu_r\Big(\int^r_{\delta}\big(\rho(s,y)-d(y)\big)\,y^{N-1}\dd y\Big)\dd r\dd s\\
&\leq\frac{\kappa\gamma\varepsilon}{8}\int^T_0\int^b_{\delta}\rho^{\gamma-2}|\rho_r|^2\,r^{N-1}\dd r\dd s+CT\int^b_{\delta}\big(|\rho_{\ast}-d(r)|+|\rho_{\ast}-d(r)|^2\big)\,r^{N-1}\dd r\\
&\quad+C(\rho_{\ast})T\big(1+\mathcal{E}^{\varepsilon,\delta,b}_0\big).
\end{split}
\end{equation}
Combining \eqref{7.222} with \eqref{6.95}--\eqref{6.97}, we obtain
\begin{equation}\label{7.3}
\begin{split}
&\int^b_{\delta}\frac{1}{2}\rho|u+\varepsilon\frac{\mu_r}{\rho}|^2\,r^{N-1}\dd r
+\frac{\kappa\gamma\varepsilon}{2}\int^t_0\int^b_{\delta}\big(\rho^{\gamma-2}\rho^2_r
+\alpha\delta\rho^{\alpha+\gamma-3}\rho^2_r\big)\,r^{N-1}\dd r\dd s\\
&\quad+\int^b_{\delta}e(\rho,\rho_{\ast})\,r^{N-1}\dd r
+\int^b_{\delta}\frac{1}{2}|\Phi_r|^2\,r^{N-1}\dd r\\
&\leq C(\rho_{\ast},T)\big(\mathcal{E}^{\varepsilon,\delta,b}_0
 +\mathcal{E}^{\varepsilon,\delta,b}_1\big)
 +CT\int^b_{\delta}\big(|\rho_{\ast}-d(r)|+|\rho_{\ast}-d(r)|^2\big)\,\dd r\nonumber\\
&\leq C\big(1+\mathcal{E}_0\big).
\end{split}
\end{equation}
This completes the proof Lemma \ref{BD}.
$\hfill\Box$

\smallskip
As shown in \cite{Chen2010}, the higher integrabilities of the density and the velocity
are important for the $L^p$ compensated compactness framework.
We now prove the higher integrability of the density.

\begin{lemma}[\bf Higher Integrability on the Density]
For given $r_1$ and $r_2$ with $\delta<r_1<r_2<b,$
any smooth solution of problem \eqref{6.1}--\eqref{6.3} satisfies
	\begin{align}\label{6.14}
	\int_0^T\int_K \rho^{\g+1}(t,r)\,\dd r\dd t\leq C(r_1,r_2,T,\mathcal{E}_0),
	\end{align}
where $K$ is any compact subset of $[r_1,r_2].$
\end{lemma}

\noindent{\bf Proof.}
Let $w(r)$ be a smooth compact support function such that $\text{supp}\,w\subseteq (r_1,r_2)$ and $w(r)=1$ for $r\in K.$
Multiplying $\eqref{6.1}_2$ by $w(r),$ we obtain
\begin{align}\label{7.32a}
&(\r u w)_t+\big((\r u^2+p)w\big)_r+\frac{N-1}{r} \r u^2 w\nonumber\\
&=\v\big((\rho+\alpha\delta\rho^{\alpha})(u_r+\frac{N-1}{r}u) w\big)_r
-\v \frac{N-1}{r} u (\rho+\delta\rho^{\alpha})_r w\nonumber\\
&\quad+\big(\r u^2+p-\v(\rho+\alpha\delta\rho^{\alpha})(u_r+\frac{N-1}{r} u)\big)w_r
-\rho w\Phi_r.
\end{align}
Integrating \eqref{7.32a} over $[r_1,r),$  we have
\begin{align}
(\r u^2+p)w&=\v (\rho+\alpha\delta\rho^{\alpha})(u_r+\frac{N-1}{r} u) w
-\v\int_{r_1}^r \frac{N-1}{y} u(\rho+\delta\rho^{\alpha})_y w\,\dd y
-\Big(\int^r_{r_1}\r u w \dd y\Big)_t \nonumber\\
&\quad-\int_{r_1}^r \frac{N-1}{y} \r u^2 w\,\dd y
+\int_{r_1}^r\Big(\r u^2+p-\v(\rho+\alpha\delta\rho^{\alpha})(u_y+\frac{N-1}{y} u)\Big)w_y\,\dd y
\nonumber\\
&\quad
-\int^r_{r_1}\rho w\Phi_y\,\dd y.\nonumber
\end{align}
Then we multiply the above equality by $\r w$ to obtain
\begin{align*}
&(\r^2 u^2+\r p)w^2\nonumber\\
&=\v \r(\rho+\alpha\delta\rho^{\alpha})(u_r+\frac{N-1}{r} u) w^2
-\v\r w\int_{r_1}^r \frac{N-1}{y} u(\rho+\alpha\delta\rho^{\alpha})_y w\,\dd y
-\r w\Big(\int_{r_1}^r\r u w\,\dd y\Big)_t \nonumber\\
&\quad-\r w\int_{r_1}^r \frac{N-1}{y} \r u^2 w\,\dd y
+\r w\int_{r_1}^r\Big(\r u^2+p-\v(\rho+\alpha\delta\rho^{\alpha})(u_y+\frac{N-1}{y} u)\Big)w_y\,\dd y\nonumber\\
&\quad-\rho w\int^r_{r_1}\rho w\Phi_y\,\dd y.
\end{align*}
Notice that
\begin{align}
\r w\Big(\int_{r_1}^r\r u w\,\dd y\Big)_t
&=-\r^2 u^2 w^2+\Big(\r w\int_{r_1}^r\r u w \dd y\Big)_t
+\Big(\r  uw\int_{r_1}^r\r u w\,\dd y\Big)_r\nonumber\\
&\quad-\r uw_r \int_{r_1}^r \r u w \,\dd y
+\frac{N-1}{r} \r u w \int_{r_1}^r \r u w\,\dd y,\nonumber
\end{align}
which yields
\begin{align}\label{7.34}
\r p w^2&=\v \r(\rho+\alpha\delta\rho^{\alpha})(u_r+\frac{N-1}{r} u) w^2
-\v\r w\int_{r_1}^r \frac{N-1}{y} u(\rho+\alpha\delta\rho^{\alpha})_y w\,\dd y
-\Big(\r w\int_{r_1}^r\r u w\,\dd y\Big)_t\nonumber\\
&\quad-\Big(\r  uw\int_{r_1}^r\r u w\,\dd y\Big)_r
+\r uw_r \int_{r_1}^r \r u w\,\dd y
-\frac{N-1}{r} \r u w \int_{r_1}^r \r u w\,\dd y\nonumber\\
&\quad-\r w\int_{r_1}^r \frac{N-1}{y} \r u^2 w\,\dd y
+\r w\int_{r_1}^r\Big(\r u^2+p-\v(\rho+\alpha\delta\rho^{\alpha})
\big(u_y+\frac{N-1}{y} u\big)\Big)w_y\,\dd y\nonumber\\
&\quad-\rho w\int^r_{r_1}\rho w\Phi_y\,\dd y:=\sum^{9}_{i=1}J_{i}.
\end{align}
To  estimate RHS of \eqref{7.34}, it follows first from Lemma \ref{bee} that
\begin{align}
&\int_{r_1}^{r_2}\r^\g\, r^{N-1}\dd r\leq C(r_2,\mathcal{E}_0),\label{7.35}\\
&\int_{r_1}^{r_2}\r\, \dd r\leq \f{C}{{r}^{N-1}_1}\int_{r_1}^{r_2}\r\,r^{N-1} \dd r
   \leq C\int_{r_1}^{r_2} \big(\rho^{\gamma}+1\big)\,r^{N-1}\dd r\leq C(r_1, r_2, \mathcal{E}_0),\label{7.36}\\
&\int_{r_1}^{r_2} \r u^2\, \dd r
 \leq \f{C}{{r}^{N-1}_1}\int_{r_1}^{r_2} \r u^2\,r^{N-1} \dd r\leq C(r_1,\mathcal{E}_0).\label{7.37}
\end{align}
We also obtain
\begin{align}
&\Big|\int^T_0\int^{r_2}_{r_1}\rho w\Big(\int^r_{r_1}\rho w\Phi_y\,\dd y\Big)
\dd r\dd t\Big|\nonumber\\
&\leq\Big|\int^T_0\int^{r_2}_{r_1}\rho w\Big(\int^r_{r_1}\rho w\big(\int^y_{\delta}(\rho-d(z))\,z^{N-1}\dd z\big)
\frac{\dd y}{y^{N-1}}\Big)\,\dd r\dd t\Big|\leq C(r_1,r_2,T,\mathcal{E}_0).
\end{align}
Now it follows from \eqref{7.36}--\eqref{7.37} that
\begin{align}\label{7.42}
\left| \int_0^T\int_{r_1}^{r_2}\big( J_3+\cdots+J_7\big)\, \dd r\dd t \right|\leq C(r_1,r_2,T, \mathcal{E}_0).
\end{align}
Next, we estimate $J_8$. Note that
\begin{align}\nonumber
\Big|\int_0^T\int_{r_1}^{r_2} \r w \Big(\int_{r_1}^r (\r u^2+p) w_y\, \dd y\Big)
\dd r \dd t\Big| \leq C(r_1,r_2,T,\mathcal{E}_0),
\end{align}
and
\begin{align}
&\v\Big|\int_0^T\int_{r_1}^{r_2} \r w \Big( \int_{r_1}^r (\rho+\alpha\delta\rho^{\alpha})(u_y+\frac{N-1}{y}u) w_y\dd y\Big) \dd r \dd t \Big|\nonumber\\
&\leq C(r_1,r_2,\mathcal{E}_0) \v\int_0^T \int_{r_1}^{r_2} \f{1}{y}
\big|(\rho+\alpha\delta\rho^{\alpha})\big(yu_y+(N-1)u\big) w_y\big| \dd y \dd t \nonumber\\
&\leq C(r_1,r_2,\mathcal{E}_0) \left(\v\int_0^T \int_{r_1}^{r_2} (\rho+\alpha\delta\rho^{\alpha})
 \big(yu_y+(N-1)u\big)^2\,\dd y\dd t
 +\int_0^T \int_{r_1}^{r_2}\v(\rho+\alpha\delta\rho^{\alpha})\,\dd y \dd t \right) \nonumber\\
&\leq C(r_1,r_2,T,\mathcal{E}_0),\nonumber
\end{align}
which yield
\begin{align}\label{7.43}
\Big|\int_0^T\int_{r_1}^{r_2} J_8\,\dd r \dd t\Big| \leq C(r_1,r_2,T,\mathcal{E}_0).
\end{align}
For $J_2$, since
\begin{align}
&\Big|\int_{r_1}^r \frac{N-1}{y} u(\rho+\delta\rho^{\alpha})_y w \dd y\Big|
=\Big|\frac{N-1}{r} \big((\rho+\delta\rho^{\alpha})u w\big)(t,r)\Big|\nonumber\\
&\,\,+(N-1)\Big|\int_{r_1}^r \Big(-(\rho+\delta\rho^{\alpha})uw
  +(\rho+\delta\rho^{\alpha}) u_y w+ (\rho+\delta\rho^{\alpha})u w_y\Big)(t,y)\,\frac{\dd y}{y}\Big|\nonumber\\
&\leq \Big|\frac{N-1}{r} \big((\rho+\delta\rho^{\alpha})u w\big)(t,r)\Big|
  +C(r_1)\int_{r_1}^{r_2} r^{N-1}(\rho+\delta\rho^{\alpha}) u_r^2 \,\dd r
   +C(r_1,r_2,\mathcal{E}_0),\nonumber
\end{align}
we have
\begin{align}\label{7.44}
\Big|\int_0^T\int_{r_1}^{r_2}J_2\,\dd r\dd t\Big|
&=\Big|\int_0^T\int_{r_1}^{r_2}\v\r w
  \Big(\int_{r_1}^r \frac{N-1}{y} u(\rho+\delta\rho^{\alpha})_y w\,\dd y\Big)
  \dd r\dd t\Big|\nonumber\\
&\leq C(r_1,r_2,T,\mathcal{E}_0)
+C(r_1,r_2,T,\mathcal{E}_0)\int_0^T\int_{r_1}^{r_2}\v r^2(\rho+\delta\rho^{\alpha}) u_r^2\,\dd r \dd t\nonumber\\
&\quad+\v\int_0^T\int_{r_1}^{r_2} \rho^2(\rho+\delta \rho^{\alpha})w^2\,\dd r\dd t\nonumber\\
&\leq C(r_1,r_2,T,\mathcal{E}_0)+\v\int_0^T\int_{r_1}^{r_2} \rho^3w^2\,\dd r\dd t,
\end{align}
since $\alpha\in(0,1).$

Finally, we estimate $J_1$. It is clear that
\begin{align}\label{7.46}
&\Big|\int_0^T\int_{r_1}^{r_2} J_1 \dd r \dd t\Big|
= \v\Big|\int_0^T\int_{r_1}^{r_2} (\r^2+\alpha\d\r^{\alpha+1})
  \big(u_r+\frac{N-1}{r} u\big) w^2  \dd r \dd t\Big|\nonumber\\
&\leq \v\int_0^T\int_{r_1}^{r_2} \r^{2}\big|u_r+\frac{N-1}{r} u\big|w^2\, \dd r \dd t
   +\alpha\d\v\int_0^T\int_{r_1}^{r_2}  \r^{\alpha+1} \big|u_r+\frac{N-1}{r} u\big| w^2 \,\dd r \dd t\nonumber\\
&\leq \v\int_0^T\int_{r_1}^{r_2} \r^{2}(\rho+\alpha\delta\rho^{\alpha}) w^2\,\dd r \dd t
  +C \v\int_0^T\int_{r_1}^{r_2}  (\r+\alpha\delta\rho^{\alpha})\big(u_r+\frac{N-1}{r} u\big)^2 w^2\,\dd r \dd t\nonumber\\
&\leq C(r_1,r_2,T,\mathcal{E}_0)+\v\int_0^T\int_{r_1}^{r_2} \rho^{3} w^2\,\dd r \dd t.
\end{align}
To close the estimate, we still need to bound the last term on RHS
of \eqref{7.44}--\eqref{7.46}.
We first consider the case: $\g\in(1,2]$. Then we have
\begin{align}\label{7.47}
\v\int_0^T\int_{r_1}^{r_2}  \rho^3  w^2\,\dd r \dd t
&\leq \v \int_0^T \int_{r_1}^{r_2} \r^{\g}\,\dd r
 \sup_{r\in[r_1,r_2]}\big(\r^{3-\g} w^2\big)\, \dd t\nonumber\\
&\leq C(r_1,r_2,\mathcal{E}_0)\int_0^T
\v\sup_{r\in[r_1,r_2]}\big(\r^{3-\g} w^2\big)\, \dd t\nonumber\\
&\leq C(r_1,r_2,\mathcal{E}_0)\int_0^T\v\int_{r_1}^{r_2}
\big|\big(\r^{3-\g} w^2\big)_r(t,r)\big|\,\dd r\dd t\nonumber\\
&\leq C(r_1,r_2,\mathcal{E}_0)\int_0^T\int_{r_1}^{r_2}
\big(\v \r^{2-\g}|\r_r| w^2+\v \r^{3-\g} w |w_r|\big)\,\dd r\dd t.
\end{align}
We now estimate each term of RHS of \eqref{7.47}.
A direct calculation shows that
\begin{align}\label{7.48}
& C(r_1,r_2,\mathcal{E}_0)
\v\int_0^T\int_{r_1}^{r_2} \r^{2-\g}|\r_r| w^2\,\dd r \dd t\nonumber\\
&\,\,\leq C(r_1,r_2,\mathcal{E}_0)
\int_0^T\int_{r_1}^{r_2} \v\r^{\g-2} \r_r^2\,\dd r \dd t
  +\f\v2\int_{0}^T\int_{r_1}^{r_2} \r^{3(2-\g)} w^2 \,\dd r \dd t\nonumber\\
&\,\,\leq C(r_1,r_2,\mathcal{E}_0)+\f\v2\int_{0}^T\int_{r_1}^{r_2}\r^{3} w^2\, \dd r \dd t,\nonumber\\[1mm]
\end{align}
and
\begin{align}\label{7.49}
& C(r_1,r_2,\mathcal{E}_0)\int_0^T\int_{r_1}^{r_2} \v \r^{3-\g} w |w_r|\,\dd r \dd t\nonumber\\
&\,\,\leq C(r_1,r_2,\mathcal{E}_0)\int_0^T\v \sup_{r}\,(\r w)(t,r) \int_{r_1}^{r_2}\r^{2-\g}|w _r|\,\dd r \dd t\nonumber\\
&\leq C(r_1,r_2,\mathcal{E}_0)\int_0^T\v \sup_{r}\,(\r w)(t,r) \, \dd t
\leq C(r_1,r_2,\mathcal{E}_0)\int_0^T\v \int_{r_1}^{r_2}\big( |\r_r| w+ \r |w_r|\big)\,\dd r \dd t\nonumber\\
&\,\,\leq C(r_1,r_2,\mathcal{E}_0)\int_0^T\int_{r_1}^{r_2} \v\r^{\g-2} \r_r^2\,\dd r \dd t
+C(r_1,r_2,\mathcal{E}_0)\int_0^T\int_{r_1}^{r_2} \r^{2-\g} w\, \dd r \dd t\nonumber\\
&\,\,\leq C(r_1,r_2,\mathcal{E}_0).\nonumber\\[1mm]
\end{align}

Combining \eqref{7.48}--\eqref{7.49}, we obtain
\begin{align}\label{7.50}
\v\int_0^T\int_{r_1}^{r_2}  \rho^3 w^2\,\dd r \dd t
\leq C(r_1,r_2,\mathcal{E}_0)\qquad\mbox{for $\g\in(1,2]$}.
\end{align}

For the case: $\g\in[2,3]$, we have
\begin{align}\label{7.51a}
\v \int_0^T\int_{r_1}^{r_2}  \rho^3  w^2\, \dd r \dd t
&\leq \v \int_0^T\sup_{r\in[r_1,r_2]}\,(\r^2 w)\int_{r_1}^{r_2} \r w\,\dd r \dd t\nonumber\\
&\leq C(r_1,r_2,\mathcal{E}_0) \v\int_0^T \int_{r_1}^{r_2}\big(\r|\r_r| w+\r^2 |w_r|\big)\,\dd r \dd t\nonumber\\
&\leq C(r_1,r_2,\mathcal{E}_0) \int_0^T \int_{r_1}^{r_2}\big(\v^2\r^{\g-2}|\r_r|^2 w+\r^2 |w_r|+\r^{4-\g} w\big)\,
\dd r \dd t\nonumber\\
&\leq C(r_1,r_2,\mathcal{E}_0).
\end{align}
For the case: $\gamma\in(3,\infty),$ we obtain
\begin{align}\label{7.511}
\v\int^T_0\int^{r_2}_{r_1}\rho^3w^2\,\dd r\dd t\leq C \int^T_0\int^{r_2}_{r_1}(\rho+\rho^{\gamma})\,\dd r\dd t\leq C(r_1,r_2,\mathcal{E}_0,T).
\end{align}
Now substituting \eqref{7.50}--\eqref{7.511} into \eqref{7.46}, we have
\begin{align}\label{7.52}
\Big|\int_0^T\int_{r_1}^{r_2}\big(J_1+J_2\big)\, \dd r \dd t\Big|\leq C(r_1,r_2,T,\mathcal{E}_0).
\end{align}
Integrating \eqref{7.34} over $[0,T]\times [r_1,r_2]$ and then using \eqref{7.42}--\eqref{7.43} and \eqref{7.52},
we conclude \eqref{6.14}.
$\hfill\Box$

\section{Uniform Higher Integrability of the Approximate Solutions}

To employ the compensated compactness framework \cite{Chen2010},
we need additional uniform integrability of the velocity for the approximate solutions.
We first define
\begin{equation}\label{M1}
\mathcal{M}_1:=\mathcal{E}_0+\rho_{\ast}+\rho^{-1}_{\ast}+\delta^{-1}+\varepsilon^{-1}+\sup_{b\geq1+\delta^{-1}}\mathcal{E}^{\varepsilon,\delta,b}_2<\infty,\quad \mathcal{M}_2:=\mathcal{M}_1+\sup_{b\geq1+\delta^{-1}}\mathcal{E}^{\varepsilon,\delta,b}_3,
\end{equation}
where
\begin{align}
\mathcal{E}_2^{\v,\d,b}:=\int_\d^b  \rho_0 \big(u_0^{8} +
\big|\frac{\mu_{0r}}{\rho_0}\big|^{8}\big)\,r^{N-1}\dd r,
\quad \mathcal{E}^{\v,\delta, b}_3:=\int_\d^b  \big(\f12 \rho_0 u_0^2+e(\rho_0,\r_{\ast})+\frac{1}{2}|\Phi_{0r}|^2\big)\,r^{2(N-1)+\vartheta}\dd r, \label{4.1d}
\end{align}
for some $\vartheta\in(0,1)$.
It follows from Lemma \ref{lem8.4} that
$\mathcal{E}_2^{\v,\d,b}$ and $\mathcal{E}_3^{\v,\delta, b}$ are uniformly bounded with respect to $b$. However,  the upper bounds
may depend on $(\v,\delta)$,
such that $\mathcal{M}_1$ and $\mathcal{M}_2$ are finite for any fixed $(\v,\delta)$, independent of $b>0$.

\begin{proposition}\label{lem6.4}
Given $r_1$ and $r_2$ with $\delta<r_1<r_2<b$, the smooth solution of \eqref{6.1}--\eqref{6.3} satisfies
that there exists $\vartheta\in(0,1)$ such that
	\begin{align}\label{6.15-1}
	&\int_0^T\int_{r_1}^{r_2}  \big(\rho|u|^3+\rho^{\g+\t}\big)(t,r)r^{N-1}\dd r\dd t
 \leq C(r_1,r_2,T,\mathcal{E}_0)+C(T,\mathcal{M}_2)\, b^{-\frac{\vartheta}{2}}.
	\end{align}
\end{proposition}

The proof of Proposition \ref{lem6.4} will be given later, after several lemmas are verified. Motivated by \cite{Chen2020}, the key point in Proposition \ref{lem6.4} is that the positive constant $C(T,\mathcal{M}_2)$ is independent of $b$,
so that this term vanishes when $b\rightarrow\infty.$
Note that CNSPEs with spherical symmetry is singular at $r=0$ and $r=\infty$.
In order to integrate the quantities in $r$ from $\infty$, we need to know the asymptotic behavior
of $\rho(t,r)$ near boundary $r=b$.
Thus, we have to obtain the lower and upper bound of $\rho$,
which should be independent of $b$.
First, we have the following lemma for the upper bound of $\rho$ (also see \cite[Lemma 4.2]{Chen2020}).

\begin{lemma}[\bf Upper Bound of the Density]\label{lem6.5}
The smooth solution of \eqref{6.1}--\eqref{6.3} satisfies
\begin{align}\label{6.16}
0< \rho(t,r)\leq C(\mathcal{M}_1)\qquad\mbox{for $t\geq0$ and $r\in[\delta,b]$},
\end{align}
\end{lemma}

\smallskip
We are now going to prove Lemma \ref{lem6.7} as a preparation to estimate the lower bound of the density.
\begin{lemma}\label{lem6.7}
The smooth solution of \eqref{6.1}--\eqref{6.3} satisfies
\begin{align}\label{6.28}
\int_{\delta}^b\rho |\rho^{-1-\frac{1}{2N}}\rho_r|^{8}\, r^{N-1}\dd r
\leq C(T,\mathcal{M}_1)\qquad \mbox{for $t\in [0,T]$}.
\end{align}
\end{lemma}

\noindent{\bf Proof.}
We write \eqref{6.93} as
\begin{equation}\label{7.6666}
\varepsilon(r^{N-1}\mu_x)_{\tau}+r^{N-1}p_x
=-u_{\tau}+\frac{1}{r^{N-1}}\Big(x-\int^{r}_{\delta}d(y)\,y^{N-1}\dd y\Big),
\end{equation}
and then integrate \eqref{7.6666} over $[0,\tau]$ to obtain
\begin{equation}\label{7.61}
\begin{split}
\varepsilon(r^{N-1}\mu_x)(\tau,x)&=\varepsilon(r^{N-1}\mu_x)(0,x)-\big(u(\tau,x)-u_0(x)\big)\\
&\quad-\int^{\tau}_0(r^{N-1}p_x)(s,x)\,\dd s
+\int^{\tau}_0\frac{1}{r^{N-1}}\Big(x-\int^r_{\delta}d(y)\,y^{N-1}\dd y\Big)\dd s.
\end{split}
\end{equation}
Multiplying \eqref{7.61} by $(r^{N-1}\mu_x)^{2k-1}$ and integrating the resultant equation,
we have
\begin{align}\label{7.622}
\varepsilon\int^{L_b}_0|r^{N-1}\mu_x|^{2k}\,\dd x
&\leq\Big(\int^{L_b}_0|r^{N-1}\mu_x|^{2k}\dd x\Big)^{\frac{2k-1}{2k}}\nonumber\\
&\quad \times \Big\{\varepsilon\Big(\int^{L_b}_0|(r^{N-1}\mu_x)(0,x)|^{2k}\,\dd x\Big)^{\frac{1}{2k}}+\|(u(\tau),u_0)\|_{L^{2k}}\nonumber\\
&\qquad\quad+C_T\Big(\int^T_0\int^{L_b}_{0}|r^{N-1}(\rho^{\gamma})_x|^{2k}\,\dd x\dd s\Big)^{\frac{1}{2k}}\nonumber\\
&\qquad\quad+C_T\Big(\int^T_0\int^{L_b}_0\Big|\frac{1}{r^{N-1}}
\int^r_{\delta}(\rho(s,y)-d(y))\dd y\Big|^{2k}\dd x\dd s\Big)^{\frac{1}{2k}}\Big\}.
\end{align}
Since
$$
|\mu_x|=|\rho_x+\delta(\rho^{\alpha})_x|=|(\frac{1}{\alpha}\rho^{1-\alpha}+\delta)(\rho^{\alpha})_x|\geq\delta(\rho^{\alpha})_x,
$$
and $(\rho^{\gamma})_x=\frac{\gamma}{\alpha}\rho^{\gamma-\alpha}(\rho^{\alpha})_x,$
it follows from \eqref{6.16} and \eqref{7.622}  that
\begin{align}\label{7.63}
&\int^{L_b}_0|r^{N-1}(\rho^{\alpha})_x|^{2k}\,\dd x\nonumber\\
&\leq C(\mathcal{E}_0,T,\varepsilon,\delta)\Big(\int^{L_b}_0
   \big(|(r^{N-1}\mu_x)(0,x)|^{2k}+|u(\tau)|^{2k}+|u_0|^{2k}\big)\dd x
   +\int^{\tau}_0\int^{L_b}_0|r^{N-1}(\rho^{\alpha})_x|^{2k}\,\dd x\dd s\nonumber\\
&\qquad\qquad\qquad\quad\,\,
+\int^{\tau}_0\int^{L_b}_0\Big|\frac{1}{r^{N-1}}\int^r_{\delta}(\rho(s,y)-d(y))\,y^{N-1}\dd y\Big|^{2k}  \dd x\dd s\Big).
\end{align}
Pulling \eqref{7.63} back to the Eulerian coordinates, we see that
\begin{align}\label{7.64}
&\int^{b}_{\delta}\rho|\rho^{-1-\frac{1}{2N}}\rho_r|^{2k}\,r^{N-1}\dd r\nonumber\\
&\leq C(\mathcal{E}_0,T,\varepsilon,\delta)\Big(\mathcal{E}^{\varepsilon,\delta,b}_2
  +\int^{b}_{\delta}\rho |u|^{2k}\,r^{N-1}\dd r+\int^t_0\int^{b}_{\delta}\rho|\rho^{-1-\frac{1}{2N}}\rho_r|^{2k}\,r^{N-1}\dd r\dd s\nonumber\\
&\qquad\qquad\qquad\quad\,\,\,+\int^t_0\int^b_{\delta}\Big|\frac{1}{r^{N-1}}
  \int^r_{\delta}\big(\rho(s,y)-d(y)\big)\,y^{N-1}\dd y\Big|^{2k}\rho\, r^{N-1}\dd r\dd s\Big).
\end{align}
Now, we consider the last term of \eqref{7.64}. Notice that
\begin{align}\label{7.65}
&\frac{1}{r^{N-1}}\Big|\int^r_{\delta}\big(\rho(s,y)-d(y)\big)\,y^{N-1}\dd y\Big|\nonumber\\
&\leq\frac{1}{r^{N-1}}\int^r_{\delta}|\rho(s,y)-\rho_{\ast}|\,y^{N-1}\dd y
  +\frac{1}{r^{N-1}}\int^r_{\delta}\big(\rho_{\ast}-d(y)\big)\,y^{N-1}\dd y\nonumber\\
&\leq \|d-\rho_{\ast}\|_{L^1}r^{-N+1}
+\frac{1}{r^{N-1}}\int^r_{\delta}|\rho(s,y)-\rho_{\ast}|\mathbf{1}_{A(t)\cup B(t)}\,y^{N-1}\dd y\nonumber\\
&\quad+\frac{1}{r^{N-1}}\int^r_{\delta}|\rho(s,y)-\rho_{\ast}|
 \big(1-\mathbf{1}_{A(t)\cup B(t)}\big)\,y^{N-1}\dd y.
\end{align}
Using \eqref{eee}, we have
\begin{align}\label{7.66}
&\frac{1}{r^{N-1}}\int^r_{\delta}|\rho(s,y)-\rho_{\ast}|\mathbf{1}_{A(t)\cup B(t)}\,y^{N-1}\dd y\nonumber\\
&\,\,\leq\frac{C(\rho_{\ast})}{r^{N-1}}\int^r_{\delta}\mathbf{1}_{B(t)}\,y^{N-1}\dd y
  +\frac{C}{r^{N-1}}\int^r_{\delta}\rho \mathbf{1}_{A(t)}\,y^{N-1}\dd y\nonumber\\
&\,\,\leq\frac{C(\rho_{\ast})}{r^{N-1}}\int^r_{\delta}e(\rho,\rho_{\ast})\,y^{N-1}\dd y
\leq C(\rho_{\ast})\mathcal{E}^{\varepsilon,\delta,b}_0r^{-N+1},\nonumber\\[1mm]
\end{align}
and
\begin{align}
&\frac{1}{r^{N-1}}\int^r_{\delta}|\rho(s,y)-\rho_{\ast}|(1-\mathbf{1}_{A(t)\cup B(t)})\,y^{N-1}\dd y\nonumber\\
&\leq\frac{1}{r^{N-1}}\Big(\int^r_{\delta}|\rho(s,y)-\rho_{\ast}|^2\big(1-\mathbf{1}_{A(t)\cup B(t)}\big)\,y^{N-1}\dd y\Big)^{\frac{1}{2}}\Big(\int^r_{\delta}y^{N-1}\dd y\Big)^{\frac{1}{2}}\nonumber\\
&\leq C(\rho_{\ast})r^{-\frac{N}{2}+1}\Big(\int^r_{\delta}e(\rho,\rho_{\ast})\,y^{N-1}\dd y\Big)^{\frac{1}{2}}
\leq C(\rho_{\ast},\mathcal{E}^{\varepsilon,\delta,b}_0)r^{-\frac{N}{2}+1}.\label{7.77}
\end{align}
Substituting \eqref{7.66}--\eqref{7.77} into \eqref{7.65}, we conclude
\begin{equation}\label{7.78}
\begin{split}
&\frac{1}{r^{N-1}}\Big|\int^r_{\delta}\big(\rho(s,y)-d(y)\big)\,y^{N-1}\dd y\Big|
\leq C(\rho_{\ast},\mathcal{E}_0)(r^{-N+1}+r^{-\frac{N}{2}+1})\leq C(\rho_{\ast}, \mathcal{E}_0,\delta)r^{-\frac{N}{2}+1}.
\end{split}
\end{equation}
Therefore, using \eqref{7.78}, we obtain
\begin{align}\label{7.79}
&\int^t_0\int^b_{\delta}\rho \Big|\frac{1}{r^{N-1}}
 \int^r_{\delta}\big(\rho(s,y)-d(y)\big)\,y^{N-1}\dd y\Big|^{2k}\,r^{N-1}\dd r\dd s\nonumber\\
&\leq C(\rho_{\ast},\mathcal{E}_0,\delta)\int^t_0\int^{b}_{\delta}\rho r^{-k(N-2)}\,r^{N-1}\dd r\dd s\nonumber\\
&=C(\rho_{\ast},\mathcal{E}_0,\delta)
\Big(\int^t_0\int^b_{\delta}2\rho_{\ast}\mathbf{1}_{\{\rho(t,r)\leq2\rho_{\ast}\}}r^{N-1-k(N-2)}\dd r\dd s\nonumber\\
&\quad\qquad\qquad\qquad
 +\int^t_0\int^b_{\delta}e(\rho,\rho_{\ast})r^{N-1-k(N-2)}\mathbf{1}_{\{\rho(t,r)>2\rho_{\ast}\}}\dd r\dd s\Big)\nonumber\\
&\leq C(\rho_{\ast},\mathcal{E}_0,T,\delta)\Big(1+\int^b_{\delta}r^{N-1-k(N-2)}\dd r\Big)\nonumber\\
&\leq C(\rho_{\ast},\mathcal{E}_0,T,\delta)\Big(1+\frac{1}{k(N-2)-N}\big(\delta^{N-k(N-2)}-b^{N-k(N-2)}\big)\Big)\nonumber\\
&\leq C(\rho_{\ast},\mathcal{E}_0,T,\delta),
\end{align}
provide that $k>\frac{N}{N-2}$, which can be satified by taking $k=4.$

In order to bound the fourth term, we multiply $\eqref{6.1}_2$ by $r^{N-1}u^{2k-1}$ to obtain
\begin{align}
&\frac{\dd }{\dd t}\int^b_{\delta}\frac{1}{2k}r^{N-1}\rho u^{2k}\,\dd r
  -\int^{b}_{\delta}p(r^{N-1}u^{2k-1})_r\,\dd r\nonumber\\
&=-\varepsilon\int^b_{\delta}\Big((\mu+\lambda)(u_r+\frac{N-1}{r}u)(r^{N-1}u^{2k-1})_r
  -(N-1)\mu(r^{N-2}u^{2k})_r\Big)\,\dd r\nonumber\\
&\quad+\int^b_{\delta}\rho u^{2k-1}
\Big(\int^r_{\delta}\big(\rho(t,y)-d(y)\big)\,y^{N-1}\dd y\Big)\dd r.\label{7.80}
\end{align}
For the last term of RHS of \eqref{7.80}, it follows from \eqref{7.78} that
\begin{align}\label{7.81}
&\int^b_{\delta}\rho u^{2k-1}\Big(\int^r_{\delta}
\big(\rho(t,y)-d(y)\big)\,y^{N-1}\dd y\Big)\dd r\nonumber\\
&=C(\rho_{\ast},\mathcal{E}_0)\int^b_{\delta}\rho u^{2k-1}r^{\frac{N}{2}}\dd r\nonumber\\
&\leq C(\rho_{\ast},\mathcal{E}_0)\Big(\int^b_{\delta}\rho u^{2k}\,r^{N-1}\dd r\Big)^{\frac{2k-1}{2k}}\Big(\int^b_{\delta}\rho r^{-N(k-1)+2k-1}\dd r\Big)^{\frac{1}{2k}}\nonumber\\
&\leq C(\rho_{\ast},\mathcal{E}_0,\varepsilon,\delta)
  \Big(\int^b_{\delta}\rho u^{2k}\,r^{N-1}\dd r\Big)^{\frac{2k-1}{2k}}
   \Big(\int^b_{\delta} r^{-N(k-1)+2k-1}\dd r\Big)^{\frac{1}{2k}}\nonumber\\
&\leq C(\rho_{\ast},\mathcal{E}_0,\varepsilon,\delta)\Big(\int^b_{\delta}\rho u^{2k}\,r^{N-1}\dd r\Big)^{\frac{2k-1}{2k}},
\end{align}
provided that $N(k-1)>2k,$ {\it i.e.}, $k>\frac{N}{N-2}.$
A direct calculation shows that
\begin{align}\label{7.82}
&(\mu+\lambda)(u_r+\frac{N-1}{r}u)(r^{N-1}u^{2k-1})_r-(N-1)\mu(r^{N-2}u^{2k})_r\nonumber\\
&=\rho\big((2k-1)r^{N-1}u^{2k-2}u^2_r+(N-1)r^{N-3}u^{2k}\big)\nonumber\\
&\quad+\delta\Big(\alpha(2k-1)r^{N-1}u^{2k-2}u^2_r+2k(N-1)(\alpha-1)r^{N-2}u^{2k-1}u_r\nonumber\\
&\qquad\quad\,+(N-1)\big(\alpha(N-1)-(N-2)\big)r^{N-3}u^{2k}\Big).
\end{align}
Since $\alpha=1-\frac{1}{2N},$ we calculate the discriminant $\triangle$ for the last term of \eqref{7.82} to have
\begin{align*}
\triangle&=4k^2(N-1)^2(\alpha-1)^2-4\alpha(2k-1)(N-1)\big(\alpha(N-1)-(N-2)\big)\nonumber\\
&=\frac{N-1}{N^2}\big(k^2(N-1)-(2k-1)(2N-1)(N+1)\big).
\end{align*}
Choosing $k=4,$ we have
\begin{equation*}\label{7.84}
\triangle=-\frac{N-1}{N^2}(14N^2-9N+9)<0,
\end{equation*}
which, together with \eqref{7.82} yields
\begin{align}\label{7.85}
&(\mu+\lambda)(u_r+\frac{N-1}{r}u)(r^{N-1}u^{2k-1})_r-(N-1)\mu(r^{N-2}u^{2k})_r\nonumber\\
&\geq \rho\big((2k-1)r^{N-1}u^{2k-2}u^2_r+(N-1)r^{N-3}u^{2k}\big).
\end{align}
For the pressure term, it follows from \eqref{6.16} that
\begin{eqnarray}
&&\Big|\int^b_{\delta}p(r^{N-1}u^{2k-1})_r\,\dd r\Big|\nonumber\\
&&=\Big|\int^b_{\delta}p\big((2k-1)r^{N-1}u^{2k-2}u_r+(N-1)r^{N-2}u^{2k-1}\big)\,\dd r\Big|\nonumber\\
&&\leq\frac{\varepsilon}{2}\int^b_{\delta}\rho\big(r^{N-1}u^{2k-2}u^2_r+r^{N-3}u^{2k}\big)\,\dd r
+ C_{\v}\int^b_{\delta}\rho^{2\gamma-1}u^{2k-2}\,r^{N-1}\dd r\nonumber\\
&&\leq \frac{\varepsilon}{2}\int^b_{\delta}\rho\big(r^{N-1}u^{2k-2}u^2_r+r^{N-3}u^{2k}\big)\,\dd r
  +C(T,\mathcal{M}_1)\int^b_{\delta}\big(\rho u^2+\rho u^{2k}\big)\,r^{N-1}\dd r\nonumber\\
&&\leq \frac{\varepsilon}{2}\int^b_{\delta}\rho\big(r^{N-1}u^{2k-2}u^2_r+r^{N-3}u^{2k}\big)\,\dd r
 +C(T,\mathcal{M}_1)\Big(1+\int^b_{\delta}\rho u^{2k}\,r^{N-1}\dd r\Big).\label{7.86}
\end{eqnarray}
Combining \eqref{7.65}, \eqref{7.81}, and \eqref{7.86} with \eqref{7.80}, we conclude
$$
\frac{\dd }{\dd t}\int^{b}_{\delta}\frac{1}{2k}\rho u^{2k}\,r^{N-1}\dd r
\leq C(T,\mathcal{M}_1)\Big(1+\int^b_{\delta}\rho u^{2k}\,r^{N-1}\dd r\Big),
$$
which, together with Gr\"{o}nwall's inequality, implies
\begin{equation}\label{7.87}
\int^{b}_{\delta}\rho u^{2k}\,r^{N-1}\dd r\leq C(T,\mathcal{M}_1)\qquad \mbox{for all $t\in[0,T]$}.
\end{equation}
Substituting  \eqref{7.79} and \eqref{7.87} into \eqref{7.64}, we have
\begin{equation*}\label{7.88}
\int^b_{\delta}\rho|\rho^{-1-\frac{1}{2N}}\rho_r|^{2k}\,r^{N-1}\dd r
\leq C(T,\mathcal{M}_1)\Big(1+\int^T_0\int^b_{\delta}\rho|\rho^{-1-\frac{1}{2N}}\rho_r|^{2k}\,r^{N-1}\dd r\dd s\Big).
\end{equation*}
Again, applying the Gr\"{o}nwall inequality, we obtain
\begin{equation}\label{7.89}
\int^b_{\delta}\rho|\rho^{-1-\frac{1}{2N}}\rho_r|^{2k}\, r^{N-1}\dd r\leq C(T,\mathcal{M}_1).
\end{equation}
Taking $k=4$, we conclude the proof.
$\hfill\Box$

%
%

\medskip
With Lemma \ref{lem6.7}, we obtain the following lower bound of the density.

\begin{lemma}[\bf Lower Bound of the Density]\label{lem6.8}
There exists $C(T,\mathcal{M}_1)>0$ depending only on $T$ and $\mathcal{M}_1$
such that
the smooth solution of \eqref{6.1}--\eqref{6.3} satisfies
\begin{align}\label{6.34}
\rho(t,r) \geq C(T,\mathcal{M}_1)^{-1}>0\qquad
\mbox{for $(t,r)\in [0,T]\times [\delta,b]$}.
\end{align}
\end{lemma}

\noindent{\bf Proof.}
Since $\rho(t,r)$ is a continuous function on $[\delta,b],$ then, for any $r\in B(t)$
with $\rho(t,r)<\frac{1}{2}\rho_{\ast},$ there exists $r_0\in B(t)$ such that
\begin{equation*}\label{2.1}
\rho(t,r_0)=\frac{\rho_{\ast}}{2}, \qquad\,\,\, |r-r_0|\leq C(\rho_{\ast},\mathcal{E}_0).
\end{equation*}
Therefore, for $0<\beta<\frac{k-N}{2kN}$ (for instance, $\beta=\frac{k-N}{4kN}$), we have
\begin{align*}
\rho(t,r)^{-\beta}&=\rho(t,r_0)^{-\beta}+\beta\int^r_{r_0}\rho^{-\beta-1}\rho_y\,\dd y\nonumber\\
&\leq C(\rho_{\ast})+\beta\Big(\int^r_{r_0}\rho|\rho^{-1-\frac{1}{2N}}\rho_y|^{2k}\,\dd y\Big)^{\frac{1}{2k}}
 \Big(\int^r_{r_0}\big(\rho^{-\beta+\frac{k-N}{2Nk}}\big)^{\frac{2k}{2k-1}}\dd y\Big)^{\frac{2k-1}{2k}}\nonumber
 \\[1mm]
&\leq C(\rho_{\ast})+\beta C(T,\mathcal{M}_1),\label{2.2}
\end{align*}
where we have used \eqref{6.16}, \eqref{6.28}, and \eqref{6.994}.
Thus, we obtain
\begin{equation*}\label{2.3}
\rho(t,r)\geq C(T,\mathcal{M}_1)^{-1} \qquad \mbox{for any $r\in B(t)$},
\end{equation*}
which yields \eqref{6.34}.
This completes the proof.
$\hfill\Box$

\medskip
Using the above lower and upper bounds of the density we obtained,
though its bounds depend on $\v^{-1}$ and $ \d^{-1},$ we can have the following useful weighted estimate.

\begin{lemma}\label{lem6.10}
For any constant $\vartheta\in(0,1)$, the smooth solution of \eqref{6.1}--\eqref{6.3}
satisfies
\begin{align}\label{6.42}
&\int_{\delta}^b \Big(\frac12\rho u^2+e(\rho,\r_{\ast})+\frac{1}{2}|\Phi_r|^2\Big)(t,r)r ^{2(N-1)+\vartheta}
  \,\dd r\nonumber\\
&\quad+\v\int_0^T\int_{\delta}^b\Big( r^{2(N-1)+\vartheta} \rho u_r^2+\alpha\d r^{2(N-1)+\vartheta} \rho^{\alpha} u_r^2\Big) \dd r\dd t\leq C(T,\mathcal{M}_2).
\end{align}
\end{lemma}

\noindent{\bf Proof.}
Let $L>0.$ Multiplying $\eqref{6.1}_2$ by $r^{N-1+L}u$ and integrating by parts, we have
\begin{equation}\label{2.31}
\begin{split}
&\frac{\dd }{\dd t}\int^b_{\delta}\frac{1}{2}\rho u^2\,r^{N-1+L}\dd r
   +\int^b_{\delta}p_ru\,r^{N-1+L}\dd r\\
&=\frac{1}{2}\int^b_{\delta}\rho u^3\,r^{N-2+L}\dd r
  -\varepsilon\int^b_{\delta}(\mu+\lambda)\big(u_r+\frac{N-1}{r}u\big)\big(ru_r+(N-1+L)u\big)\,r^{N-2+L}\dd r\\
&\quad+\varepsilon(N-1)\int^b_{\delta}\mu\big(2ruu_r+(N-2+L)u^2\big)\,r^{N-3+L}  \dd r\\
&\quad+\int^b_{\delta}\rho u\Big(\int^r_{\delta}\big(\rho(t,z)-d(z)\big)\,z^{N-1}\dd z\Big)\, r^L\dd r.
\end{split}
\end{equation}
For the last term on RHS of \eqref{2.31}, we notice from \eqref{phib} that,
for any $r\in[\delta,b],$
\begin{equation*}\label{2.32}
r^{N-1}\Phi_{rt}=-\int^r_{\delta}(y^{N-1}\rho(t,y))_t\,\dd y
=\int^r_{\delta}(y^{N-1}\rho u)_y\,\dd y
=\rho ur^{N-1},
\end{equation*}
which yields
\begin{equation}\label{2.33}
\Phi_{rt}=\rho u\qquad \mbox{for any $t\in[\delta,b]$}.
\end{equation}
Using \eqref{2.33}, we have
\begin{align*}
&\int^b_{\delta}\rho u\Big(\int^r_{\delta}\big(\rho(t,y)-d(y)\big)\,y^{N-1}\dd y\Big)\,r^L\dd r\nonumber\\
&=-\int^b_{\delta}\Phi_{rt}\Phi_r\,r^{N-1+L}\dd r
=-\frac{1}{2}\frac{\dd}{\dd t}
\int^b_{\delta}|\Phi_r|^2\,r^{N-1+L}\dd r.
\end{align*}
The other terms can be dealt by similar arguments as in \cite{Chen2020}. Then we conclude
\begin{align}\label{2.35}
&\frac{\dd }{\dd t}\int^b_{\delta}\Big(\frac{1}{2}\rho u^2+e(\rho,\rho_{\ast})+\frac{1}{2}|
  \Phi_r|^2\Big)\,r^{N-1+L}\dd r
   +\frac{\varepsilon}{2}\int^b_{\delta}\big(\rho u^2_r+\alpha\delta \rho^{\alpha}u^2_r\big)\, r^{N-1+L}\dd r
     \nonumber\\
&\leq C(T,M_1)\bigg(\int^b_{\delta}\big(\rho u^2+e(\rho,\rho_{\ast})\big)\,r^{N-2+L}\dd r
+\Big(\int^b_{\delta}\rho u^2\,r^{N-2+L}\dd r\Big)^{\frac{4}{3}}+\int^b_{\delta}\rho u^2_r\,r^{N-1}\dd r\bigg).
\end{align}
Taking $L=1$ in \eqref{2.35} and integrating it over $[0,t]$, it follows from from \eqref{6.11} that
\begin{equation*}\label{2.36}
\begin{split}
&\int^b_{\delta}\Big(\frac{1}{2}\rho u^2+e(\rho,\rho_{\ast})+\frac{1}{2}|\Phi_r|^2\Big)\,r^{N}\dd r
+\frac{\varepsilon}{2}\int^t_0\int^{b}_{\delta}\big(\rho u^2_r+\alpha\delta \rho^{\alpha}u^2_r)\,r^N\dd r\dd s\\
&\leq\int^b_{\delta}\Big(\frac{1}{2}\rho_0u^2_0+e(\rho_0,\rho_{\ast})+\frac{1}{2}|\Phi_{0r}|^2\Big)\,r^N\dd r+C(T,\mathcal{M}_2)\leq C(T,\mathcal{M}_2).
\end{split}
\end{equation*}
Then,  taking $L=2,3, \cdots, N-1$ in \eqref{2.35} step by step, we have
\begin{align}\label{6.50}
&\int_\d^b \big(\f12 \rho u^2+e(\rho,\r_{\ast})+\frac{1}{2}|\Phi_r|^2\big)\,r^{2N-2}\dd r
  +\frac{\v}{2}\int_0^t\int_\d^b (\rho +\alpha\delta \rho^{\alpha}) u_r^2\, r^{2N-2}\dd r\dd s\nonumber\\
&\leq \int_\d^b \big(\f12 \rho_0 u_0^2+e(\rho_0,\r_{\ast})+\frac{1}{2}|\Phi_{0r}|^2\big)\,r^{2N-2}\dd r
  +C(T,\mathcal{M}_2)
\leq C(T,\mathcal{M}_2)
 \qquad\,\, \mbox{for any $t\in[0,T]$}.
\end{align}
Finally, taking $L=N-1+\vartheta$ with $\vartheta\in(0,1)$ in \eqref{2.35}
and integrating it over $[0,t]$, we conclude from \eqref{6.50} that
\begin{align*}
&\int_\d^b  \big(\f12 \rho u^2+e(\rho,\r_{\ast})+\frac{1}{2}|\Phi_{r}|^2\big)\,r^{2N-2+\vartheta}\dd r
  +\v\int_0^t\int_\d^b (\rho +\alpha\delta \rho^{\alpha}) u_r^2\, r^{2N-2+\vartheta}\dd r\dd s\nonumber\\
&\leq \int_\d^b \big(\f12 \rho_0 u_0^2+e(\rho_0,\r_{\ast})+\frac{1}{2}|\Phi_{0r}|^2\big)\,r^{2N-2+\vartheta}\dd r
   +C(T,\mathcal{M}_2)
\leq C(T,\mathcal{M}_2)\qquad \mbox{for any $t\in[0,T]$}.
\end{align*}
Then the proof is completed.
$\hfill\Box$

\medskip
Using Lemmas \ref{bee}, \ref{BD}, \ref{lem6.5}, \ref{lem6.8}, and \ref{lem6.10},
by similar arguments as in \cite{Chen2020}, we obtain
the following decay estimates.

\begin{lemma}[\cite{Chen2020}, Lemma 4.7]\label{lem6.11}
The smooth solution of \eqref{6.1}--\eqref{6.3} satisfies that, for $r\in[1,b]$,
\begin{align*}
&|(\rho-\r_{\ast})(t,r)|\leq C(T,\mathcal{M}_2)\,r^{-\frac{3}{4}N+\frac{3}{4}-\frac{\vartheta}{4}},\\
&\int_0^T \big(|u(t,r)|+ |u(t,r)|^3\big)\, r^{N-1}\dd t\leq C(T,\mathcal{M}_2)\,r^{-\frac{\vartheta}{2}}.
\end{align*}
\end{lemma}

\medskip
Taking $\psi(s)=\frac{1}{2}s|s|$ in \eqref{weakentropy}, then the corresponding entropy and entropy flux pair
can be represented as
\begin{align*}
\begin{cases}
\eta^{\#}(\r,\r u)=\f12 \r \int_{-1}^1 (u+\r^{\t} s) |u+\r^{\t}s| [1-s^2]_+^{\l}\,\dd s,\\[2mm]
q^{\#}(\r, \r u)=\f12 \r \int_{-1}^1 (u+\theta\r^{\t}s)(u+\r^{\t} s) |u+\r^{\t}s| [1-s^2]_+^{\l}\,\dd s,
\end{cases}
\end{align*}
where $\theta=\frac{\gamma-1}{2}$. Then it is direct to show that
\begin{align}\label{6.59}
|\eta^{\#}(\r,\r u)|\lesssim \r |u|^2+\r^{\g}, \qquad\, q^{\#}(\r,\r u)\gtrsim \r |u|^3+\r^{\g+\t}.
\end{align}
Since
\begin{align*}
\begin{cases}
\partial_\r\eta^{\#}=\int_{-1}^1 \big(-\f12 u+(\t+\f12)\r^{\t} s\big)\,|u+\r^{\t}s| [1-s^2]_+^{\l}\,\dd s,\\[1mm]
\partial_m\eta^{\#}=\int_{-1}^1|u+\r^{\t}s| [1-s^2]_+^{\l}\,\dd s,
\end{cases}
\end{align*}
then it is direct to see that, for some constant $C=C(\g)>0,$
\begin{align*}
&|\eta^{\#}_m|\leq C\big(|u|+\rho^\theta\big),\qquad\, |\eta^{\#}_\rho|\leq C\big(|u|^2+\rho^{2\theta}\big),\\
&\eta^{\#}_\rho(\rho,0)=0,\quad  \eta^{\#}_m(\rho,0)=2{\rho}^{\theta} \int_{0}^1 s[1-s^2]^{\lambda}\dd s.
\end{align*}
%
%
%
%
%

Define the relative-entropy pair:
\begin{align*}
\begin{cases}
\tilde{\eta}(\rho,\rho u)=\eta^{\#}(\r,\r u)-\eta^{\#}(\r_{\ast},0)- \eta^{\#}_m(\r_{\ast},0) \rho u,\\[1mm]
\tilde{q}(\rho,\rho u)=q^{\#}(\r,\r u)-q^{\#}(\r_{\ast},0)
  -\eta^{\#}_m(\r_{\ast},0)\big(\rho u^2+p(\rho)-p(\rho_{\ast})\big).
\end{cases}
\end{align*}

Now, we recall a key lemma established in \cite{Chen2020}, which will be used to prove Proposition \ref{lem6.4}.
\begin{lemma}[\cite{Chen2020}, Lemma 4.8]\label{lem6.14}
There exists $C_\gamma(\rho_{\ast})>0$ depending only on $\gamma$ such that
\begin{align}\label{6.78}
-\tilde{q}(\rho,\rho u)+\rho u  \partial_\rho \tilde{\eta}(\rho,\rho u)
+\rho u^2 \partial_m \tilde{\eta}(\rho,\rho u)
\leq C_\gamma(\rho_{\ast}) \big(\rho u^2+ e(\rho,\rho_{\ast})\big).
\end{align}
\end{lemma}

\medskip
Now, we are ready to prove  Proposition \ref{lem6.4}.

\medskip
\noindent{\bf Proof of Proposition \ref{lem6.4}:}
A direct calculation shows that
\begin{align}\label{3.1}
&(r^{N-1}\tilde{\eta})_t+(r^{N-1}\tilde{q})_r
+(N-1)r^{N-2}\big(-\tilde{q}+\rho u\partial_{\rho}\tilde{\eta}+\rho u^2\partial_m\tilde{\eta}\big)\nonumber\\
&=\varepsilon r^{N-1}\partial_m\tilde{\eta}\Big\{\big((\rho+\alpha\delta\rho^{\alpha})(u_r+\frac{N-1}{r}u)\big)_r
 -\frac{N-1}{r}(\rho+\delta\rho^{\alpha})_ru\Big\}+\partial_m\tilde{\eta}\rho r^{N-1}\Phi_r.
\end{align}

Let $y\in[b-1,b]$ and $r\in[r_1,r_2].$ Integrating \eqref{3.1} over $[r,y]$ and then over $[0,T]\times[b-1,b]\times[r_1,r_2],$ we have
\begin{align}\label{3.3}
&\int^T_0\int^{r_2}_{r_1}\tilde{q}(t,r)\,r^{N-1}\dd r\dd t\nonumber\\
&=(N-1)\int^T_0\int^b_{b-1}\int^{r_2}_{r_1}\int^y_r
\big(-\tilde{q}+\rho u\partial_{\rho}\tilde{\eta}+\rho u^2\partial_{m}\tilde{\eta}\big)
\,z^{N-2}\dd z\dd r\dd y\dd t\nonumber\\
&\quad+\int^b_{b-1}\int^{r_2}_{r_1}\int^y_r\big(\tilde{\eta}(t,z)-\tilde{\eta}(0,z)\big)\,z^{N-1}\dd z\dd r\dd y
+(r_2-r_1)\int^T_0\int^b_{b-1}\tilde{q}(t,y)\,y^{N-1}\dd y\dd t\nonumber\\
&\quad-\varepsilon\int^T_0\int^b_{b-1}\int^{r_2}_{r_1}\int^y_r
\partial_m\tilde{\eta}
\Big\{\big((\rho+\alpha\delta\rho^{\alpha})(u_z+\frac{N-1}{z}u)\big)_z
  -\frac{N-1}{z}(\rho+\delta\rho^{\alpha})_zu\Big\}\,z^{N-1}\dd z\dd r\dd y\dd t\nonumber\\
&\quad-\int^T_0\int^b_{b-1}\int^{r_2}_{r_1}\int^y_r\partial_m\tilde{\eta}\rho \Phi_z\,z^{N-1}\dd z\dd r\dd y\dd t:=\sum^{5}_{i=1}K_{i}.
\end{align}
It follows from \eqref{6.78} that
\begin{equation*}\label{3.4}
K_1\leq C(\rho_{\ast})\frac{r_2}{r_1}\int^T_0\int^b_{r_1}\big(\rho u^2+e(\rho,\rho_{\ast})\big)(t,y)\,y^{N-1}\dd y\dd t\leq C(\rho_{\ast})\frac{r_2}{r_1}(\mathcal{E}_0+1).
\end{equation*}
Similar to the argument in \cite[Proposition 4.1]{Chen2021}, we can obtain
$$
|\tilde{\eta}(\rho,\rho u)|\leq C_{\gamma}\big(\frac{m^2}{\rho}+e(\rho,\rho_{\ast})\big),
$$
and then
\begin{equation*}\label{3.5}
|K_2|\leq C_{\gamma}r_2(\mathcal{E}_0+1).
\end{equation*}
For the third term on RHS of \eqref{3.3}, since
\begin{equation*}\label{3.6}
|q^{\#}(\rho,\rho u)-q^{\#}(\rho,0)|\leq C_{\gamma}\big(\rho |u|^3+\rho^{1+2\theta}|u|\big),
\end{equation*}
we have
\begin{equation}\label{3.7}
\tilde{q}(\rho,\rho u)\leq C(T,\mathcal{M}_2)\big(|\rho-\rho_{\ast}|^2+|u|^3+|u|\big).
\end{equation}
It follows from \eqref{3.7} and Lemma \ref{lem6.11} that
\begin{equation}\label{3.8}
\begin{split}
|K_3|\leq C(T, \mathcal{M}_2)r_2\int^b_{b-1}\int^T_0
\big(|\rho-\rho_{\ast}|^2+|u|^3+|u|\big)(t,y)\,y^{N-1}\dd y\dd t
\leq C(T,\mathcal{M}_2)\,b^{-\frac{\vartheta}{2}}\\
\end{split}
\end{equation}
To bound the viscous term, we regard $\tilde{\eta}_m(\rho,\rho u)$ as a function of $(\rho, u)$,
then
\begin{equation}\label{3.9}
\begin{cases}
|\partial_m\tilde{\eta}(\rho,\rho u)|\leq C_{\gamma}\big(|u|+|\rho^{\theta}-\rho^{\theta}_{\ast}|\big),\\[1mm]
|\partial_{mu}\tilde{\eta}(\rho,\rho u)|
+\rho^{1-\theta}|\partial_{m\rho}\tilde{\eta}(\rho,\rho u)|\leq C_{\gamma},
\end{cases}
\end{equation}
It follows from integration by parts that
\begin{align}\label{4.0}
K_4&\leq \varepsilon\int^T_0\int^b_{b-1}\int^{r_2}_{r_1}\Big|\int^y_r
 \partial_m\tilde{\eta}\Big\{\big((\rho+\alpha\delta\rho^{\alpha})u_z\big)_z
   +(\rho+\alpha\delta\rho^{\alpha})(\frac{N-1}{z}u)_z\nonumber\\
&\qquad\qquad\qquad\qquad\qquad\quad\quad\,\,\,+(\alpha-1)\delta(\rho^{\alpha})_z\frac{N-1}{z}u
\Big\}\,z^{N-1}\dd z\Big|\,\dd r\dd y\dd t\nonumber\\
&\leq\varepsilon\int^T_0\int^b_{b-1}\int^{r_2}_{r_1}\int^y_r(\rho+\alpha\delta\rho^{\alpha})|u_z(z^{N-1}\partial_m\tilde{\eta})_z|\,\dd z\dd r\dd y\dd t\nonumber\\
&\quad+\varepsilon\int^T_0\int^b_{b-1}\int^{r_2}_{r_1}\int^y_r\delta\rho^{\alpha}\Big|(z^{N-1}\partial_m\tilde{\eta})_z\frac{u}{z}\Big|\,\dd z\dd r\dd y\dd t\nonumber\\
&\quad+\varepsilon\int^T_0\int^b_{b-1}\int^{r_2}_{r_1}\int^y_r(\rho+\delta\rho^{\alpha})|\partial_m\tilde{\eta}(\frac{u}{z})_z|\,z^{N-1}\dd z\dd r\dd y\dd t\nonumber\\
&\quad+C\varepsilon\int^T_0\int^{r_2}_{r_1}
\big(r|(\partial_m\tilde{\eta}(\rho+\delta\rho^{\alpha})u_r)(t,r)|+\delta|\rho^{\alpha}u\partial_m\tilde{\eta}(t,r)|
\big)\,r^{N-2}\dd r\dd t\nonumber\\
&\quad+C\varepsilon\int^T_0\int^b_{b-1}\big(y(\rho+\delta\rho^{\alpha})|\partial_m\tilde{\eta}u_y|
   +\delta|(\rho^{\alpha}u\partial_m\tilde{\eta})(t,y)|\big)\,y^{N-2} \dd y\dd t.
\end{align}
In order to estimate the terms on RHS of \eqref{4.0}, we notice that
$$
e(\rho,\rho_{\ast})\,\mathbf{1}_{B(t)}(r)\geq\frac{1}{C(\rho_{\ast})}  \qquad\mbox{for $r\in B(t)$}.
$$
Then
\begin{align}\label{4.1}
&\int^b_{r_1}\rho^{\alpha}(\rho^{\theta}-\rho^{\theta}_{\ast})^2\,r^{N-1}\dd r\nonumber\\
&\leq C(\rho_{\ast})\int^b_{r_1}\mathbf{1}_{B^c(t)}(r)\rho(\rho^{\theta}-\rho^{\theta}_{\ast})^2\,r^{N-1}\dd r
+\int^b_{r_1}\mathbf{1}_{B(t)}(r)\rho^{\alpha}(\rho^{\theta}-\rho^{\theta}_{\ast})^2(t,r)\,r^{N-1}\dd r\nonumber\\
&\leq C(\rho_{\ast})\int^b_{r_1}e(\rho,\rho_{\ast})\,r^{N-1}\dd r
+C(\rho_{\ast})\int^b_{r_1}\mathbf{1}_{B(t)}(r)\,r^{N-1}\dd r\nonumber\\
&\leq C(\rho_{\ast})\int^b_{r_1}e(\rho,\rho_{\ast})r^{N-1}\dd r
+C(\rho_{\ast})\int^b_{r_1}e(\rho,\rho_{\ast})\,r^{N-1}\dd r\nonumber\\
&\leq C(\rho_{\ast},\mathcal{E}_0).
\end{align}

Using \eqref{3.9} and \eqref{4.1},
the first to third terms on RHS of \eqref{4.0} can be bounded by
\begin{align*}\label{4.2}
&C(\frac{1}{r_1})\varepsilon\int^T_0\int^b_{b-1}\int^{r_2}_{r_1}\int^y_r\Big\{\rho u^2_z+\rho^{\gamma-2}\rho^2_z
+\delta \rho^{\alpha}u^2_z+\delta \rho^{\gamma+\alpha-3}\rho^2_z+\delta z^{-2}\rho^{\alpha}u^2\nonumber\\
&\quad\qquad\qquad\qquad\qquad\qquad\,\,\,
 +\rho u^2+e(\rho,\rho_{\ast})+\delta z^{-2}\rho^{\alpha}(\rho^{\theta}-\rho^{\theta}_{\ast})^2\Big\}\,z^{N-1}\dd z\dd r\dd y\dd t\nonumber\\
&\leq C(\rho_{\ast},\mathcal{E}_0,r_1,r_2,T).
\end{align*}
The fourth term on RHS of \eqref{4.0} is bounded by
\begin{align*}
&C(r_2,\frac{1}{r_1})\int^T_0\int^{r_2}_{r_1}
\Big\{\varepsilon \rho u^2_r+\varepsilon\delta \rho^{\alpha}u^2_r+\varepsilon\delta r^{-2}\rho^{\alpha}u^2
+\rho u^2+e(\rho,\rho_{\ast})+\rho^{\alpha}(\rho^{\theta}-\rho^{\theta}_{\ast})^2\Big\}\,r^{N-1}\dd r\dd t
\nonumber\\[1mm]
&\leq C(\rho_{\ast},\mathcal{E}_0,r_1,r_2,T).
\end{align*}
The last term of \eqref{4.0} is bounded by
\begin{eqnarray}
&&r_2\varepsilon\int^T_0\int^b_{b-1}(\rho+\delta\rho^{\alpha})|u_y|
  \big(|u|+|\rho^{\theta}-\rho^{\theta}_{\ast}|\big)\,y^{N-1}\dd y\dd t\nonumber\\
&&\leq Cr_2\int^T_0\int^b_{b-1}\big(\varepsilon (\rho+\delta\rho^{\alpha})u^2_y
 +(\rho u^2+e(\rho,\rho_{\ast}))\big)\,y^{N-1}\dd y\dd t\nonumber\\
&&\quad+C(T,\mathcal{M}_2)\int^b_{b-1}\int^T_0|u(t,y)|^2\,y^{N-1}\dd t\dd y\nonumber\\
&&\leq C(\rho_{\ast},\mathcal{E}_0,r_1,r_2,T)+C(T,\mathcal{M}_2)b^{-\frac{\vartheta}{2}}.\label{4.4}
\end{eqnarray}
Finally, for $K_5,$ it follows from \eqref{7.78} that
\begin{equation}\label{4.5}
|\Phi_r|=\frac{1}{r^{N-1}}\Big|\int^r_{\delta}\big(\rho(t,y)-d(y)\big)\,y^{N-1}\dd y\Big|
\leq C(\rho_{\ast},\mathcal{E}_0)\big(r^{-N+1}+r^{-\frac{N}{2}+1}\big).
\end{equation}
Using \eqref{4.5}, we have
\begin{align}\label{4.6}
&\Big|\int^T_0\int^b_{b-1}\int^{r_2}_{r_1}\int^y_r\partial_m\tilde{\eta}\, \rho \,\Phi_z
   z^{N-1}\dd z\dd r\dd y\dd t\Big|\nonumber\\
&\leq\int^T_0\int^b_{b-1}\int^{r_2}_{r_1}\int^y_r|z^{N-1}\Phi_z|
  \rho\big(|u|+|\rho^{\theta}-\rho^{\theta}_{\ast}|\big)\,\dd z\dd r\dd y\dd t\nonumber\\
&\leq\int^T_0\int^b_{b-1}\int^{r_2}_{r_1}\int^y_r
\big(\rho u^2+|\Phi_z|^2\rho+\rho|\rho^{\theta}-\rho^{\theta}_{\ast}|^2\big)\,z^{N-1}\dd z\dd r\dd y\dd t\nonumber\\
&\leq Cr_2\int^T_0\int^b_{r_1}\big(\rho u^2+e(\rho,\rho_{\ast})\big)\,y^{N-1}\dd y\dd t
  +Cr_2\int^T_0\int^b_{r_1}|\Phi_y|^2\rho\,y^{N-1}\dd y\dd t\nonumber\\
&\leq C(T,r_2,\mathcal{E}_0)
  +Cr_2\int^T_0\int^b_{r_1}|\Phi_y|^2\rho_{\ast}\mathbf{1}_{\{\rho(t,y)\leq2\rho_{\ast}\}}\,y^{N-1}\dd y\dd t\nonumber\\
 &\quad+Cr_2\int^T_0\int^b_{r_1}|\Phi_y|^2\rho \mathbf{1}_{A(t)}\,y^{N-1}\dd y\dd t\nonumber\\
&\leq C(r_2,T,\mathcal{E}_0)+C(\rho_{\ast},\mathcal{E}_0,r_1)\int^T_0\int^b_{r_1}\rho \mathbf{1}_{\{\rho>2\rho_{\ast}\}}\,y^{N-1}\dd y\dd t\nonumber\\
&\leq C(r_2,T,\mathcal{E}_0)+C(\rho_{\ast},\mathcal{E}_0,r_1)\int^T_0\int^b_{r_1}e(\rho,\rho_{\ast})\,y^{N-1}\dd y\dd t\nonumber\\
&\leq C(\rho_{\ast},\mathcal{E}_0,r_1,r_2,T).
\end{align}
Combining all the estimates, we obtain
\begin{equation}\label{4.7}
\begin{split}
\int^T_0\int^{r_2}_{r_1}\tilde{q}(t,r)\,r^{N-1}\dd r\dd t
\leq C(\rho_{\ast},r_1,r_2,T,\mathcal{E}_0)+C(T,\mathcal{M}_2)\,b^{-\frac{\vartheta}{2}}.
\end{split}
\end{equation}
Then \eqref{6.15-1} follows from  \eqref{6.11}, \eqref{6.59}, and \eqref{4.7}. This completes the proof.
$\hfill\Box$


\bigskip

\begin{lemma}\label{lem6.40}
The smooth solution of \eqref{6.1}--\eqref{6.3} satisfies that, for any $t\in[0,T],$
\begin{equation}\label{6.43}
\|u_r(t)\|^2_{L^2}+\int^T_0\|u_t(t)\|^2_{L^2}+\|u_{rr}(t)\|^2_{L^2}\dd t
\leq C(T,\|u_{0r}\|_{L^2},\mathcal{M}_2).
\end{equation}
\end{lemma}

\noindent{\bf Proof.}
It follows from $\eqref{6.1}_1$ that
\begin{equation}\label{6.44}
-\varepsilon((\mu+\lambda)u_r)_r+\rho u_t=\mathcal{H},
\end{equation}
where $\mathcal{H}:=-\rho uu_r-p_r+\varepsilon(\mu+\lambda)(\frac{N-1}{r}u)_r+\varepsilon\frac{N-1}{r}u\lambda_r-\rho\Phi_r.$
Multiplying \eqref{6.44} by $u_t$ and integrating it over $[\delta,b],$ we have
\begin{equation}\label{6.55}
\frac{\varepsilon}{2}\frac{\dd }{\dd t}\int^b_{\delta}(\mu+\lambda)|u_r|^2\,\dd r
+\int^b_{\delta}\rho u^2_t\,\dd r
=\frac{\varepsilon}{2}\int^b_{\delta}(\mu+\lambda)_t|u_r|^2\,\dd r
+\int^b_{\delta}\mathcal{H}u_t\,\dd r.
\end{equation}
Using \eqref{6.11}, \eqref{6.1222}, \eqref{6.16}, \eqref{6.34},
and the Sobolev inequalities, we obtain
$$
\|u_r\|_{L^{\infty}}\leq C\big(\|u_r\|_{L^2}+\|u_r\|^{\frac{1}{2}}_{L^2}\|u_{rr}\|^\frac{1}{2}_{L^2}\big),
$$
and
\begin{align}\label{6.56}
\frac{\varepsilon}{2}\int^b_{\delta}(\mu+\lambda)_t|u_r|^2\,\dd r
&\leq C(T,\mathcal{M}_2)\int^b_{\delta}\big(|\rho_ru|+|u_r|+|u|\big)|u_r|^2\,\dd r\nonumber\\
&\leq C(T,\mathcal{M}_2)\Big(\|u\|^{\frac{1}{2}}_{L^2}\|\rho_r\|_{L^2}
  \big(\|u_r\|^2_{L^2}\|u_{rr}\|^{\frac{1}{2}}_{L^2}+\|u_r\|^{\frac{5}{2}}_{L^2}\big)\nonumber\\
&\qquad\qquad\qquad\,\,+\|u\|^{\frac{1}{2}}_{L^2}\|u_r\|^{\frac{5}{2}}_{L^2}+\|u_r\|^2_{L^2}
\big(\|u_r\|_{L^2}+\|u_r\|^{\frac{1}{2}}_{L^2}\|u_{rr}\|^{\frac{1}{2}}_{L^2}\big)
\Big)\nonumber\\
&\leq C(T,\mathcal{M}_2)\Big(\big(\|u_r\|^2_{L^2}+\|u_r\|^{\frac{5}{2}}_{L^2}\big)\|u_{rr}\|^{\frac{1}{2}}_{L^2}+\|u_r\|^{3}_{L^2}+1\Big).
\end{align}
Notice that
\begin{equation}\label{6.57}
\int^b_{\delta}\rho^{-1}|\rho\Phi_r|^2\,\dd r\leq C(T,\mathcal{M}_2).
\end{equation}
Then, using \eqref{6.57}, we obtain
\begin{align}\label{6.58}
\Big|\int^b_{\delta}\mathcal{H}u_t\dd r\Big|
&\leq\frac{1}{8}\int^b_{\delta}\rho|u_t|^2\,\dd r
+C\int^b_{\delta}\rho^{-1}|\mathcal{H}|^2\,\dd r\nonumber\\
&\leq C(T,\mathcal{M}_2)\big((\|u\|^2_{L^{\infty}}+1)\|(\rho_r,u_r)\|^2_{L^2}+\|u\|^2_{L^2}\big)
+\frac{1}{8}\int^b_{\delta}\rho|u_t|^2\,\dd r\nonumber\\
&\leq \frac{1}{8}\int^b_{\delta}\rho|u_t|^2\,\dd r+ C(T,\mathcal{M}_2)\big(\|u_r\|^3_{L^2}+1\big).
\end{align}
In order to close estimate \eqref{6.58},
we use \eqref{6.1222}, \eqref{6.16}, \eqref{6.34}, and \eqref{6.44} to obtain
\begin{align}\label{6.5999}
\|u_{rr}\|^2_{L^2}
&\leq C(T,\mathcal{M}_2)
\Big(\|\sqrt{\rho}u_t\|^2_{L^2}+\|\rho_r\|^2_{L^2}\|u_r\|_{L^2}\|u_{rr}\|_{L^2}
+\|\mathcal{H}\|^2_{L^2}\Big)\nonumber\\
&\leq C(T,\mathcal{M}_2)\Big(\|\sqrt{\rho}u_t\|^2_{L^2}
+\|u_r\|_{L^2}\|u_{rr}\|_{L^2}+\|u_r\|^3_{L^2}+1\Big)\nonumber\\
&\leq C(T,\mathcal{M}_2)\Big(\|\sqrt{\rho}u_t\|^2_{L^2}+\|u_r\|^3_{L^2}+1\Big).
\end{align}
Using \eqref{6.55}--\eqref{6.5999}, we have
\begin{equation}\nonumber
\frac{\dd }{\dd t}\int^b_{\delta}(\mu+\lambda)|u_r|^2\,\dd r
+\int^b_{\delta}\rho u^2_t\,\dd r
\leq C(T,\mathcal{M}_2)\Big(1+\|u_r\|^2_{L^2}\int^b_{\delta}(\mu+\lambda)|u_r|^2\,\dd r\Big).
\end{equation}
Applying the Gr\"{o}nwall inequality, we obtain
\begin{equation}\nonumber
\int^b_{\delta}(\mu+\lambda)|u_r|^2\,\dd r+\int^t_0\int^b_{\delta}\rho u^2_t\,\dd r\dd s\leq C(T,||u_{0r}||_{L^2},\mathcal{M}_2).
\end{equation}
together with \eqref{6.5999}, yields \eqref{6.43}.
$\hfill\Box$

\section{Limits of the Approximate solutions for the CNSPEs}

In this section, we first pass to the limit, $b\rightarrow\infty$, to obtain
global strong solutions $(\rho^{\v,\d}, u^{\v,\delta})$ of CNSPEs with required uniform bounds.
Then,  taking the limit, $\delta\rightarrow0+$, we obtain global spherically symmetric
solutions of
CNSPEs \eqref{nsp}
with some desired uniform bounds on $[0,T]\times [0, \infty)$,
which are important for us to utilize the compensated compactness framework in \S 6.

\subsection{Passage the limit: $b\rightarrow \infty$}
$\,\,$ In this subsection,
we fix parameters $(\v, \d)$ and denote the solution of \eqref{6.1}--\eqref{6.3}
as $(\rho^{\v,\d,b}, u^{\v,\d,b})$.
Similar to \cite{Chen2021,Chen2020} (see Appendix for details),  we can construct
that there exist sequences of smooth approximate initial data functions $(\rho_0^{\v,\d,b}, u_0^{\v,\d,b})$
and $(\rho_0^{\v,\d}, u_0^{\v,\d})$ satisfying \eqref{upperlower} and
\begin{align}\label{5.6}
\begin{cases}
\displaystyle(\rho^{\v,\d,b}_0, m^{\v,\d,b}_0)(r)\rightarrow (\rho_0^{\v,\delta}, m_0^{\v,\d})(r) \qquad
  &\mbox{in $L^{p}_{\rm loc}([\delta, \infty))\times L^1_{\rm loc}([\delta, \infty))\, $
  as $b\rightarrow \infty$},\\[2mm]
\displaystyle (\mathcal{E}_0^{\v,\d,b},  \mathcal{E}_1^{\v,\d,b})\rightarrow (\mathcal{E}_0^{\v,\d}, \mathcal{E}_1^{\v,\d})\qquad &\mbox{as $b\rightarrow \infty$},\\[2mm]
\displaystyle \mathcal{E}_2^{\v,\d,b}+ \mathcal{E}_3^{\v,\d,b}+ \|u^{\v,\d,b}_{0r}\|_{L^2} \quad   &\mbox{is uniformly bounded with respect to $b$},
\end{cases}
\end{align}
where $p=\max\{\gamma,\frac{2N}{N+2}\},$
\begin{align}\label{1.53-1}
\begin{split}
\mathcal{E}_0^{\v,\d}:&=\int_\d^\infty\big(\bar{\eta}^{\ast}(\rho_0^{\v,\d},m_0^{\v,\d})+|\Phi^{\varepsilon,\delta}_{0r}|^2\,\big)\,r^{N-1} \dd r<\infty,\\
\mathcal{E}_1^{\v,\d}:&= \v^2\int_\d^\infty \big(1+2\alpha\d (\rho_0^{\v,\delta})^{\alpha-1}+\alpha^2\d^2 (\rho_0^{\v,\delta})^{2\alpha-2}\big)
  \big|(\sqrt{\rho_{0}^{\v,\delta}})_r\big|^2\, r^{N-1} \dd r <\infty.
\end{split}
\end{align}

\smallskip
From  \eqref{6.11}, \eqref{6.1222}, \eqref{6.16}, \eqref{6.34}, and \eqref{6.43},
there exists a constant $\tilde{C}>0$ that may depend on $(\v, \delta, T)$, but is independent of $b$, so that
\begin{align}
&0<\tilde{C}^{-1}\leq \rho^{\v,\d,b}(t,r) \leq \tilde{C}, \label{5.1-1}\\
&\sup_{t\in[0,T]}\Big(\big\|(\rho^{\v,\d,b}-\rho_{\ast},  u^{\v,\d,b})(t)\big\|^2_{H^1([\d,b])}
  +\big\|\rho^{\v,\d,b}_t(t)\big\|^2_{L^2([\d,b])}\Big)
  +\int_0^T\big\|(u^{\v,\d,b}_t, u^{\v,\d,b}_{rr})(t)\big\|^2_{L^2([\d,b])}\dd t\leq \tilde{C}.\label{5.1}
\end{align}
We extend $\rho^{\v,\d,b}(t,r)$ and $u^{\v,\d,b}(t,r)$ to $[0,T]\times[\d,\infty)$
by defining $\rho^{\v,\d,b}(t,r)=\rho_{\ast}$ and $u^{\v,\d,b}(t,r)=0$ for all $(t,r)\in [0,T]\times(b,\infty)$.
Then it follows from \eqref{5.1} and  the  Aubin-Lions lemma that
\begin{align}\nonumber
(\rho^{\v,\d,b}, u^{\v,\d,b}) \qquad \mbox{is compact in $C([0,T]; L^p_{\rm loc}[\d,\infty))$ with $ p\in[1,\infty)$}.
\end{align}

To be more precise, we obtain
\begin{lemma}\label{lem5.1}
There exist functions $(\rho^{\v,\delta},u^{\v,\delta})(t,r)$ such that,  as $b\rightarrow\infty$ $($up to a subsequence$)$,
\begin{align}\nonumber
(\rho^{\v,\d,b}, u^{\v,\d,b})\rightarrow (\rho^{\v,\delta},u^{\v,\delta})  \qquad \mbox{strongly in $C([0,T]; L^p_{\rm loc}[\d,\infty))$ for all $p\in[1,\infty)$.}
\end{align}
In particular,  as $b\to\infty$ $($up to a subsequence$)$,
\begin{align}\nonumber
(\rho^{\v,\d,b}, u^{\v,\d,b})\rightarrow (\rho^{\v,\delta},u^{\v,\delta})  \qquad \mbox{a.e.\, $(t,r)\in  [0,T]\times[\d,\infty)$}.
\end{align}
\end{lemma}

Using Lemma \ref{lem5.1}, it is straightforward to show that $(\rho^{\v,\delta},u^{\v,\delta})$ is a global weak solution of the initial-boundary value problem (IBVP)
of CNSPEs \eqref{6.1}:
\begin{align}\label{5.3}
\begin{cases}
(\rho, u)(0,r)=(\rho_0^{\v,\delta}, u_0^{\v,\delta})(r) &\qquad \mbox{for $r\in [\d, \infty)$},\\
u|_{r=\d}=0 &\qquad\mbox{for $t\ge 0$}.
\end{cases}
\end{align}
Moreover,  it yields from \eqref{5.1-1}--\eqref{5.1} and the lower semicontinuity of approximate solutions that
\begin{align}
&0<\tilde{C}^{-1}\leq \rho^{\v,\delta}(t,r) \leq \tilde{C},\label{5.5}\\
&\sup_{t\in[0,T]}\big(\|(\rho^{\v,\delta}-\rho_{\ast},  u^{\v,\delta})(t)\|^2_{H^1([\d,\infty))}
 + \|\rho_t^{\v,\d}(t)\|^2_{L^2([\d,\infty))}\big)
+\int_0^T\|(u_t^{\v,\delta}, u_{rr}^{\v,\delta})(t)\|^2_{L^2([\d,\infty))}\dd t\leq \tilde{C}.\label{5.4}
\end{align}
These facts show that the weak solution $(\rho^{\v,\delta}, u^{\v,\delta},\Phi^{\v,\delta})$ of \eqref{5.3}
is indeed a strong solution. The uniqueness of the strong solution $(\rho^{\v,\delta},u^{\v,\delta},\Phi^{\v,\delta})$
is guaranteed by properties \eqref{5.5}--\eqref{5.4} and the $L^2$--energy basic estimate.
This indicates that the whole sequence $(\rho^{\v,\d,b},u^{\v,\d,b})$
converges to $(\rho^{\v,\delta},u^{\v,\delta})$ as $b\to\infty$.

\vspace{1.5mm}
Then it is straightforward to see that
$(\rho^{\v,\d}, \mathcal{M}^{\v,\d},\Phi^{\v,\delta})(t,\textbf{x})=(\rho^{\v,\d}(t,r),m^{\v,\d}(t,r)\, \frac{\textbf{x}}{r},\Phi^{\v,\delta}(t,r))$
with $\rho^{\v,\d}(t,\textbf{x})>0$
is a strong solution of the initial-boundary problem of system \eqref{nsp} with $(\mu,\lambda)$ determined by
\eqref{3.3a}
for  $(t,\textbf{x})\in [0,\infty)\times \big(\mathbb{R}^{N}\backslash B_{\d}(\mathbf{0})\big)$
with initial-boundary data as follows:
\begin{align*}\label{5.5-2}
\begin{cases}
(\rho^{\v,\d}, \mathcal{M}^{\v,\d})(0,\textbf{x})=(\rho^{\v,\d}_0(r), m^{\v,\d}_0(r)\,\frac{\textbf{x}}{r}),\\[1.5mm]
\mathcal{M}^{\v,\d}(t,\textbf{x})|_{\textbf{x}\in \partial B_\d(\mathbf{0})}=\mathbf{0}.
\end{cases}
\end{align*}

\medskip
From  Lemma \ref{lem5.1}, \eqref{6.11}, \eqref{6.1222}, \eqref{6.99}, \eqref{6.14}, \eqref{6.15-1}--\eqref{6.16}, \eqref{5.6},
Fatou's lemma, and the lower semicontinuity, we have

\begin{proposition}\label{prop5.1}
Under assumption \eqref{5.6}, for any fixed $(\v, \delta)$,  there exists a unique strong solution $(\rho^{\v,\delta}, u^{\v,\delta},\Phi^{\v,\delta})$
of IBVP \eqref{5.3}. Moreover, $(\rho^{\v,\delta}, u^{\v,\delta})$ satisfies \eqref{5.5} and, for $t>0$,
\begin{align*}
&\int_\d^\infty  \Big(\frac12\rho^{\v,\delta} |u^{\v,\delta}|^2+e(\rho^{\v,\delta},\rho_{\ast})+\frac{1}{2}|\Phi^{\varepsilon,\delta}_r|^2 \Big)(t,r)\,r^{N-1}\dd r\nonumber\\
&
\quad+\v\int_0^T\int_\d^\infty \Big(\rho^{\v,\delta}|u_r^{\v,\delta}|^2+\rho^{\v,\delta} \frac{|u^{\v,\delta}|^2}{r^2}\Big)(s,r)\,r^{N-1}\dd r \dd s\nonumber\\
&\quad +c_N\v\delta \int_0^T\int_\d^\infty  \Big((\rho^{\v,\delta})^{\alpha}\big(|u_r^{\v,\delta}|^2+ \frac{|u^{\v,\delta}|^2}{r^{2}}\big)\Big)(s,r)\,r^{N-1}\dd r \dd s
\leq \mathcal{E}_0^{\v,\delta}\le C(\mathcal{E}_0+1),
  \\[3mm]
&\v^2 \int_\d^\infty\Big(\big|(\sqrt{\rho^{\v,\d}})_r\big|^2 +\d (\rho^{\v,\delta})^{\alpha-2}|\rho_r^{\v,\delta}|^2
  +\d^2 (\rho^{\v,\delta})^{2\alpha-3}|\rho_r^{\v,\delta}|^2\Big)(t,r)\, r^{N-1}\dd r\nonumber\\
&\quad+\v\int_0^T\int_\d^\infty \left(\big|\big((\rho^{\v,\d})^{\frac{\gamma}{2}}\big)_r\big|^2
  +\d (\rho^{\v,\delta})^{\g+\alpha-3} |\rho_r^{\v,\delta}|^2 \right)(s,r)\, r^{N-1}\dd r\dd s \le C (\mathcal{E}_0+1),
\\[2mm]
&\int_0^T\int_{r_1}^{r_2} (\rho^{\v,\delta})^{\g+1}(t,r)\,\dd r \dd t\leq C(r_1,r_2,T,\mathcal{E}_0),
\\[3mm]
&\int_0^T\int_{r_1}^{r_2}  \big(\rho^{\v,\delta}|u^{\v,\delta}|^3+(\rho^{\v,\delta})^{\gamma+\theta}\big)(t,r)\,r^{N-1}\dd r\dd t
 \leq C(r_1,r_2,T,\mathcal{E}_0),
\end{align*}
for any fixed $T>0$ and any compact subset $[r_1,r_2]$ of $[\d,\infty)$.
\end{proposition}

\subsection{Passage the limit: $\delta\rightarrow0+$}$\,\,$
In this subsection, for fixed $\v>0$, we consider the limit: $\delta\rightarrow0+$ to obtain the weak solution of CNSPEs.
It follows from Lemma \ref{lem8.4} in the appendix that
\begin{equation}\label{5.11}
\begin{cases}
(\rho_0^{\v,\delta}, m_0^{\v,\d})(r)\rightarrow (\rho_0^\v, m_0^\v)(r) \,\,\,&\mbox{in $L^{p}_{\rm loc}([0,\infty); r^{N-1}\dd r)\times L^1_{\rm loc}([0,\infty); r^{N-1}\dd r)\,$
   as $\delta\rightarrow 0+$},\\[1.5mm]
(\mathcal{E}_0^{\v,\d}, \mathcal{E}_1^{\v,\d}) \rightarrow (\mathcal{E}_0^\v, \mathcal{E}_1^\v) \qquad &\mbox{as $\d\rightarrow0+$},
\end{cases}
\end{equation}
where $p=\max\{\gamma,\frac{2N}{N+2}\}.$
To take the limit, we need to be careful, since the weak solution may contain the vacuum state.
Here, we adopt similar compactness arguments as in \cite{MV,Guo-Jiu-Xin-2} to
consider the limit process: $\delta\rightarrow0+$.
We first extend our solution $(\rho^{\v,\delta}, u^{\v,\delta})$ as the zero extension of $(\rho^{\v,\delta}, u^{\v,\delta})$
outside  $[0,T]\times[\d,\infty)$.
In the following lemma, we collect some results that are used for the limits for self-containedness.

\begin{lemma}\label{5.3a}
The extended solution sequence $(\rho^{\v,\delta}, m^{\v,\delta})$
satisfies the following{\rm :}
\begin{enumerate}
\item[\rm (i)] There exists a function $\rho^{\v}(t,r)$ such that, as $\delta\rightarrow0+$ $($up to a subsequence$)$,
\begin{align}
(\rho^{\v,\delta},\sqrt{\rho^{\v,\delta}}) \rightarrow (\rho^{\v},\sqrt{\rho^{\v}})
 \qquad \mbox{a.e. and strongly in $C(0,T; L^{q}_{\rm loc})$ for any $q\in [1,\infty)$},\label{5.12}
\end{align}
where $L^q_{\rm loc}$ means $L^q(\mathcal{K})$ for any $\mathcal{K}\Subset (0,\infty)$.

\item[\rm (ii)]
The sequence of pressure function  $p(\rho^{\v,\d})$ is uniformly bounded in $L^\infty(0,T; L^q_{\rm loc})$
for all $q\in[1,\infty]$ and, as $\delta\to 0+$ $($up to a subsequence$)$,
\begin{align}\label{5.16}
p(\rho^{\v,\d})\rightarrow p(\rho^{\v}) \qquad  \mbox{strongly in $L^q(0,T; L^q_{\rm loc})$ for all $q\in[1,\infty)$}.
\end{align}

\item[\rm (iii)]
As $\delta\rightarrow 0+$ $($up to a subsequence$)$,
$m^{\v,\delta}$
converges strongly in $L^2(0,T; L^q_{\rm loc})$
to some function $m^{\v}(t,r)$ for all $q\in[1,\infty)$, which yields
\begin{align}
m^{\v,\d}(t,r)=(\rho^{\v,\delta} u^{\v,\delta})(t,r)\rightarrow m^{\v}(t,r) \qquad  \mbox{a.e. in $[0,T]\times(0,\infty)$}. \nonumber
\end{align}
\end{enumerate}
\end{lemma}
In Lemma \ref{5.3a},  the convergence results in (i)--(iii) are from
Lemma 5.3, Corollary 5.4, and Lemma 5.5 in \cite{Chen2020}.

\medskip
The proof of the following lemma is similar to that in \cite[Lemma 4.4]{Chen2021}.
\begin{lemma}\label{lem5.6}
$m^{\v}(t,r)=0$ a.e. on $\{(t,r)\, :\, \rho^{\v}(t,r)=0\}$. Furthermore,
there exists a function $u^{\v}(t,r)$ so that $m^{\v}(t,r)=\rho^{\v}(t,r) u^{\v}(t,r)$ a.e.,
$\, u^\v(t,r)=0$ a.e. on $\{(t,r)\, :\, \rho^{\v}(t,r)=0\}$,
and
\begin{align}
& m^{\v,\d} \rightarrow m^{\v} &&\mbox{strongly in $L^2(0,T; L^p_{\rm loc})$ for $p\in[1,\infty)$},\nonumber\\
&  \frac{m^{\v,\d}}{\sqrt{\rho^{\v,\d}}}  \rightarrow \frac{m^\v}{\sqrt{\rho^\v}}=\sqrt{\rho^\v}u^\v &&\mbox{strongly in $L^2(0,T; L^2_{\rm loc})$}.\nonumber
\end{align}
\end{lemma}

\medskip
Then we can obtain the following theorem for the global weak solutions of CNSPEs:

\begin{theorem}\label{thm5.10}
Let $(\rho^\v, m^\v)$ be the limit of $(\rho^{\v,
\delta}, m^{\v,\delta}),$  as $\delta\rightarrow 0+.$ Let $(\rho_0^\v, m_0^\v)$ be the initial data satisfying \eqref{ip}--\eqref{Ebound}.
For each $\v>0$, there exists a global spherical symmetry weak solution
$$
(\rho^\v, \M^\v,\Phi^{\v})(t,\mathbf{x}):=\left(\rho^\v(t,r), m^\v(t,r)\frac{\mathbf{x}}{r},\Phi^{\v}(t,r) \right)
$$
of CNSPEs \eqref{nsp} in the sense of Definition {\rm \ref{weaknsp}}.
Moreover,
$$
(\rho^\v, m^\v,\Phi^{\v})(t,r)=(\rho^\v(t,r), \rho^\v(t,r) u^\v(t,r),\Phi^{\v}(t,r)),
$$
with
$u^\v(t,r):=\frac{m^\v(t,r)}{\rho^\v(t,r)}$ a.e. on  $\{(t,r)\,:\,\rho^\v(t,r)\ne 0\}$ and $u^\v(t,r):=0$ a.e. on $\{(t,r)\,:\, \rho^\v(t,r)=0\,\, \mbox{or $\,r=0$}\}$,
satisfies the following bounds{\rm :}
\begin{align}
&\rho^\v(t,r)\ge 0 \,\,\, a.e., \qquad \big(\frac{m^\v}{\sqrt{\rho^\v}}\big)(t,r)=\sqrt{\rho^\v(t,r)}u^\v(t,r)=0\,\,\, \mbox{a.e. on $\{(t,r)\,:\, \rho^\v(t,r)=0\}$},\\[1mm]
&\int_0^\infty
\Big(\frac12 \Big|\frac{m^{\v}}{\sqrt{\rho^\v}}\Big|^2 +e(\rho^{\v},\r_{\ast})+\frac{1}{2}|\Phi^{\v}_{r}|^2 \Big)(t,r)\,r^{N-1} \dd r\nonumber\\
&\quad+\v\int_{\mathbb{R}_+^2} \Big|\frac{m^{\v}}{\sqrt{\rho^\v}}\Big|^2(s,r)\, r^{N-3}\dd r \dd s\leq \mathcal{E}_0^\v\le \mathcal{E}_0+1\quad\,\,\mbox{for $t>0$},\label{5.55}\\
&\v^2 \int_0^\infty \big|\big(\sqrt{\rho^{\v}(t,r)}\big)_r\big|^2\,r^{N-1}\dd r
  +\v\int_{\mathbb{R}_+^2} |\big((\rho^{\v}(s,r))^{\frac{\gamma}{2}}\big)_r|^2\,r^{N-1}\dd r\dd s
\leq C (\mathcal{E}_0+1)\quad\mbox{for $t>0$},\label{5.56}\\
&\int_0^T\int_{r_1}^{r_2} (\rho^{\v})^{\g+1}(t,r)\,\dd r\dd t\leq C(r_1,r_2,T,\mathcal{E}_0),\label{5.57}\\
&\int_0^T\int_{r_1}^{r_2}   \big(\rho^\v |u^\v|^3+(\rho^{\v})^{\gamma+\theta}\big)(t,r)\,r^{N-1}\dd r\dd t
\leq C(r_1,r_2,T, \mathcal{E}_0),\label{5.58}
\end{align}
for any fixed $T>0$ and any compact subset $[r_1,r_2]\Subset (0,\infty)$.
In addition, the following energy inequality holds{\rm :}
\begin{equation*}\label{ei2}
\begin{split}
&\int_{\R^N}\Big(\frac{1}{2}\Big|\frac{\mathcal{M}^{\varepsilon}}{\sqrt{\rho^{\varepsilon}}}\Big|^2+e(\rho^{\varepsilon},\rho_{\ast})+\frac{1}{2}|\nabla_{\mathbf{x}}\Phi^{\varepsilon}|^2\Big)(t,\mathbf{x})\dd\mathbf{x}\\
&\leq\int_{\R^N}\Big(\frac{1}{2}\Big|\frac{\mathcal{M}^{\varepsilon}_0}{\sqrt{\rho^{\varepsilon}_0}}\Big|^2+e(\rho^{\varepsilon}_0,\rho_{\ast})+\frac{1}{2}|\nabla_{\mathbf{x}}\Phi^{\varepsilon}_0|^2\Big)(\mathbf{x})\dd\mathbf{x} \qquad \text{ for } t\geq0.
\end{split}
\end{equation*}
For any entropy pair $(\eta, q)$ defined in \eqref{weakentropy} for any smooth
compactly supported function $\psi(s)$ on $\mathbb{R}$,
\begin{align}
\partial_t\eta(\rho^\v,m^\v)+\partial_rq(\rho^\v,m^\v) \qquad  \mbox{is compact in $H^{-1}_{\rm loc}(\mathbb{R}^2_+)$}.\nonumber
\end{align}
\end{theorem}

\smallskip
\noindent{\bf Proof.}
We divide the proof into four steps:

\medskip
\noindent
{\bf 1.} It follows from Fatou's lemma, the lower semicontinuity,  and Proposition \ref{prop5.1}
that,
{\it under assumption \eqref{5.11}, for any fixed $\v$ and  $T>0$, the limit functions $(\rho^\v, m^\v)=(\rho^\v, \rho^\v u^\v)$
satisfy
\begin{align*}
&\rho^{\v}(t,r)\geq0 \,\,\,\, a.e.,\\
&u^\v(t,r)=0, \,\,\big(\frac{m^\v}{\sqrt{\rho^\v}}\big)(t,r)=\sqrt{\rho^\v}(t,r)u^\v(t,r)=0 \qquad\, a.e.\,\, \mbox{on $\{(t,r)\,:\, \rho^\v(t,r)=0\}$}, 
\\[1mm]
&\int_0^\infty \Big(\frac12\Big|\frac{m^{\v}}{\sqrt{\rho^\v}}\Big|^2+e(\rho^{\v},\r_{\ast})+\frac{1}{2}|\Phi^{\varepsilon}_r|^2 \Big)(t,r)\, r^{N-1} \dd r\nonumber\\
&
  \quad+\v\int_{\mathbb{R}_+^2} \Big|\frac{m^{\v}}{\sqrt{\rho^\v}}\Big|^2(s,r)\,r^{N-3}\dd r \dd s\leq \mathcal{E}_0^\v\le \mathcal{E}_0+1\quad\,\, \mbox{for $t\geq0$},
  \\
&\v^2\int_0^\infty \big|(\sqrt{\rho^{\v}(t,r)})_r\big|^2\,r^{N-1}\dd r
  +\v\int_{\mathbb{R}_+^2} \big|\big((\rho^{\v}(s,r))^{\frac{\gamma}{2}}\big)_r\big|^2\,r^{N-1}\dd r\dd s
    \leq C(\mathcal{E}_0+1)\quad\,\,\mbox{for $t\geq 0$},
    \\
&\int_0^T\int_{r_1}^{r_2} (\rho^{\v})^{\g+1}(t,r)\,\dd r\dd t\leq C(r_1,r_2,T,\mathcal{E}_0),
\\[1mm]
&\int_0^T\int_{r_1}^{r_2}  \big(\rho^{\v}|u^\v|^3+(\rho^{\v})^{\gamma+\theta}\big)(t,r)\,r^{N-1} \dd r\dd t
  \leq C(r_1,r_2,T,\mathcal{E}_0),
\end{align*}
where $[r_1,r_2]\Subset (0,\infty)$.
}

\medskip
\noindent
{\bf 2.} Similar to \cite[Lemma 4.6]{Chen2021}, we can obtain the convergence of the potential function $\Phi^{\varepsilon,\delta}$:
{\it For fixed $\varepsilon>0,$ there exists a function $\Phi^{\varepsilon}(t,\mathbf{x})=\Phi^{\varepsilon}(t,r)$ such that, as $\delta\rightarrow0+$ $($up to a subsequence$)$,
\begin{align}
&\Phi^{\varepsilon,\delta}\rightharpoonup\Phi^{\varepsilon}
\mbox{ weak-$\ast$ in } L^{\infty}(0,T;H^1_{\rm loc}(\R^N)) \mbox{ and weakly in } L^2(0,T; H^{1}_{\rm loc}(\R^N)),\label{philimit}\\
&\Phi^{\varepsilon,\delta}_r(t,r)r^{N-1}\rightarrow\Phi^{\varepsilon}_r(t,r)r^{N-1}=-\int^r_0\big(\rho^{\varepsilon}(t,z)-d(z)\big)\,z^{N-1}dz \,\,\,\,\,\mbox{ in $C_{\rm loc}([0,T]\times[0,\infty))$,}\label{philimit2}\\
&\|\Phi^{\varepsilon}(t)\|_{L^{\frac{2N}{N-2}}(\R^N)}+\|\nabla\Phi^{\varepsilon}(t)\|_{L^2(\R^N)}\leq C(\mathcal{E}_0) \qquad  \mbox{ for $t\geq0$}.\label{pestimate}
\end{align}
}

\medskip
Similar to \cite[Lemma 4.8]{Chen2021}, we can also obtain the following energy inequality:
{\it For $\gamma>1,$  the following energy inequality holds{\rm :}
\begin{equation*}\label{ei3}
\begin{split}
&\int^{\infty}_0\Big(\frac{1}{2}\Big|\frac{m^{\varepsilon}}{\sqrt{\rho^{\varepsilon}}}\Big|^2
+e(\rho^{\varepsilon},\rho_{\ast})\Big)(t,r)\,r^{N-1}\dd r
+\frac{1}{2}\int^{\infty}_0|\Phi^{\varepsilon}(t,r)|^2\,r^{N-1}\dd r\\
&\leq\int^{\infty}_0\Big(\frac{1}{2}\Big|\frac{m^{\varepsilon}_0}{\sqrt{\rho^{\varepsilon}_0}}\Big|^2+e(\rho^{\varepsilon}_0,\rho_{\ast})\Big)(r)\,r^{N-1}\dd r
+\frac{1}{2}\int^{\infty}_0|\Phi^{\varepsilon}_0(r)|^2\,r^{N-1}\dd r.\\
\end{split}
\end{equation*}
}

\medskip
\noindent
{\bf 3.}
We now show that
\begin{align*}\label{3.35-1}
(\rho^\v, \M^\v,\Phi^{\v})(t,\mathbf{x})=(\rho^\v(t,r), m^\v(t,r)\frac{\mathbf{x}}{r},\Phi^\v(t,r))
\end{align*}
is a weak solution of the Cauchy problem \eqref{nsp} and \eqref{appinitial}
in $\mathbb{R}^N$ in the sense of Definition \ref{weaknsp}, following the arguments as for Lemmas 4.9--4.10
in \cite{Chen2021}.
More precisely,
{\it the function sequence $( \rho^\v, \M^\v,\Phi^{\v})(t,\mathbf{x})$ satisfies
the following properties:
\begin{enumerate}
\smallskip
\item[\rm (i)]
Let $0\leq t_1<t_2\leq T$,
and let $\zeta(t,\mathbf{x})\in C^1([0,T]\times\mathbb{R}^N)$ be any smooth function with compact support.
Then
\begin{align*}
\int_{\mathbb{R}^N} \rho^{\v}(t_2,\mathbf{x}) \zeta(t_2,\mathbf{x})\, \dd\mathbf{x}
=\int_{\mathbb{R}^N} \rho^{\v}(t_1,\mathbf{x}) \zeta(t_1,\mathbf{x})\, \dd\mathbf{x}
 +\int_{t_1}^{t_2} \int_{\mathbb{R}^N} \big(\rho^{\v} \zeta_t + \M^{\v}\cdot\nabla\zeta\big)\,\dd\mathbf{x}\dd t.
 \end{align*}

\smallskip
\item[\rm (ii)] Let  $\psi(t,\mathbf{x})\in \left(C^2_0([0,\infty)\times \mathbb{R}^N)\right)^N$ be any smooth function with $\mbox{supp}\, \psi \Subset [0, T)\times \mathbb{R}^N$ for some fixed $T\in (0, \infty)$.
Then, it has
\begin{align*}
&\int_{\mathbb{R}_+^{N+1}} \Big\{\M^{\v}\cdot\partial_t\psi
+\frac{\M^{\v}}{\sqrt{\rho^{\v}}} \cdot \big(\frac{\M^{\v}}{\sqrt{\rho^{\v}}}\cdot \nabla\big)\psi + p(\rho^{\v})\,\mbox{\rm div}\,\psi\Big\}\, \dd\mathbf{x}\dd t\nonumber\\
&\quad-\int_{\mathbb{R}_+^{N+1}}\rho^{\v}\nabla_\mathbf{x}\Phi^{\varepsilon}\cdot\psi d\mathbf{x}dt
   +\int_{\mathbb{R}^N} \M_0^\v(\mathbf{x})\cdot \psi(0,\mathbf{x})\,\dd\mathbf{x} \nonumber\\
&=-\v\int_{\mathbb{R}_+^{N+1}}
\Big\{\frac{1}{2}\M^{\v}\cdot \big(\Delta \psi+\nabla\mbox{\rm div}\,\psi \big)+ \frac{\M^{\v}}{\sqrt{\rho^{\v}}} \cdot \big(\nabla\sqrt{\rho^{\v}}\cdot \nabla\big)
+ \nabla\sqrt{\rho^{\v}}  \cdot \big(\frac{\M^{\v}}{\sqrt{\rho^{\v}}}\cdot \nabla\big)\Big\} \psi\, \dd\mathbf{x}\dd t 
\\
&=\sqrt{\v}\int_{\mathbb{R}_+^{N+1}}
\sqrt{\rho^{\v}} \Big\{\mathrm{V}^{\v}  \frac{\mathbf{x}\otimes\mathbf{x}}{r^2}
  +\frac{\sqrt{\v}}{r}\frac{m^\v}{\sqrt{\rho^\v}}\big(I_{N\times N}-\frac{\mathbf{x}\otimes\mathbf{x}}{r^2}\big)\Big\}: \nabla\psi\, \dd\mathbf{x}\dd t,
\end{align*}
with $\mathrm{V}^{\v}(t,\mathbf{x})\in L^2(0,T; L^2(\mathbb{R}^N))$ as a function satisfying
$\displaystyle\int_0^T\int_{\mathbb{R}^N}
|\mathrm{V}^{\v}(t,\mathbf{x})|^2\,\dd\mathbf{x}\dd t\leq C\mathcal{E}_0 $ for some $C>0$, independent of $T>0$.

\smallskip
\item[\rm (iii)]
It follows from \eqref{philimit2} that $\Phi^{\v}$ satisfies the Poisson equation in the classical sense except for the origin:
$$
\Delta\Phi^{\v}=\rho^{\v}(t,\mathbf{x}), \quad\quad \text{ for } t\geq0,\quad \mathbf{x}\in \R^N\setminus\{\mathbf{0}\}.
$$
Moreover, for any smooth function $\xi(\mathbf{x})\in C^1_0(\R^N)$ with compact support, it has
 \begin{equation*}\label{pdef}
\int_{\R^N}\nabla\Phi^{\varepsilon}(t,\mathbf{x})\cdot\xi(\mathbf{x})d\mathbf{x}
=\int_{\R^N}\big(\rho(t,\mathbf{x})-d(\mathbf{x})\big)\xi(\mathbf{x})\,\dd\mathbf{x}.
\end{equation*}
\end{enumerate}
}
\bigskip

\noindent
{\bf 4.} We also need the $H_{\rm loc}^{-1}$--compactness of weak entropy pairs:
{\it
Let $(\eta, q)$ be a weak entropy pair defined in \eqref{weakentropy}
for any smooth compact supported function $\psi(s)$ on $\mathbb{R}$.
Then
\begin{align*}\label{7.1}
\partial_t\eta(\rho^\v,m^\v)+\partial_rq(\rho^\v,m^\v) \qquad   \mbox{is compact in $ H^{-1}_{\rm loc}(\mathbb{R}^2_+)$}.
\end{align*}
}

The proof is similar with \cite[Lemma 4.12]{Chen2021}. Hence, we omit the proof for brevity.
We also refer to \cite{F. Murat,L. Tartar} as an introduction to the method of compensated compactness, especially div-curl lemma.

\medskip
In summary, Theorem \ref{thm5.10} is proved by combining
Steps 1--4 above together.
$\hfill\Box$

\medskip
\section{Proof of the Main Theorems}
In this section, we give a complete proof of Main Theorem II: Theorem \ref{theorem2},
which leads to Main Theorem I: Theorem \ref{existence1}.
Since the proof is similar to that in \cite{Chen2021},
we sketch it briefly for self-containedness.
The proof is divided into four steps.

\medskip
\noindent\textbf{1.}
The uniform estimates and compactness properties obtained in Theorem \ref{thm5.10}
imply  that the weak solutions $(\rho^{\v}, m^{\v})$ of CNSPEs \eqref{1.3}
satisfy the compensated compactness framework in Chen-Perepelitsa \cite{Chen2010} for the general case $\gamma>1$; also see LeFloch-Westdickenberg \cite{LeFloch} for $\gamma\in (1,\frac{5}{3}]$.
Then the compactness theorem in \cite{Chen2010}
shows that there exists  $(\rho,m)(t,r)$ such that
\begin{equation}\label{strongdensity}
(\rho^\v,m^{\v})\rightarrow (\rho,m) \qquad  \mbox{{\it a.e.}\, $(t,r)\in \mathbb{R}^2_+$\, as $\v\rightarrow0+$ (up to a subsequence)}.
\end{equation}

By similar arguments as in the proof of \cite[Lemma 4.4]{Chen2021}, we obtain that
$m(t,r)=0$ {\it a.e.} on $\{(t,r) \,:\, \rho(t,r)=0\}$.
We can define the  limit velocity $u(t,r)$ by denoting $u(t,r):=\frac{m(t,r)}{\rho(t,r)}$ {\it a.e.} on $\{(t,r)\,:\,\rho(t,r)\neq0\}$
and $u(t,r):=0$ {\it a.e.} on $\{(t,r)\,:\, \rho(t,r)=0\, \mbox{ or $r=0$}\}$.
Then we obtain that
\begin{align}
m(t,r)=\rho(t,r) u(t,r).\nonumber
\end{align}
Similarly, we can define $(\frac{m}{\sqrt{\rho}})(t,r):=\sqrt{\rho(t,r)} u(t,r)$, which is $0$ a.e. on the vacuum states $\{(t,r) \ :\ \rho(t,r)=0\}$.
Moreover, one can obtain that, as $\v\rightarrow0+$,
\begin{equation}\label{7.0-9}
\frac{m^{\v}}{\sqrt{\rho^{\v}}}\equiv\sqrt{\rho^{\v}} u^{\v}\rightarrow \frac{m}{\sqrt{\rho}}\equiv\sqrt{\rho} u \qquad
\mbox{strongly in } \ L^2([0,T]\times[0,r_2], r^{N-1}\dd r \dd t).
\end{equation}

Notice that
$|m|^{\frac{3(\gamma+1)}{\gamma+3}}\leq C\big(\frac{|m|^3}{\rho^2}+\rho^{\gamma+1} \big)$,
which, together with \eqref{5.57}--\eqref{5.58}, yields that
\begin{align}\label{7.51}
(\rho^\v,m^{\v})\rightarrow (\rho,m) \qquad \mbox{in $L^p_{\rm loc}(\mathbb{R}^2_+)\times L^q_{\rm loc}(\mathbb{R}^2_+)$ as $\v\rightarrow 0+$}
\end{align}
for $p\in[1,\gamma+1)$ and $q\in[1,\frac{3(\gamma+1)}{\gamma+3})$,
where $L^q_{\rm loc}(\mathbb{R}^2_+)$ represents $L^q([0,T]\times \mathcal{K})$ for any $T>0$ and $\mathcal{K}\Subset (0,\infty)$.

From \eqref{7.0-9}--\eqref{7.51}, we also obtain the convergence of the relative mechanical energy as $\v\rightarrow0+$:
\begin{align}\nonumber
\bar{\eta}^{\ast}(\rho^\v,m^\v)\rightarrow \bar{\eta}^{\ast}(\rho,m)\qquad \mbox{in $L^1_{\rm loc}(\mathbb{R}_+^2)$}.
\end{align}
By passing the limit in \eqref{5.55} and noticing that $\bar{\eta}^{\ast}(\rho,m)$ is a convex function,
we have
\begin{align}\nonumber
\int_{t_1}^{t_2} \int_0^\infty(\bar{\eta}^{\ast}(\rho,m)+\frac{1}{2}|\Phi_r|^2)(t,r)\,  r^{N-1}\dd r\dd t
 \leq (t_2-t_1) \int_0^\infty(\bar{\eta}^{\ast}(\rho_0,m_0)+\frac{1}{2}|\Phi_{0r}|^2)(r)\,  r^{N-1}\dd r,
\end{align}
which yields that
\begin{align}\label{7.83}
\int_0^\infty(\bar{\eta}^{\ast}(\rho,m)+\frac{1}{2}|\Phi_r|^2)(t,r)\, r^{N-1}\dd r\leq  \int_0^\infty(\bar{\eta}^{\ast}(\rho_0,m_0)+\frac{1}{2}|\Phi_{0r}|^2)(r)\, r^{N-1}\dd r
\qquad\, \mbox{for {\it a.e.} $t\geq0$}.
\end{align}
This indicates that no concentration (Dirac mass) is formed in the density $\rho$ at the origin $r=0$.
\medskip


\noindent\textbf{2.}
For the convergence of the electric potential function $\Phi^{\varepsilon}(t,r),$ by similar calculation in \cite{Chen2021}, as $\varepsilon\rightarrow0+$ (up to a subsequence), we obtain
\begin{equation}\label{limitp}
\Phi^{\varepsilon}_r(t,r)r^{N-1}
=-\int^r_0\big(\rho^{\varepsilon}(t,y)-d(y)\big)\,y^{N-1}\dd y
\rightarrow-\int^r_0\big(\rho(t,y)-d(y)\big)\,y^{N-1}\dd y
\,\,\,\text{ a.e. $(t,r)\in \R^2_+$}.
\end{equation}
Then \eqref{pestimate} implies that there exists a function $\Phi(t,\mathbf{x})=\Phi(t,r)$ such that,
as $\varepsilon\rightarrow0+$ (up to a subsequence),
\begin{align*}
&\Phi^{\varepsilon}\rightharpoonup \Phi \text{ weak-$\ast$ in } L^{\infty}(0,T;H^1_{\rm loc}(\R^N)) \text{ and weakly in } L^2(0,T;H^1_{\rm loc}(\R^N)),
\\[1mm]
&\|\Phi(t)\|_{L^{\frac{2N}{N-2}}(\R^N)}+\|\nabla\Phi(t)\|_{L^2(\R^N)}\leq C(\mathcal{E}_0) \qquad  a.e.\,\,  t\geq0. 
\end{align*}
It follows from \eqref{limitp} and the uniqueness of the limit that
\begin{equation*}\label{unique}
\Phi_r(t,r)r^{N-1}=-\int^r_0\big(\rho(t,y)-d(y)\big)\,y^{N-1}\dd y \qquad  a.e. \,\, (t,r)\in\R^2_+.
\end{equation*}

\smallskip
\noindent
{\bf 3.}
Define the functions:
\begin{align*}\label{3.113-1}
(\rho,\M,\Phi)(t,\mathbf{x}):=(\rho(t,r), m(t,r) \frac{\mathbf{x}}{r},\Phi(t,r))=(\rho(t,r), \rho(t,r) u(t,r)\frac{\mathbf{x}}{r},\Phi(t,r)).
\end{align*}
Then we obtain that, for {\it a.e.} $t\geq0,$
\begin{equation}\nonumber
\begin{split}
&\int_{\R^N}\Big(\frac{1}{2}\Big|\frac{\mathcal{M}}{\sqrt{\rho}}\Big|^2+e(\rho,\rho_{\ast})+\frac{1}{2}|\nabla_{\mathbf{x}}\Phi|^2\Big)(t,\mathbf{x})\,\dd \mathbf{x}\\
&\leq\int_{\R^N}\Big(\frac{1}{2}\Big|\frac{\mathcal{M}_0}{\sqrt{\rho_0}}\Big|^2+e(\rho_0,\rho_{\ast})+\frac{1}{2}|\nabla_{\mathbf{x}}\Phi_0|^2\Big)(\mathbf{x})\,\dd \mathbf{x},
\end{split}
\end{equation}
This implies \eqref{finiteenergy} in Definition \ref{weakep}.
From \eqref{7.83}, we know that $\frac{\M}{\sqrt{\rho}}=\sqrt{\rho}u\,\frac{\mathbf{x}}{r}$ is well-defined and in $L^2$ for {\it a.e.} $t>0$.

\smallskip
\noindent
{\bf 4.}
Using Theorem \ref{thm5.10} and applying the similar argument as in \cite[Section 5]{Chen2021},
we can finally conclude that $(\rho,\M,\Phi)(t,\mathbf{x})$ is a global weak solution of the Cauchy problem
of CEPEs \eqref{1.1} in $\mathbb{R}^N$ and \eqref{1.2}--\eqref{initial2}
in the sense of Definition \ref{weakep}. The proof of Theorem \ref{existence1} is completed.
$\hfill\Box$

\smallskip
\appendix
\section{Construction of Approximate Initial Data}

For completeness, in this appendix, we present the construction of the approximate
initial data functions with desired estimates and regularity.
From \eqref{initial}, there exists a positive constant $\hat{r}\gg1$ so that
\begin{align*}
0<\frac12 \rho_{\ast}\leq \rho_0(r)\leq \frac32 \rho_{\ast}\qquad\,\, \mbox{for $r\geq \hat{r}$}.
\end{align*}
The density function $\rho_0(r)$ is cut off as follows:
\begin{equation*}
\tilde{\rho}_0^\v(r)=
\begin{cases}
(\beta \v)^{\frac14} \qquad &\mbox{if $\rho_0(r)\leq (\beta\v)^{\frac14} $},\\
\rho_0(r)  &\mbox{if $(\beta \v)^{\frac14} \leq \rho_0(r)\leq (\beta \v)^{-\frac12} $},\\
(\beta \v)^{-\frac12}  &\mbox{if $\rho_0(r)\geq (\beta \v)^{-\frac12} $},
\end{cases}
\end{equation*}
where $\v\in(0,1]$, and $0<\beta\ll1$ is a fixed small positive constant, which
is to guarantee $(\beta \v)^{\frac14} \ll (\beta \v)^{-\frac12} $ for all $\v\in(0,1]$. It is direct to check that
\begin{equation}\label{8.3-1}
\tilde{\rho}_0^\v(r)\leq \rho_0(r)+1,
\qquad\tilde{\rho}_0^\v(r)\rightarrow \rho_0(r) \,\,\,\mbox{as $\v\to 0+\,\,$ {\it a.e.} $r\in \mathbb{R}_+$}.
\end{equation}

To keep the $L^p$--properties of mollification, it is more convenient to smooth out the initial data in the original coordinate $\mathbb{R}^N$;
therefore we do not distinguish functions $(\rho_0,m_0)(r)$ from $(\rho_0,m_0)(\mathbf{x})=(\rho_0,m_0)(|\mathbf{x}|)$ for simplicity.

\medskip
It yields from \eqref{6.10}--\eqref{finiteenergy} that $\rho_0(\mathbf{x})\in L^\gamma_{\rm loc}(\mathbb{R}^N)$. We also assume that $\rho_0(\mathbf{x})-d(\mathbf{x})\in (L^{\frac{2N}{N+2}}\cap L^1)(\mathbb{R}^N).$
Denote $$\Delta\Phi_0=-(\rho_0-d(\mathbf{x})),$$
which is well-defined under the above assumption.
Using the convexity of $e(\rho,\rho_{\ast})$, we have
\begin{align}\label{8.4}
0\leq e(\tilde{\rho}_0^\v(\mathbf{x}),\rho_{\ast})\leq e(\rho_0(\mathbf{x}),\rho_{\ast}).
\end{align}

Using \eqref{finiteenergy}, \eqref{8.3-1}--\eqref{8.4} and the Lebesgue dominated convergence theorem, we obtain
\begin{align*}
\lim_{\v\rightarrow0+}  \int_{\mathcal{K}}\Big(|\tilde{\rho}_0^\v(\mathbf{x})-\rho_0(\mathbf{x})|^{\gamma}
+|\sqrt{\tilde{\rho}_0^\v(\mathbf{x})}-\sqrt{\rho_0(\mathbf{x})}|^{2\gamma}\Big)\,\dd\mathbf{x}=0
\qquad \mbox{for any $\mathcal{K}\Subset\mathbb{R}^N$}.\label{8.4-1}
\end{align*}

Since a better decay property for approximate initial data is needed,
we take the following cut-off density function  $\tilde{\rho}_0^\v(\mathbf{x})$ at the far-field:
\begin{equation*}
\hat{\rho}_0^\v(\mathbf{x})=
\begin{cases}
\tilde{\rho}_0^\v(\mathbf{x}) \quad &\mbox{if $|\mathbf{x}|\leq (\beta \v)^{-\frac1{2N}} $},\\[1mm]
\rho_{\ast} \quad &\mbox{if $|\mathbf{x}|> (\beta \v)^{-\frac1{2N}} $},
\end{cases}
\end{equation*}
where we further take $\beta$ small enough such that $|\mathbf{x}|\geq (\beta \v)^{-\frac1{2N}} \geq \hat{r}+2$ for all $\v\in(0,1]$.
It is obvious that  $\hat{\rho}_0^{\v}(\mathbf{x})$ is not a smooth function such that
we need to  mollify $\hat{\rho}_0^{\v}(\mathbf{x})$.
Let  $J(\mathbf{x})$ be the standard mollifier and
$J_{\sigma}(\textbf{x}):=\frac{1}{\sigma^N}J(\frac{\mathbf{x}}{\sigma})$ for $\sigma\in(0,1)$.
Here, we take $\sigma=\v^{\frac14}$ and define  $\rho_{0}^{\v}(\mathbf{x})$  as
\begin{align*}\label{8.6}
\rho_{0}^{\v}(\mathbf{x})
:=\Big(\int_{\mathbb{R}^N} \sqrt{\hat{\rho}_{0}^\v(\mathbf{x-y})} J_\sigma(\mathbf{y}) \dd\mathbf{y}\Big)^2.
\end{align*}
Then $\rho_{0}^{\v}(\mathbf{x})$ is still a spherically symmetric function,
{\it i.e.}, $\rho_{0}^{\v}(\mathbf{x})=\rho_{0}^{\v}(|\mathbf{x}|)$. By similar arguments as in \cite{Chen2021,Chen2020}, we can obtain the following Lemma (the details are omitted here for simplicity of presentation).

\begin{lemma}\label{lem8.4}
The following three statements hold{\rm :}
\begin{enumerate}
\item[\rm (i)]  	As $\v\rightarrow0+$,
\begin{align*}
&(\mathcal{E}^\v_0, \mathcal{E}_1^\v)\rightarrow (\mathcal{E}_0,0), \qquad\\
&(\rho_0^{\v}, m_0^{\v})(r)\rightarrow (\rho_0, m_0)(r) \,\,\,\, \mbox{in $L_{\rm loc}^p([0,\infty); r^{N-1}\dd r)\times L_{\rm loc}^1([0,\infty); r^{N-1}\dd r)$}.
\end{align*}
where $\mathcal{E}^\v_0, \mathcal{E}_1^\v$, and $\mathcal{E}_0$ are defined in \eqref{E0}, \eqref{E1}, and \eqref{ife}, separately, and $p\in\{\gamma,\frac{2N}{N+2}\}$.

\smallskip
\item[\rm (ii)]  For any fixed $\v\in(0,1]$, as $\d\rightarrow0+$,
\begin{align*}
&(\mathcal{E}_0^{\v,\d},\mathcal{E}_1^{\v,\d})\rightarrow (\mathcal{E}_0^\v, \mathcal{E}_1^\v), \quad\,\,\\
&
(\rho_0^{\v,\d}, m_0^{\v,\d})(r)\rightarrow (\rho_0^\v, m^\v_0)(r) \,\,\,\, \mbox{in $L_{\rm loc}^p([0,\infty); r^{N-1}\dd r)\times L_{\rm loc}^1([0,\infty); r^{N-1}\dd r)$},
\end{align*}
where $\mathcal{E}_0^{\v,\d}$ and $\mathcal{E}_1^{\v,\d}$ are defined in \eqref{1.53-1}, and
$p\in\{\gamma,\frac{2N}{N+2}\}$.

\smallskip
\item[\rm (iii)] For any fixed $(\v,\d)$, as $b\rightarrow\infty$,
\begin{align*}
&(\mathcal{E}_0^{\v,\d,b}, \mathcal{E}_1^{\v,\d,b})\rightarrow (\mathcal{E}_0^{\v,\d}, \mathcal{E}_1^{\v,\d}),
\\
&(\rho_0^{\v,\d,b}, m_0^{\v,\d,b})(r)\rightarrow (\rho_0^{\v,\d}, m_0^{\v,\d})(r) \,
 \mbox{in $L_{\rm loc}^p((\delta,\infty); r^{N-1}\dd r)\times L_{\rm loc}^1((\delta,\infty); r^{N-1}\dd r)$},
\end{align*}
where $\mathcal{E}_0^{\v,\d,b}, \mathcal{E}_1^{\v,\d,b}, \mathcal{E}_2^{\v,\d,b}$, and $\mathcal{E}_3^{\v,\d,b}$
are defined
in {\rm Lemmas \ref{bee}}--{\rm \ref{BD}}
and \eqref{4.1d}, and $p\in\{\gamma,\frac{2N}{N+2}\}$.
In addition, the upper bounds of $\mathcal{E}_0^{\v,\d,b}, \mathcal{E}_1^{\v,\d,b}, \mathcal{E}_2^{\v,\d,b}$, and
$\mathcal{E}_3^{\v,\d,b}$
are independent of $b$ {\rm (}but may depend on $\v, \d${\rm )}, and
\begin{align*}
&\qquad \mathcal{E}_0^{\v,\d,b}+ \mathcal{E}_1^{\v,\d,b}\le C(\mathcal{E}_0+1), 
\\
&\qquad \mathcal{E}_3^{\v,\d,b}:=\int_\delta^b\bar{\eta}^{\ast}(\rho_0^{\v,\delta,b},m_0^{\v,\d,b}) (1+r)^{N-1+\vartheta}\,r^{N-1} \dd r
\leq C \mathcal{E}_0 \big(\delta^{-N+1-\vartheta}+\v^{-\frac{N-1+\vartheta}{2N}}\big),
\end{align*}
for some $C>0$ independent of $(\v,\d,b)$, where $\vartheta\in(0,1)$ is arbitrary fixed constant.
\end{enumerate}
\end{lemma}

\bigskip
\noindent{\bf Acknowledgments.}
The research of Gui-Qiang G. Chen was supported in part by the UK
Engineering and Physical Sciences Research Council Awards
EP/L015811/1, EP/V008854, and EP/V051121/1. The research of Lin He was supported in part by the National Natural Sciences Foundation of China No.12001388 and 12371223, the Sichuan Youth Science and Technology Foundation No. 2021JDTD0024,
National Key R\&D Program of China No.2022YFA1007700, and
Natural Science Foundation of Sichuan Province, Grant No. 2023NSFSC1367.
The research of Yong Wang was supported in part
by the National Natural Sciences Foundation of China Grants No. 12022114 and  12288201,
Youth Innovation Promotion Association of Chinese Academy of Sciences No. 2019002,
and CAS Project for Young Scientists in Basic Research Grant No. YSBR-031.
The research of Difan Yuan was supported in part by the National Natural Sciences Foundation of China Grant No. 12001045, the UK
Engineering and Physical Sciences Research Council Award EP/V051121/1.

\medskip

\noindent
{\bf Declarations}

\noindent
{\bf Conflicts of interest/Competing interests}: The authors declare that they have no conflict

$\qquad\qquad\qquad\qquad\qquad\qquad\qquad\qquad\qquad\qquad\,\,$ of interest/Competing interest.

\newpage

\bigskip
\smallskip
%

\end{document}